\documentclass[12pt]{article}
\usepackage{graphics}
\usepackage{amsmath}
\usepackage{amssymb}
\usepackage{graphicx}
\usepackage{makeidx}
\usepackage{mathrsfs}
\usepackage{amsfonts}
\usepackage[mathscr]{eucal}
\usepackage{amsmath}

\usepackage{graphics}

\DeclareGraphicsExtensions{.eps,.pdf,.jpg,.png}

\def\e{{\rm e}}
\def\e{\hbox{e}}
\def\d{\hbox {d}}
\def\dt{\hbox{dt}}
\def\ds{\displaystyle}
\def\RR{\vbox {\hbox to 8.9pt {I\hskip-2.1pt R\hfil}}}
\def\NN{{\rm I\hskip-2pt N}}
\def\CC{{\rm C\hskip-4.8pt \vrule height 6pt width 12000sp\hskip 5pt}}
\def\q{\quad}  \def\qq{\qquad}
\def\cen{\centerline}
\def \rec#1{{\frac{1}{#1}}}
\def\pni{\par\noindent}
\def\vsh{\smallskip}
\def\vs{\medskip}
\def\vvs{\bigskip}
\def\vvvs{\bigskip\medskip} %% {\vskip 1.5truecm}
\def\vsp{\vsh\pni} %% ie. \smallskip + \par
\def\vsn{\vsh\pni}

\DeclareGraphicsExtensions{.eps,.pdf,.jpg,.png}

\begin{document}

%\cen{{\bf FRACALMO PRE-PRINT: \ http://www.fracalmo.org}}
%% \cen{{\bf Pre-print submitted }}
%
%%\cen{{\bf Special issue:   Perspectives on Fractional Dynamics and Control}}
%%
%% \cen{{\bf Guest Editors: Changpin LI  and Francesco MAINARDI }}
%\vsh
%\hrule
% \end{center}
%%%%%%%%%%%%%%%%%%%%%%%%%%%%%%%%%%%%%%%%%%%%%%%%%%%%%%%%%%%%%%%%%%%%%%%%
% \vskip 0.5truecm
\font\title=cmbx12 scaled\magstep2
\font\bfs=cmbx12 scaled\magstep1
\font\little=cmr10
\begin{center}
{\title Parametric subordination }
\\ [0.25truecm]
{\title in fractional diffusion processes\footnote{%%
E-print based on an updated version of Chapter 10, pp 227-261 in  S.C. Lim, J. Klafter and  R. Metzler (Editors):  {\it Fractional Dynamics, Recent Advances},  
World Scientific, Singapore 2012 (ISBN 9814340588, 9789814340588),}}
 \\  [0.25truecm]
 Rudolf GORENFLO$^{(1)}$  and Francesco MAINARDI$^{(2)}$ 
\\ [0.25truecm]
$\null^{(1)}$ {\little Department of Mathematics and Informatics, Free University of Berlin,}
 \\ {\little  Arnimallee  3, D-14195 Berlin, Germany}
 \\ {\little E-mail: gorenflo@mi.fu-berlin.de}
\\  [0.25truecm]
$\null^{(2)}$ {\little Department of Physics and Astronomy, University of Bologna, and INFN}
\\ {\little Via Irnerio 46, I-40126 Bologna, Italy}
\\ {\little Corresponding Author.   E-mail: francesco.mainardi@bo.infn.it}
\end{center}

\noindent 
{\it MSC}: 26A33, 33E12, 33C60, 44A10,
45K05, 60G18, 60G50, 60G52, 60K05,  76R50.
\\
{\it Key Words and Phrases}: Fractional derivatives and  integrals,
fractional diffusion, Mittag-Leffler function, Wright function, 
random walks, subordination,  self- similar stochastic processes, stable distributions, infinite divisibility.
%%%%%%%%%%%%
\begin{abstract}
\noindent
We consider simulation of spatially one-dimensional space-time fractional diffusion. 
Whereas in an earlier paper of ours (Eur.  Phys. J. Special Topics, Vol. 193,  119--132 (2011);
E-print: {\tt  http://arxiv.org/abs/1104.4041}),
 we have developed the basic theory of what we call parametric subordination via three-fold splitting applied to continuous time random walk with subsequent passage to the diffusion limit, here we go the opposite way. 
 Via Fourier-Laplace manipulations of the relevant fractional partial differential equation of evolution we obtain the subordination integral formula that teaches us how a particle path can be constructed by first generating the operational time from the physical time and then generating in operational time the spatial path. By inverting the generation of operational time from physical time we arrive at the method of parametric subordination. Due to the infinite divisibility of the stable subordinator we can simulate particle paths by discretization where the generated points of a path are precise snapshots of a true path. By refining the discretization more and more fine details of a path become visible.
\end{abstract}
\def\indice{\leaders\hbox to 1 em {\hss.\hss}\hfill}
\def\hb#1{\bf {\hbox to 1.25 truecm{ \hss#1}}}
%%\vskip 1 truecm
%%\centerline{\bf Contents}
\newpage 
\tableofcontents
 {{\bf Acknowledgements} \indice \hb{34}}

 {{\bf References} \indice \hb{34}}

%%%%% FONTS

% ==================================================

\def\pni{\par \noindent}
\def\vsh{\vskip 0.25truecm\noindent}
\def\vs{\vskip 0.5truecm}
\def\vvs{\vskip 1.0truecm}
\def\vvvs{\vskip 1.5truecm}
\def\vsp{\vsh\pni}
\def\vsn{\vsh\pni}
\def\cen{\centerline}
\def\ra{\item{(a)\ }} \def\rb{\item{(b)\ }}   \def\rc{\item{(c)\ }}
\def\eg{{e.g.}\ } \def\ie{{i.e.}\ }
\def\sg{\hbox{sign}\,}
\def\sgn{\hbox{sign}\,}
\def\sign{\hbox{sign}\,}
\def\e{\hbox{e}}
\def\exp{\hbox{exp}}
%%%% MATHEMATICS
\def\ds{\displaystyle}
\def\dis{\displaystyle}
\def\q{\quad}    \def\qq{\qquad}
\def\lan{\langle}\def\ran{\rangle}
\def\l{\left} \def\r{\right}
\def\lt{\left} \def\rt{\right}
\def\lra{\Longleftrightarrow}
\def\arg{{\rm arg}}
\def\argz{{\rm arg}\, z}
\def\argG{{x^2/ (4\,a\, t)}}
\def\d{\partial}
 \def\dr{\partial r}  \def\dt{\partial t}
\def\dx{\partial x}   \def\dy{\partial y}  \def\dz{\partial z}
\def\rec#1{\frac{1}{#1}}
\def\log{{\rm log}\,}
\def\erf{{\rm erf}\,}     \def\erfc{{\rm erfc}\,}

%%%%%%%% SETS of NATURAL, REAL, COMPLEX NUMBERS : \NN, \RR, \CC
\def\NN{{\rm I\hskip-2pt N}}
\def\MM{{\rm I\hskip-2pt M}}
\def\RR{\vbox {\hbox to 8.9pt {I\hskip-2.1pt R\hfil}}\;}
\def\CC{{\rm C\hskip-4.8pt \vrule height 6pt width 12000sp\hskip 5pt}}
%%% IDENTITY
\def\II{{\rm I\hskip-2pt I}}
%%%%%%%%%%%%%%%%%%%%%%%%%%
%% DEFINITIONS of ERROR and EXPONENTIAL FUNCTIONS %%%%%
\def\erf{{\rm erf}\,}   \def\erfc{{\rm erfc}\,}
\def\exp{{\rm exp}\,} \def\e{{\rm e}}
\def\ss{{s}^{1/2}}   %% for LAPLACE TRANSFORMS
%%%%%%%%%%%%%%%%%%%%%%%%%%%%%%%%%%%%%%%%%%%%%%%%%%%%%%%
%%%%%%%%%%%%%%%%%%%%%%%%%%
\def\N{\bar N}  %%%%%%%%%%%%%
\def\ss{{s}^{1/2}} %%%%%%%
%%%%%%%%%%
\def\Re{\hbox{Re}}
\def\Im{\hbox{Im}}
%%%%%%%%%%%%%%%%%%%%%
%\def\sss{{s}^{1/2}}   %% for LAPLACE TRANSFORMS
\def\stt{{\sqrt t}}
\def\lst{{\lambda \,\stt}}
\def\Et{{E_{1/2}(\lst)}}
\def\u{\widetilde{u}}
\def\ul{\widetilde{u}} %%% Laplace Transform  (LT)
\def\uf{\widehat{u}} %%% Fourier Transform  (FT)
\def\A{{\mathcal {A}}}
%%  justaposition symbols for Fourier, Laplace, Mellin transforms
\def\L{{\cal L}} %%% Laplace Transform !!!!
\def\F{{\cal F}} %%% Fourier Transform !!!!
\def\M{{\cal M}}  %%% Mellin Transform
\def\Fdiv{\,\stackrel{{\cal F}} {\leftrightarrow}\,}
  \def\Ldiv{\,\stackrel{{\cal L}} {\leftrightarrow}\,}
  \def\Mdiv{\,\stackrel{{\cal M}} {\leftrightarrow}\,}
\def\barr{\widetilde}
%% GREEN FUNCTIONS
\def\G{{\cal {G}}}
\def\Gc{{\cal {G}}_c}	\def\Gcs{\barr{\Gc}} %% CAUCHY PROBLEM
\def\Gs{{\cal {G}}_s}	\def\Gss{\barr{\Gs}} %% SIGNALLING PROBLEM
%%%%%% LAPLACE TRANSFORM
\def\f{\bar{f}}
\def\g{\bar{g}}
\def\u{\bar{u}}

\newcommand{\ka}{I\!\!K}
\newcommand{\rgr}{{\rm grad}}
\newcommand{\ce}{I\!\!\!\!C}
\newcommand{\re}{I\!\!R}

\newpage
\section{ Introduction}
% We summarize our earlier and recent contributions on random walks
% related to fractional diffusion processes including parametric subordination and "particle tracking".
%%%%%%%
The purpose of this chapter is to describe our method of parametric subordination
to produce particle trajectories for the so-called fractional diffusion processes.
\vsp
  By replacing in the common diffusion equation the first order time derivative
  and the second order space derivative by appropriate fractional derivatives
  we obtain a fractional diffusion equation whose solution describes the temporal evolution
  of the density of an extensive quantity, e.g. of the sojourn probability of a diffusing particle.
  \vsp
  After giving a survey on analytic methods for determination of the solution
  (this is the macroscopic aspect) we turn attention to the problem of simulation of particle
  trajectories (the microscopic aspect). By some authors such simulation is called
  "particle tracking", see e.g.\cite{Zhang-Meerschaert-Baeumer_PRE08}.
  \vsp
  As an approximate method among physicists the so-called  Continuous Time Random Walk ($CTRW$)
is very popular. On the other hand, it is possible to produce a sequence of precise snapshots
of a true trajectory.
This is achieved by a change from the "physical time" to an "operational time" in which the
simulation is carried out.
By two Markov processes happening in operational time the running of physical time and
the motion in space are produced.
Then, elimination of the operational time yields a picture of the desired trajectory.
It is remarkable that so by combination of two Markov processes a non-Markovian process is generated.
\vsp
The two Markov processes can be obtained and analyzed in two ways:
\pni (a) from the $CTRW$ model by a well-scaled passage to the "diffusion limit",
\pni (b) directly from an integral representation of the fundamental solution
     of the fractional diffusion equation.
\vsp
%% If we decide to leave away the CTRW splitting technique to arrive at parametric
%% subordination we can add to the end of the Introduction (after mentioning the two ways
%% (a) and (b)) the sentence:
%% ################################################
We have developed way (a) in our 2007 paper\cite{GorMaiViv_CSF07}  via passage to the diffusion
limit in the Cox-Weiss solution formula for $CTRW$ and by the technique of splitting the
$CTRW$ into three separate walks and passing in each of these to the diffusion limit
in our more recent recent paper\cite{Gorenflo-Mainardi_EPJ-ST11}.	 
For another access (more oriented towards measure-theoretic theory of stochastic
processes) see the recent papers \cite{Hahn-Umarov_FCAA2011, Hahn-Umarov_JTB2012}.
In \cite{Hahn-Umarov_JTB2012} the authors also treat the problem
of subordination for diffusion with distributed orders of time-fractional differentiation.
	 \vsp
	 The plan of our chapter is as follows. 
%	\vsp 
In Section 2 we provide for the reader's convenience
some preliminary notions and notations as a mathematical background for our further analysis.
% \vsp
In Section 3  we introduce the  {\it space-time fractional diffusion equation}, 
based on the {\it Riesz-Feller} and {\it Caputo} fractional derivatives, and we present the fundamental solution..
% \vsp
In Section 4 we provide the stochastic interpretation of the space-time  fractional diffusion equation 
discussing the concepts of subordination, the main goal of this chapter.
% \vsp
Finally, in Section 5 we show some graphical representations along with  conclusions. 
 % \body
\section{Notions and Notations}
In this Section we survey some preliminary notions including  
 Fourier and Laplace transforms, special functions of Mittag-Leffler and Wright type 
and L\'evy stable probability distributions.
%% that provide a necessary mathematical background for our   analysis.
%% that is  based on the inversion of the fractional integrals,
\vsp
Since in what follows we shall meet only real or
complex-valued	functions of a real variable that
are defined and  continuous in a given open interval
$I= (a,b)\,,$	$\,-\infty \le a < b \le +\infty\,,$
except, possibly, at isolated points where these
functions can be infinite,
we restrict our presentation
of the integral transforms to the class of functions
for which the Riemann improper integral on $I$
absolutely converges.
In so doing we follow Marichev \cite{Marichev_83} and
we denote this class by $L^c(I)$
or $L^c(a,b)\,.$
%%%%%%%%%%%%
\subsection{The Fourier transform}
% \vskip 0.25truecm \noindent
Let
$$  \widehat f(\kappa)	=
{\cal F} \l\{ f(x);\kappa \r\}
  = \int_{-\infty}^{+\infty} \e^{\,\ds +i\kappa x}\,f(x)\, dx\,,
  \q \kappa \in \re\,, \eqno(2.1a)$$
be the	Fourier  transform  of a %%% sufficiently well-behaved
function  $f(x) \in L^c (\RR)$, and let
$$ f(x) =
{\cal F}^{-1} \l\{ \widehat f(\kappa );x \r\}
  = \rec{2 \pi}\, \int_{-\infty}^{+\infty} \e^{\,\ds -i\kappa x}\,
 \widehat f(\kappa )\, d \kappa\,,	\q x \in \re\,, \eqno(2.1b)$$
be the inverse Fourier transform\footnote{%%
%%%%%%%%%%%%%%%%%%%%%% FOOTNOTE (1) on the FOURIER TRANSFORM %%%
%% For our purposes $f(x)$ is assumed to be continuous everywhere on
%% the real axis (except, possibly, for a finite number of points),
%% and for which the finite improper integral
%% $\int_{-\infty}^{+\infty} |f(x)|\,dx\,$ exists.
If $f(x)$ is piecewise differentiable, then the formula (2.1b)
holds true at all points where $f(x)$ is continuous and
the integral in it must be understood in the sense of the Cauchy
principal value.}.
%%%%%%% THE END OF THE FOOTNOTE (1) %%%%
\vsp
Related to the Fourier transform is the notion of pseudo-differential operator.
%%%%%%% FOOTNOTE (2) on PSEUDODIFFERENTIAL OPERATORS
Let us recall that
a generic pseudo-differential operator $A$,
acting with respect to the variable $x \in \RR\,,$
%% (of which the $n$-th derivative operator is a special case)
is defined through its Fourier representation, namely
  $$    \int_{-\infty}^{+\infty}\!\!
  \e ^{\, i\kappa x} \,  A \,[ f(x)] \, dx =
 \widehat A(\kappa )\, \widehat f (\kappa )\,,	\eqno(2.2)$$
  where
%% $f(x)$ denotes a sufficiently well-behaved function in $\RR\,,$  and
$\widehat A(\kappa)\,$ is referred to as  symbol of $A\,,$
  formally  given as
 ${\ds \widehat A (\kappa ) = \left( A\, \e^{\, -i\kappa x}\right)\,
  \e^{\, +i\kappa x}\,.} $
%%%%%%%%%%% THE END of the FOOTNOTE (2) %%%%%
\vsp
\noindent
\subsection{The Laplace transform}
% \vskip 0.25truecm
Let
$$ \widetilde f(s) =
{\cal L} \left\{ f(t);s\right\}
 = \int_0^{\infty} \e^{\ds \, -st}\, f(t)\, dt\,, \q
\Re\,(s) > a_f\,,\eqno(2.3a) $$
be the Laplace transform of   a %%% sufficiently well-behaved
function $f(t) \in \L^c(0,T)\,,\, \forall T>0\, $ and let
$$ f(t) =
  {\cal L}^{-1} \l\{ \widetilde f(s);t\r\}
 = \rec{2\pi i}\, \int_{\gamma -i\infty} ^{\gamma + i\infty}
\e^{\ds \, st}\, \widetilde f(s) \, ds\,, \q
\Re\,(s) = \gamma > a_f\,,\eqno(2.3b) $$
with $t>0\,,$
be the inverse Laplace transform\footnote{%%}
%%%%%% FOOTNOTE (4) on the LAPLACE TRANSFORM
A sufficient condition of the existence of the Laplace transform
is that the original function is of exponential
order as $t \to \infty\,. $  This means that some constant
$a_f $	exists such that the product
$ \e^{-a _f t}\,  |f(t)|$ is bounded for all $t$ greater
than some $T\,. $
Then $\tilde{f}(s)$	exists and
is analytic in the half plane $\Re (s) > a_f\,. $
If $f(t)$ is piecewise differentiable, then the formula (2.3b)
holds true at all points where $f(t)$ is continuous and
the (complex) integral in it must be understood in the sense of the Cauchy
principal value.}
%%%%%%% END OF THE FOOTNOTE (4) on the LAPLACE TRANSFORM  %%%%%%%%%
\subsection{The auxiliary functions of Mittag-Leffler type}
The Mittag-Leffler functions, that we denote by $E_\alpha(z)$, $E_{\alpha, \beta}(z)$ are so named  
in honour of G\"osta Mittag-Leffler, the eminent Swedish mathematician, who introduced and 
investigated these functions
 in a series of notes starting from 1903 in the framework of the theory of entire 
 functions\cite{Mittag-Leffler_03,Mittag-Leffler_04,Mittag-Leffler_05}.  
The functions are defined by the series representations, convergent in the whole complex plane $\CC$ 
$$
E_\alpha (z) := \sum_{n=0}^\infty \frac{z^n}{\Gamma (\alpha n+1)}
\,,\q \hbox{Re}(\alpha) > 0\,; \eqno(2.4)$$
$$
E_{\alpha, \beta  } (z ) := \sum_{n=0}^\infty
  \frac{z^n}{ \Gamma(\alpha n+ \beta   )}\,,
\quad \hbox{Re} (\alpha) >0\,,\; \beta  \in \CC\,.
   \eqno(2.5)$$
   Originally Mittag-Leffler assumed only the parameter $\alpha$ and assumed it as positive,
    but soon later the generalization 
   with two complex parameters was considered
   by Wiman\cite{Wiman_05a}.
   In both cases the Mittag-Leffler functions are entire of order 
   $1/\hbox{Re}(\alpha)$. Generally $E_{\alpha,1}(z)= E_\alpha(z)$.
   \vsp
Using their series representations it is easy to recognize
$$
\left\{
\begin{array}{lll}
  & E_{1,1}(z)= E_1(z)  = \e^z\,, \q &E_{1,2}(z)= {\ds\frac{\e^z - 1}{z}}\,,\\
  & E_{2,1}(z^2)= \cosh(z)\,,\q &E_{2,1}(-z^2)= \cos(z)\,, \\
   & E_{2,2}(z^2)= {\ds \frac{\sinh(z)}{z}}\,, \q &E_{2,2}(-z^2)={\ds \frac{\sin(z)}{z}}\,,
   \end{array} \right . \eqno(2.6) $$
and more generally
$$
\left\{  
\begin{array}{ll} 
&E_{\alpha,\beta}(z) + E_{\alpha,\beta}(-z)= 2\, E_{2\alpha,\beta}(z^2)\,,\\
&E_{\alpha,\beta}(z) - E_{\alpha,\beta}(-z)= 2z\, E_{2\alpha,\alpha+ \beta}(z^2)\,.
\end{array}\right . \eqno(2.7)
$$
 We note that in  Chapter 18 of Vol. 3 of the  handbook of the 1955 Bateman Project
\cite{Erdelyi_HTF} devoted to Miscellaneous Functions,
 we find a valuable survey of these functions, which  later were recognized as belonging to the more general class
 of Fox $H$-functions introduced after 1960. 
 \vsp
 For our purposes relevant roles are played by the following
 auxiliary functions of the Mittag-Leffler type on support  $\RR^+$ defined as follows, where $\lambda>0$,
  along with their Laplace transforms
$$ e_\alpha (t;\lambda ) := E_\alpha \left(-\lambda \, t^\alpha \right)
   \div
\frac{s^{\alpha -1}}{ s^\alpha +\lambda}\,,
  \eqno(2.8)$$
  %%%%%%%%%%%%%%%%%
 $$ e_{\alpha,\beta }  (t;\lambda):=
   t^{\beta   -1}\,  E_{\alpha,\beta  } \left(-\lambda \, t^\alpha\right)
   \div
 \frac{s^{\alpha -\beta  }}{ s^\alpha +\lambda}
           \,, \eqno(2.9)$$
	$$ e_{\alpha,\alpha}  (t;\lambda):=
   t^{\alpha  -1}\,  E_{\alpha,\alpha} \left(-\lambda \, t^\alpha\right)
   = \frac{d}{dt} e_\alpha(-\lambda\, t^\alpha)
=    \div
 -\frac{\lambda}{ s^\alpha +\lambda}
\,, \eqno(2.10)$$	   
		   Here we have used the sign $\div$ for the juxtaposition of
a function depending on $t$ with its Laplace transform depending on $s$.
Later we use this sign also for juxtaposition of a function depending on $x$ 
with its Fourier transform depending on $\kappa$. 
\pni {\bf Remark}: We outline that the above auxiliary functions (for restricted values of the parameters)
turn out to be {\it completely monotone} (CM) functions so that they enter in some types of
relaxation phenomena of physical relevance.
\pni
We recall that a  function $f(t)$ is CM in $\RR^+$ if $(-1)^n\, f^{n}(t) \ge 0 $.
The function  $\e^{-t}$ is the prototype of a CM function.
For a Bernstein theorem, more generally they are expressed in terms of a (generalized) real Laplace transform
of a  positive measure
$$
 f(t)= \int_0^\infty \e^{-rt}\, K(r)\, {\mbox{d}}r \,, \q K(r) \ge 0\,.\eqno(2.11)$$
\pni
Restricting attention  to the auxiliary function in two parameters, 
we can prove for $\lambda >0$ that
$$e_{\alpha,\beta}(t;\lambda):=
 t^{\beta-1}\,{\ds E_{\alpha,\beta}\left( - \lambda t^\alpha\right)} \quad \hbox{CM}
 \quad \hbox{iff} \quad  0<\alpha \le \beta \le  1\,. \eqno(2.12)
  $$
  Using the Laplace transform we can prove, following Gorenflo and Mainardi
  \cite{GorMai_CISM97},  that  for $0<\alpha<1$ and $\lambda=1$
$$
 E_\alpha \left(-  t^\alpha \right) \simeq
\left\{
\begin{array}{ll}
 {\ds 1-  \frac{t^{\alpha}}{\Gamma{(\alpha +1)}}
 + \frac{t^{2\alpha}}{\Gamma{(2\alpha +1)}} \cdots} & t\to 0^+\,,\\ 
 {\ds  \frac{t^{-\alpha}}{\Gamma(1-\alpha)}
 - \frac{t^{-2\alpha}}{\Gamma(1-2 \alpha)} \cdots} & t \to +\infty \,,
\end{array} \right. \eqno(2.13)$$
and 
$$
E_\alpha \left(-  t^\alpha \right)=
\int_0^\infty \!\!\e^{-rt}\, K_\alpha(r) \, {\mbox{d}}r \eqno(2.14)$$
 with
 $$ \! K_\alpha (r)
% = -\,\rec{\pi}\,    \hbox{Im}\,
%  \left\{ \left. \frac{s^{\alpha -1}} { s^\alpha +1}\right \vert_{s=r\,\e^{i\pi}} \right\}
 \!=\! \rec{\pi}\,
   \frac{ r^{\alpha -1}\, \sin (\alpha \pi)}
    {r^{2\alpha} + 2\, r^{\alpha} \, \cos  (\alpha \pi) +1}
\!=\! 	\rec{\pi\,r}\,
   \frac{  \sin (\alpha \pi)}
    {r^{\alpha} + 2 \cos  (\alpha \pi) + r^{-\alpha}}
	 > 0\,. \eqno(2.15)
    $$
%%  More generally we can  prove with $\lambda>0 $ that
%% $$  e_\alpha(t;\lambda):= E_\alpha \left(-\lambda   t^\alpha \right)=
%% \int_0^\infty \e^{-rt}\, K_\alpha(r;\lambda) \, {\mbox{d}}r  \,,$$
%% where
%% $$ K_\alpha (r;\lambda) = {\lambda}^{1/\alpha} \, K_\alpha (\lambda^{1/\alpha}\, r)\,.$$
%%%%%%%%%
\subsection{The auxiliary functions of the Wright type}
The Wright function, that we denote by  $W_{\lambda,\mu,}(z)$, is so named in honour
of E. Maitland Wright, the eminent British mathematician, who introduced
and investigated this function
%% (with the restriction $\lambda \ge 0$
in a series of notes starting from 1933 in the framework of the asymptotic theory of partitions,
see  \cite{Wright_33,Wright_35a,Wright_35b}.
The function is defined by the series representation,
convergent in the whole $z$-complex plane $\CC$,
  $$ W_{\lambda ,\mu }(z ) :=
   \sum_{n=0}^{\infty}\frac{z^n}{n!\, \Gamma(\lambda  n + \mu )}\,,
 \q \lambda  >-1\,, \; \mu \in \CC\,. \eqno(2.16)$$   %%% (F.1)
Originally, Wright assumed  $\lambda \ge 0$, and,
only  in 1940 \cite{Wright_40}, he considered %% the case
$-1<\lambda <0$.
 We note that in  Chapter 18 of Vol. 3 of the  handbook of the 1955 Bateman Project
\cite{Erdelyi_HTF} devoted to Miscellaneous Functions,
 we find an earlier analysis of these functions, which, similarly with the Mittag-Leffler functions,  
  were later recognized as belonging to the more general class
 of Fox $H$-functions introduced after 1960. 
However, in that Chapter, presumably for a misprint,  the parameter $\lambda $ 
of the Wright function is restricted to be non negative.
It is possible to prove that the Wright function is entire of order
$1/(1+\lambda)\,, $ hence it is
of exponential type only if $\lambda \ge 0$. For this reason
we propose to distinguish the Wright functions in two kinds according to
$\lambda \ge 0$ ({\it first kind})
and $-1<\lambda<0$ ({\it second kind}). 
Both kinds of functions are related to the Mittag-Leffler function via  Laplace transform pairs:
in fact we have,
see for details the appendix F of the recent book by Mainardi \cite{Mainardi_BOOK10},
for the case $\lambda>0$ (Wright functions of the first kind)
$$   W_{\lambda ,\mu } (\pm r) \,\div \,
    \rec{s}\, E_{\lambda ,\mu }\left( \pm \rec{s}\right) \,,
  \q \lambda > 0\,,  \eqno(2.17)$$ 
and for  the case $\lambda
=-\nu \in (-1,0)$ (Wright functions of the second kind), 
$$     W_{-\nu ,\mu }(-r) \,\div\,
   E_{\nu  , \mu +\nu  }(-s)\,, \q 0<\nu  <1\,. \eqno(2.18)$$
For our purposes relevant roles are played by the following
 auxiliary functions of the Wright  type (of the second kind)
 $$ F_\nu (z) :=   W _{-\nu , 0}(-z)= 
 {\ds \sum_{n=1}^{\infty}
\frac{(-z)^n}{  n!\, \Gamma(-\nu n)}}
\,, \; 0<\nu<1\,, \eqno(2.19)$$
and
$$ M_\nu (z) :=  W _{-\nu , 1-\nu }(-z)
 ={\ds \sum_{n=0}^{\infty}
 \frac{(-z)^n }{  n!\, \Gamma[-\nu n + (1-\nu )]} }\,,
\; 0<\nu<1 \,,\eqno(2.20)$$
interrelated through
$ F_\nu (z) = \nu  \, z \, M_\nu (z )$.
The relevance of these functions  was pointed out by Mainardi in his 
former analysis of the time fractional diffusion equation via  Laplace transform.   
Restricting our attention to the $M$-Wright functions on support  $\RR^+$ 
 we point out  the Laplace transforms pairs
 $$  M_\nu (r) \,\div\,  E_\nu (-s)\,, \q 0<\nu<1\,.  \eqno(2.21)$$
$$ 
   \frac{\nu }{  r^{\nu +1}}\,  M_\nu \left( 1/{r^\nu } \right)\,\div\,
    \e^{\ds \,-s^\nu}\,, \q  0<\nu <1\,. \eqno(2.22)$$
%%%%%%%%%%%%%%%%%%%%%%%%%
$$	
   \frac{1}{  r^{\nu}}\,  M_\nu \left( 1/{r^\nu } \right)\,\div\,
    \frac{\e^{\ds\, -s^\nu}}{s^{1-\nu}}\,, \q  0<\nu <1\,. \eqno(2.23)$$
It was also proved in \cite{Mainardi_BOOK10} that the $M$-Wright function on support  $\RR^+$
is a probability density function (pdf) )that in the literature is sometimes known  as the density related 
to Mittag-Leffler probability distribution. Its absolute moments  of order $\delta>-1$ in $\RR^+$ 
   are finite  and turn out to be
   $$   \int_0^\infty \!\! r^\delta M_\nu(r)\, dx =
   \frac{\Gamma(\delta+1)}{  \Gamma(\nu \delta+1)}\,,
    \q\delta >-1\,,\q 0\le \nu <1\,. \eqno(2.24)$$
	We  point out that in the  limit $\nu \to 1^-$ the function $M_\nu(r)$, for $r \in \RR^+$, tends to the
   Dirac generalized function $\delta(r-1)$.
   \vsp
For our next purposes 
it is worthwhile to introduce the 
{\it function  in two variables} 
  $$M_\nu(x,t):= t^{-\nu}\, M_\nu(xt^{-\nu})\,,\q 0<\nu < 1\,,\q x,t \in \RR^+ \,,\eqno(2.25)$$
  which defines a  spatial probability density in $x$ evolving in 
  time $t$ with self-similarity exponent $H=\nu$.
  Of course for $x\in \RR$ we can consider the symmetric version obtained from 
  (2.25) multiplying by $1/2$ and replacing $x$ by $|x|$.
  %% \vsp 
   Hereafter we provide
  a list of the main properties of this density, 
  which can be derived from the Laplace and Fourier transforms 
  for the corresponding Wright $M$-function 
  in one variable.  %%  presented above.
 \vsp
   From Eq. (2.23) we   derive the Laplace transform of $M_\nu(x,t)$ with respect to $t \in\RR^+$,
   $$\L\left\{M_\nu (x,t); t\to s \right\}= s^{\nu-1}\, \e^{\ds \, -xs^\nu}\,.\eqno(2.26)$$
   %%%%%%%%%%%%%%%%%%%
    From Eq. (2.21) we   derive the Laplace transform of $M_\nu(x,t)$ with respect to $x\in \RR^+$,
	$$\L\left\{M_\nu(x,t); x\to s \right\}= E_{\nu}\left( -s t^\nu \right)\,.\eqno(2.27)$$
	%%%%%%%%%%%%%%%%%%%%%%%%%%%%%%%%%%%%%%%%%%%%%%%%%%%%%%%%%%%%%%%%%%%%%%%%%%%%%%%
	     From the recent book by Mainardi\cite{Mainardi_BOOK10} we   recall 
		 the Fourier transform of $M_\nu(|x|,t)$ with respect to $x\in \RR$,
	$$\F\left\{M_\nu(|x|,t); x\to \kappa \right\}= 2E_{2\nu}\left( -\kappa^2 t^{2\nu} \right)\,,\eqno(2.28)$$
	and, in particular,
	%%%%%%%%%%%%%%%%%%%%%%%%%%%%%%%%%%%%%%%%%%%%%%%%%%%%%%%%%%%%%%%%
$$ \left\{ \begin{array}{ll}
\int_0^\infty \cos (\kappa x) \, M_\nu(x,t) \, dx &= E_{2\nu, 1} (-\kappa^2 \,t^{2\nu})\,,\\
\int_0^\infty \sin (\kappa x) \, M_\nu(x,t) \, dx &= \kappa \,t^\nu \, E_{2\nu, 1+\nu} (-\kappa^2 \,t^{2\nu})\,.
\end{array} \right. \eqno(2.29)$$
 It is worthwhile to note that for $\nu=1/2$ 
 we recover the Gaussian density evolving with time with variance $\sigma^2=2t$
$$ \frac{1}{2\,}M_{1/2}(x,t)= \rec{2\sqrt{\pi}t^{1/2}}\, \e^{\ds\, -x^2/(4t)} \,. \eqno(2.30)$$

 %%%%%%%%%%%%%%%%%%%%%%%%%%%  SECTION 5
% \newpage
\subsection{The L\'evy stable distributions}
\begin{figure}
\begin{center} %% \centering
\includegraphics[width=0.49\textwidth]{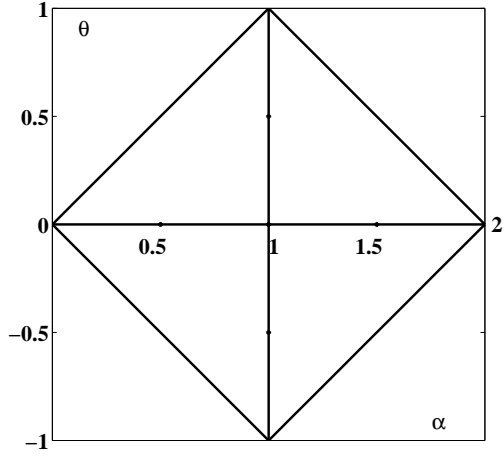}
\end{center} % \vskip 0.25truecm
\vspace{-0.2truecm}
 \caption{The Feller-Takayasu diamond.}
\label{fig:DIAMOND}
%\cen{{\bf Fig. 1} The Feller-Takayasu diamond}
\end{figure}

 The term stable has been assigned by the French  mathematician Paul L\'evy,
 who, in the twenties of the last century, started a systematic research
in order to generalize the celebrated {\it Central Limit Theorem}\index{Central limit theorem} to
probability distributions  with infinite variance.
For stable distributions we can assume the following
 {\sc Definition}:
%% \\
{\it If two independent real random variables
with the same shape or type of distribution are combined linearly with positive coefficients and
the distribution of the resulting random variable has  the same shape,
the common distribution (or its type, more precisely) is said to be
stable}.
\vsp
The restrictive condition of stability enabled L\'evy (and then other authors) to derive
the {\it canonic form} for the characteristic function of the densities of these distributions.
Here we follow the parameterization by Feller\cite{Feller_1971}
 revisited in \cite{GorMai_FCAA98} and  in  \cite{Mainardi-Luchko-Pagnini_FCAA01}.
 Denoting by  $L_\alpha^\theta(x) $ a generic stable density  in $\RR$,
 where $\alpha$ is the {\it index of stability} and
 and $\theta$ the asymmetry parameter, improperly called {\it skewness},
 its characteristic function reads:
$$\left\{
\begin{array}{ll}
{\ds L_\alpha^\theta(x)}  \div {\ds \widehat{L}_\alpha ^\theta(\kappa)}  =
  {\ds \exp \left[- \psi_\alpha ^\theta(\kappa )\right]} \,, \;
   {\ds \psi_\alpha  ^\theta(\kappa )} =
   {\ds |\kappa|^{\ds \alpha } \, \e^{\ds  i (\sgn \kappa)\theta\pi/2}}\,,\\
 0<\alpha  \le 2\,, \;
 |\theta| \le  \,\hbox{min}\, \{\alpha  ,2-\alpha  \}\,.
\end{array} \right. \eqno(2.31) $$
 %%%%%%%%%%%%% \newpage
 %%%%%%%%%%%
 %%\vsp
 We note that the allowed region for the  real
parameters $\alpha  $ and $\theta$
turns out to be
 a {diamond} in the plane $\{\alpha  , \theta\}$
with vertices in the points
$(0,0)\,, \,(1,1)\,, \, (1,-1) \,,\,(2,0)$,
that we call the {\it Feller-Takayasu diamond}\index{Feller-Takayasu diamond},
see Figure~\ref{fig:DIAMOND}. 
For values of $\theta$ on the border of the diamond
(that is $\theta = \pm \alpha $ if $0<\alpha  < 1$, and $\theta = \pm (2-\alpha )$ if $1<\alpha  <2$)
we obtain the so-called {\it  extremal stable densities}.

%%%%%%%%
 \vsp
 We note the {\it symmetry relation}
$L_\alpha ^\theta (-x)=	L_\alpha ^{-\theta} (x)$, so that a stable density with $\theta=0$ is symmetric.
  \vsp
  Stable distributions have noteworthy properties on which the interested reader
  can be informed from the relevant existing literature.
 Here-after we recall some  peculiar {\sc Properties}:
 \pni
- {\it Each  stable density $L_\alpha^\theta$ possesses a {\it domain of attraction}}, 
see \eg \cite{Feller_1971}.
\pni
- {\it Any stable  density  is {unimodal} and indeed {bell-shaped}}, \ie
its $n$-th derivative has exactly $n$ zeros in $\RR$, see \cite{Gawronski_1984}.
 \pni
- {\it The stable distributions are  {self-similar} and  {infinitely divisible}}.
\vsp
These properties derive from the canonic form (2.31) through the scaling property of
the Fourier transform.
\pni
{\it Self-similarity}\index{Self-similarity} means
$$L_\alpha ^\theta (x,t) \div \exp \left [-t \psi _\alpha  ^\theta(\kappa)\right]
\Longleftrightarrow  L_\alpha ^\theta (x,t)  = t^{-1/\alpha }\,L_\alpha ^\theta (x/t^{1/\alpha} ) ]\,,
 \eqno (2.32)  $$
where $t$ is a positive parameter.
If $t$ is time, then  $L_\alpha ^\theta (x,t)$ is a spatial density evolving in time with self-similarity.
 \pni
   {\it Infinite divisibility}\index{Infinite divisibility} means that  for every positive integer
$n$,  the characteristic function  can be expressed as the $n$th power of some characteristic function,
so that
%% Equivalently we can say that for every positive integer $n$
any stable distribution can be expressed as the	$n$-fold convolution of a
stable distribution of the same type.
Indeed, taking in (2.31) $\theta=0$, without loss of generality, we have
$$\e^{-t|\kappa|^\alpha } = \left[\e^{-(t/n)|\kappa|^\alpha }\right]^n
\Longleftrightarrow  L_\alpha ^0 (x,t)  =
\left[L_\alpha ^0 (x,t/n)\right]^{*n} \,, \eqno(2.33)$$
where
$$\left[L_\alpha ^0 (x,t/n)\right]^{*n} :=
L_\alpha ^0 (x,t/n) * L_\alpha ^0 (x,t/n) *  \dots * L_\alpha ^0 (x,t/n) \eqno(2.34) $$
is the multiple Fourier convolution in $\RR$ with $n$
identical terms.
%%%%%%%%%%
% \newpage
\vsp
Only in special cases the inversion of the Fourier transform in (2.31)
can be carried out using standard tables, and provides well-known probability distributions.
\vsp
For $\alpha  =2$ (so $\theta =0$), we recover the  {\it Gaussian pdf}, that turns out to be the
only stable density with finite variance, and more generally with finite  moments of any order
$\delta \ge 0$. In fact
$$ L_2^0(x) = \frac{1}{2\sqrt{\pi}}\e^{\,\ds -x^2/4} %%%\,\div \, \e^{\,\ds -\kappa^2}
\,.\eqno(2.35)$$
All the other stable densities have finite absolute moments  of order $\delta \in [-1, \alpha )$
as we will later show.
\vsp
%%  with $-1<\delta <\alpha $.
For $\alpha  =1 $  and  $|\theta| <1$, we get
$$L_1^\theta (x) =
 \frac{1}{\pi} \, \frac{  \cos (\theta \pi/2)}
 {[x+   \sin (\theta \pi/2)]^2 +[ \cos (\theta \pi/2)]^2 }\,,\eqno(2.36)$$
which  for $\theta=0$ includes the  {\it Cauchy-Lorentz pdf},
$$ L_1^0(x) = \frac{1}{\pi} \frac{1}{1+x^2}   %%%\,\div \, \e^{\,\ds -|\kappa|}
\,.\eqno(2.37)$$
In the limiting cases   $\theta = \pm 1$ for $\alpha =1$ we obtain
 the {\it  singular Dirac pdf's}
 $$ L_1^{\pm 1}(x)=\delta(x \pm 1)\,.\eqno(2.38)$$
\vsp
In general, we must recall the power series expansions
provided in \cite{Feller_1971}.
We restrict our attention to $x>0$
since the evaluations for $x<0$  can be obtained using the symmetry relation.
%%%%%%%%%%%%%%
% \newpage
%
The convergent expansions of $L_\alpha ^{\theta} (x)$  ($x>0$) turn out to be;
 \par
for $ 0<\alpha   <1\,,\q |\theta| \le \alpha   \,:$
 $$L_\alpha  ^\theta (x) =
\frac{1}{\pi\,x}\,  \sum_{n=1}^{\infty}
   (-x^{-\alpha  })^n \, 
   \frac{\Gamma (1+ n\alpha  )}{n!}\,
  \sin \left[ \frac{n\pi}{2}(\theta -\alpha  )\right]\,;
    \eqno(2.39)  $$
\par
for $ 1<\alpha   \le 2\,, \q |\theta| \le 2-\alpha  \,:$
$$ L_\alpha  ^\theta (x)=
\frac{1}{\pi\,x}\,  \sum_{n=1}^{\infty}
   (-x)^{n} \, \frac{\Gamma (1+ n/\alpha  )}{n!}\,
  \sin \left[\frac{ n\pi}{ 2\alpha }(\theta -\alpha  )\right]\,.
 \eqno(2.40) $$
%%%%%%%%%
From the series in (2.39) and the  symmetry relation
we note that {\it the extremal stable densities for $0<\alpha   <1$ are
unilateral}, precisely vanishing for $x>0$ if $\theta =\alpha $,
vanishing for $x<0$ if $\theta =-\alpha $.
In particular the unilateral extremal densities $L_\alpha ^{-\alpha }(x)$ with $0<\alpha <1$
have support  $\RR^+$ and Laplace transform
$\exp (-s^\alpha )$. For $\alpha=1/2$ we obtain the so-called {\it L\'evy-Smirnov} $pdf$:
$$ L_{1/2}^{-1/2} (x) =
  \frac{\,x^{-3/2}}{2\sqrt{\pi}}\, \,\e^{\ds\, - 1/(4x)}\,,
 \q x \ge 0	\,. \eqno(2.41)$$
%%%%%%%%%%%%%
It is worth to note that the Gaussian $pdf$ (2.35) and the L\'evy-Smirnov $pdf$ (2.41)
 are well known in the treatment of the Brownian motion: the former as the spatial density on an infinite 
 real line, the latter    as the first passage time density on a semi-infinite line, see e.g.\cite{Feller_1971}. 
\vsp
As a consequence of the convergence of the series
in (2.39)-(2.40) and of the symmetry relation
%% $L_\alpha^\theta (-x)=	L_\alpha^{-\theta} (x)\,,$
we recognize that
the stable $pdf$'s with $1< \alpha \le 2$
are  entire functions, whereas
 with $0< \alpha <1$
have the form
$$L_\alpha^\theta (x) = \begin{cases}
   (1/x) \,\Phi_1(x^{-\alpha })
	& \hbox{for} \; x>0\,, \\
  (1/|x|) \,\Phi_2(|x|^{-\alpha })
	& \hbox{for} \; x<0\,,
	\end{cases}
	  \eqno(2.42)		 $$
where $\Phi_1(z)$ and $\Phi_2(z)$ are distinct {entire functions}. 
The  case $\alpha =1$ with  $|\theta|< 1$ must be considered  in the limit for
$\alpha \to 1$ of  (2.39)-(2.40), because the corresponding
series reduce to power series akin with geometric series
in $1/x$ and $x$, respectively, with a finite radius of convergence.
The corresponding stable densities are no longer represented by
entire functions, as can be noted directly from their explicit expressions (2.36)-(2.37).
\vsp
From a comparison between  the series expansions in (2.39)-(2.40) and in (2.19)-(2.20),
we recognize that for $x>0$ our
{\it auxiliary functions of the Wright type are related to the extremal stable densities
as follows}, see\cite{Mainardi-Tomirotti_GEO97},
$$ L_\alpha ^{-\alpha } (x)
 =  \frac{1}{x}\,  F_\alpha  (x^{-\alpha }) =
 \frac{\alpha} {x^{\alpha  +1}}\,  M_\alpha (x^{-\alpha }) \,,\q
 0<\alpha  <1\,, \eqno(2.43)$$
$$ L_\alpha ^{\alpha  -2}(x)
 = \frac{1}{x}\,  F_{1/\alpha }(x) =
 \frac{1}{\alpha }\,  M_{1/\alpha }(x) \,,\q
 1<\alpha  \le 2\,.\eqno(2.44)$$
 In the above equations, for $\alpha=1$, the skewness parameter turns out to be $\theta =-1$,
 so  we get the singular limit
 $$L_1^{-1}(x)= M_1(x)= \delta(x-1)\,.\eqno(2.45)$$
%%%%%%%%%% \vsp
We do not provide here the asymptotic representations of the stable densities referring the 
interested reader to\cite{Mainardi-Luchko-Pagnini_FCAA01}. However, based on asymptotic representations, 
we  can state the following:
 For $0<\alpha <2$ the stable densities  exhibit
{\it fat tails} in such a way that their  absolute moment
of  order $\delta$ is finite only if $-1 < \delta <\alpha $.
More precisely, one can show that for non-Gaussian,
not extremal, stable densities
the asymptotic decay of the tails  is
$$ L_\alpha^\theta (x )= O\left(|x|^{-(\alpha +1)}\right)\,, \q
	       x \to \pm \infty\,. \eqno(2.46)$$
For the extremal densities with $\alpha \ne 1$
this is valid only for one
tail (as $|x|\to \infty$), the other (as $|x|\to \infty$) being of exponential order.
%% For $0<\alpha <1$ we have one-sided $pdf$'s:
% for $\theta = -\alpha $ the support is $\RR^+$ and the $pdf$ tends exponentially to zero
% as $x \to 0^+\,;$
%for $\theta = +\alpha $ the support is $\RR^-$and the $pdf$ tends exponentially to zero
% as $x \to 0^-\,.$
For $1<\alpha <2$ the extremal $pdf$'s are two-sided and exhibit
an exponential left tail  (as $x \to -\infty)$
if $\theta  =+(2-\alpha)\,,  $
or  an exponential right tail  (as $x \to +\infty $)
if $\theta  =-(2-\alpha)\,.$
%% less than  the index of stability $\alpha \,,$
%% \vsp
Consequently, the Gaussian $pdf$ is the unique
 stable density with finite variance.
Furthermore,  when $0<\alpha \le 1$,
the first absolute moment  is infinite  so
we should use the  median instead of the non-existent expected value
in order to characterize the corresponding $pdf$.
%%%%%
\vsp
Let us also recall a relevant identity between stable densities
 with index $\alpha $ and $1/\alpha$  
  (a sort of reciprocity relation) pointed out in\cite{Feller_1971},
that is, assuming $x>0$,
$$ \rec{x^{\alpha +1}}\, L_{1/\alpha}^\theta (x^{-\alpha} )
  =L_\alpha^{\theta^*} (x)\,,  \;1/2\le \alpha\le 1\,,\;
  \theta ^*=\alpha(\theta +1)-1 \,. \eqno(2.47)$$
 The condition $1/2\le \alpha \le 1$ implies $1\le  1/\alpha \le 2$. A check shows
 that $\theta^*$ falls within the prescribed range
$|\theta ^*|\le\alpha$ if $|\theta |\le 2-1/\alpha $.
We leave as an exercise for the interested reader the verification of this reciprocity relation
in the limiting cases $\alpha=1/2$ and $\alpha=1$.
%% \vsp
%%%%%%%%%%%%%%%%%%%%%%%%%
For more details on L\'evy stable
densities we refer the reader to  specialized treatises,
as\cite{Feller_1971,Janicki-Weron_94,Samo-Taqqu_94,Sato_99,Uchaikin-Zolotarev_99,Zolotarev_86},
where different notations are adopted.
%% however, the notations for the asymmetry
%% parameter are different form ours here adopted.
We like to refer also to the  %% noteworthy
1986 paper by Schneider\cite{Schneider_LNP86},
where he first provided  the   Fox $H$-function
representation of the  stable distributions (with $\alpha  \ne 1$)
and	to the 1990 book by Takayasu\cite{Takayasu_FRACTALS},
where he first gave the diamond  representation in the plane
$\{\alpha ,\theta\}$.
%%%%%%%%%%%%
\section{The Space-Time Fractional Diffusion}
%%%%%%%%
 We now consider the Cauchy problem for the
(spatially one-dimensional) {\it space-time fractional diffusion} (STFD) equation.
$$  {\, _t}D_{*}^{\, \beta }\, u(x,t)
 \, = \,
 {\, _x}D_{ \theta}^{\,\alpha} \,u(x,t)\,,
\quad  u(x,0) = \delta (x)\,, \quad x \in \RR,\quad t \ge 0\,, \eqno(3.1) $$
where
%% $ -\infty<x<+\infty\,,$ $\,t\ge 0\,, $
   \{$\alpha \,,\,\theta\,,\, \beta $\} are real parameters
 restricted to the ranges
$$ 0<\alpha\le 2\,,\quad  |\theta| \le \min \{\alpha, 2-\alpha\}\,,
  \quad 0<\beta\leq 1\,.\eqno(3.2) $$
Here
 $ {\,_x}D_{\,\theta}^{\,\alpha}$
denotes
the {\it  Riesz-Feller fractional derivative}
of order $\alpha $ and
skewness $\theta$,   acting on the space variable $x$,
and 
${ \,_t}D_*^{\,\beta}  $
denotes  the
{\it  Caputo fractional derivative}
of order $\beta $, acting on the time variable $t$.
We recall the definitions of these fractional derivatives
%%For the sake of the reader's convenience  here we present an
%% introduction to the {\it Riesz-Feller} and {\it Caputo} fractional
based on their representation in the Fourier and Laplace transform domain, respectively.
So doing we avoid
%% the more rigorous approach to the fractional derivatives
%% that is  based on the inversion of the fractional integrals,
the subtleties lying in the inversion of the corresponding fractional integrals,
see e.g. the 2001 survey by Mainardi et al.\cite {Mainardi-Luchko-Pagnini_FCAA01}.
For general information on fractional integrals and derivatives 
we recommend the books\cite{Kilbas-Srivastava-Trujillo_BOOK06,Podlubny_99,SKM_93}.  
%%%%%%%%%
\subsection{The Riesz-Feller  space-fractional derivative}
We define the {\it Riesz-Feller}
derivative as
the pseudo-differential operator whose symbol
%% $- \psi_\alpha ^\theta(\kappa ) $
is the logarithm of the characteristic function of a
general {\it L\'evy strictly stable} probability density
with {\it index of stability} $\alpha $ and asymmetry parameter
$\theta$ (improperly called {\it skewness}).
As a consequence of Eq. (2.31), for a sufficiently well-behaved function $f(x)$, we
define the {\it Riesz-Feller} space-fractional derivative
of order $\alpha  $ and skewness $\theta$ via the Fourier transform
$$ 
\left\{
\begin{array}{ll}
{\cal F} \left\{\, _xD_\theta^\alpha\, f(x);\kappa \right\} =
  - \psi_\alpha ^\theta(\kappa ) \,
  \, \widehat f(\kappa) \,, \;
   \psi_\alpha ^\theta(\kappa ) =
|\kappa|^\alpha \, i^{\ds  i \theta \sgn \kappa)}
\,,\\
 0<\alpha  \le 2\,, \q
 |\theta| \le  \,\hbox{min}\, \{\alpha ,2-\alpha \}\,.
\end{array} \right.
\eqno(3.3)$$
Notice that
$i ^{\,\ds \theta \,\sgn \kappa}= \exp [ i \,(\sgn \kappa)\,\theta\,\pi/2 ]$.
For $\theta=0$	we have a symmetric operator
with respect to $x\,,$ which can be interpreted as
%%% the negative of the
%% $\alpha/2$ power of the (positive definite) operator
%% $ - D^2 = -{d^2\over dx^2}\,,$  namely
$$
%%% {d^\alpha  \over d|x|^\alpha} =
  _xD_0^\alpha = - \left ( - \frac{d^2} {dx^2}  \right) ^{\alpha /2}\,,
\eqno(3.4) $$
as can be formally deduced by writing
$- |\kappa| ^\alpha = - (\kappa^2)^{\alpha /2} \,.   $
 We thus recognize that the  operator $D_0^\alpha $ is related
to a power of the positive definitive operator
  $\, - _xD^2= -\frac{d^2}{ dx^2}$
and  must not be confused with a
power  of the first order differential operator
  $\,_xD= \frac{d} {dx}$
for which the symbol is $-i\kappa \,. $
An alternative illuminating  notation for the symmetric
fractional derivative	is due to Zaslavsky\cite{Saichev-Zaslavsky_97}, and reads
$$ _xD_0^\alpha =     \frac{d^\alpha}  { d |x|^\alpha}\,.
\eqno(3.5) $$
 For $0<\alpha <2$ and
$|\theta| \le \hbox{min} \,\{\alpha, 2-\alpha \} $
the {\it Riesz-Feller} derivative  reads
$$
\begin{array}{ll}
_xD_\theta^\alpha \,f(x)
 &= {\ds {\Gamma(1+\alpha) \over \pi }} \,
 \left\{ \sin \,[(\alpha+\theta) \pi/2] \,
 {\ds \int_0^\infty
 {f(x+\xi)- f(x)  \over {\xi}^{1+\alpha}}\, d \xi} \right.
\\
 &+ \left.{ \sin \,[(\alpha-\theta) \pi/2]} \,
 {\ds \int_0^\infty
 {f(x-\xi)- f(x)  \over {\xi}^{1+\alpha}}\, d \xi} \right\}\,.
 \end{array}
 \eqno(3.6)$$
%%%%
%% \newpage
%%%%%%%%
\subsection{ The Caputo fractional derivative}
% \vskip 0.25truecm
For a sufficiently well-behaved function $f(t)$ we
 define the {\it Caputo} time-fractional derivative  of
order $\beta$ with $0 <\beta \le 1$ 
through its Laplace transform
%%% $$ {\cal L} \l\{ \,_tD_*^\beta \,f(t) ;s\r\} =
%% s^\beta \, \tilde f(s) - s^{\beta-1}\, f(0^+)\,, \q 0<\beta \le 1\,.  %%
$$ {\cal L} \left\{ _tD_*^\beta \,f(t) ;s\right\} =
      s^\beta \,  \widetilde f(s)
   -	s^{\beta  -1}\, f(0^+) \,,
  \q 0<\beta	\le 1 \,. \eqno(3.7)$$
This leads us to define
$$
    _tD_*^\beta \,f(t) :=
\left\{ \begin{array}{ll}
    {\ds \rec{\Gamma(1-\beta )}}\,{\ds\int_0^t
 {\ds {f^{(1)}(\tau)\, d\tau \over (t-\tau )^{\beta}}}} \,,
  & \; 0<\beta  <1\,, \\
     {\ds {d\over dt}} f(t)\,,
    & \; \beta  =1\,. \\
\end{array}\right.	
   \eqno(3.8) $$
   For the essential properties of the Caputo derivative, 
   see\cite{Diethelm_10,GorMai_CISM97,Podlubny_99}. 
%%%%%%%%%%%%%%%%%%%%%%
\subsection{The fundamental solution of the space-time fractional diffusion equation}
Let us note that the solution $u(x,t)$  of the
 Cauchy problem (3.1)--(3.2), known as the {\it Green function} or fundamental solution
of the space-time fractional diffusion equation,
is a
probability density in the spatial variable $x$, evolving in time $t$.
In the case $\alpha =2$ and $\beta =1$  we recover
the standard diffusion equation for which the fundamental solution
is the Gaussian density with variance $\sigma ^2 =2t$.
Sometimes, to point out the parameters, we may denote the fundamental solution as
 $$ u(x,t) = G_{\alpha,\beta}^\theta (x,t)\,, \eqno(3.9)$$
%% \vsp
For our purposes let us here confine ourselves to recall the
representation in the Laplace-Fourier domain of the (fundamental) solution
 as it results from the application of the transforms
of Laplace and Fourier to Eq. (3.1). Using $\widehat \delta (\kappa ) \equiv 1$
we have:
 $$  s^{\, \ds \beta}\,\widehat{\widetilde{u}}(\kappa ,s) - s^{\, \ds \beta -1}
    = -|\kappa|^{\, \ds \alpha} \,  i ^{\,\ds \theta \,\sgn \kappa}
   \, \widehat{\widetilde{u}}(\kappa ,s) \,, 
$$
hence
$$  \widehat{\widetilde{u}}(\kappa ,s) = \widehat{\widetilde{G_{\alpha,\beta}^\theta}}(\kappa ,s)
    =  \frac{ s^{\, \ds \beta -1}}
{s^{\,\ds \beta} + |\kappa |^{\, \ds \alpha} \,
i^{\,\ds \theta \,\sgn \kappa} }\,.
   \eqno(3.10)$$
For explicit expressions and plots of  the fundamental solution of (3.1)
in the space-time domain
we refer the reader to Mainardi, Luchko and Pagnini\cite{Mainardi-Luchko-Pagnini_FCAA01}.
There, starting from the fact  that the Fourier transform
of the fundamental solution can be written as a Mittag-Leffler function
with complex argument,
$$\widehat{u}(\kappa ,t)=
\widehat{G_{\alpha,\beta}^\theta}(\kappa ,t)=
E_\beta \left(-|\kappa|^{\, \ds \alpha} \,  i ^{\,\ds \theta \,\sgn \kappa}\,t^\beta\right) 
\,,\eqno(3.11)$$
 these authors
have  derived a Mellin-Barnes integral representation
of $u(x,t)= G_{\alpha,\beta}^\theta(x ,t)$  with which they have proved the non-negativity
of the solution for values of the parameters
$\{\alpha,\, \theta, \,  \beta \}$ in the range (3.2)
and analyzed the evolution in time of its  moments.
The representation of $u(x,t)$ in terms of Fox $H$-functions
can be found  in  Mainardi, Pagnini and Saxena\cite{Mainardi-Pagnini-Saxena_JCAM05}, see also
Chapter 6 in the recent book by Mathai, Saxena and Haubold\cite{Mathai-Saxena-Haubold_BOOK-H-2010}.
%%%%%%%%%%%%%%
We note, however,  that the solution of the STFD Equation (3.1) and
its variants has been investigated by several authors;
let us only mention some of them%%
 \cite{Barkai_PRE01,Barkai_ChemPhys02,Fulger-et-al_PRE08,Germano-et-al_PRE09,%
 Hilfer_HONNEF08,Hilfer_FRACTALS95,Hilfer-Anton_PRE95,M3_PRE02sub_PRE02,%
 M3_PRE02sol_PRE02,Meerschaert-Scheffler_04,%%
 Metzler-Barkai-Klafter_EPL99,Metzler-Klafter_JPCB00,%%
 Metzler-Klafter-Sokolov_PRE98,Metzler-Klafter_PhysRep00,%
 Metzler-Klafter_JPhysics04,Saichev_PhysA05,Scalas_PRE04},
where the connection with the $CTRW$ was also pointed out.
\vsp
In particular
the fundamental solution
 for the {\it space fractional diffusion} 
  $\{0<\alpha <2, \, \beta=1\}$ 
is expressed in terms
of a stable density of order $\alpha$ and skewness $\theta$,
$$
 G_{\alpha,1}^\theta(x,t) = t^{-1/\alpha }\,L_\alpha ^\theta(x/t^{1/\alpha })\,,
\q -\infty<x<+\infty\,, \q t\ge 0
\,. \eqno(3.12) $$
%%%%%%%%%%%
whereas for the {\it time fractional diffusion}
$\{\alpha =2, \, 0<\beta<1\}$ 
in terms of a (symmetric) $M$-Wright function of order $\beta/2$,
$$ G_{2,\beta}^0 (x,t) =
{1\over 2} t^{-\beta /2}\, M_{\beta/2} \left(|x|/t^{\beta /2}\right)\,,
\q -\infty < x <+\infty\,,  \q t\ge 0\,,
\eqno(3.13)$$
For  the {\it standard diffusion}
$\{\alpha =2\,, \,\beta =1\}$ we recover
the Gaussian density
$$ G_{2,1}^0 (x,t) = \frac{1}{2\sqrt{\pi t}} \,\e^{\ds -x^2/(4t)} =
t^{-1/2}\,L_2 ^0 (x/t^{1/2})= \frac{1}{2} t^{-1 /2}\, M_{1/2} \left(|x|/t^{1 /2}\right)
\,.$$ 
%%%%%%%%%%%%%%%%%%%%%
\vsp
Let us finally recall that the $M$-Wright function does appear also in the fundamental solution of 
the {\it rightward time fractional drift equation},
%% Writing this equation in non-dimensional form and adopting  the Caputo derivative we have 
$$_tD^{\beta}_* u(x,t) = -  
\frac{\partial }{\partial x} u(x,t) \,,\q -\infty<x< +\infty\,, \; t\ge 0\,. \eqno(3.14)
$$
% where $0<\beta<1$ and $u(x,0^+) = u_0(x)$.
% When $u_0(x)= \delta(x)$ we  obtain the fundamental solution (Green function) that we denote by.
Denoting by $\G_\beta^*(x,t)$ this fundamental solution, we have  
$$\G_\beta^*(x,t) = 
 \left\{
  \begin{array}{ll}
  {\ds t^{-\beta}\, M_\beta\left(\frac{x}{t^\beta}\right)}\,, & x>0\,,\\
  0\,, & x<0 \,,
  \end{array}
  \right.
  \eqno(3.15) 
$$
that for $\beta=1$ reduces to the right running pulse $\delta(x-t)$ for $x>0$. For details
see\cite{Gorenflo-Mainardi_HONNEF08,Mainardi-Mura-Pagnini_IJDE10}.
%%%%%%%%%
\subsection{ Alternative forms of the space-time fractional diffusion equation}
We note that in the literature  there exist other forms alternative and  equivalent to Eq. (3.1) 
with  initial condition 
 $u(x,0) = u_0(x)$ including the case $u_0(x)= \delta(x)$.
 For this purpose we must briefly recall the definitions of fractional integral and fractional derivative
 according to Riemann-Liouville.
 \vsp
  The Riemann-Liouville fractional integral for a sufficiently well behaved function $f(t)$ ($t\ge 0$) 
 is defined for any order $\mu>0$ as   
$$  \,_tJ^\mu f(t) :=  \frac{1}{\Gamma(\mu)}\int_0^t (t-\tau)^{\mu-1} f(\tau) \, d\tau\,.
\eqno (3.16)$$
We note the convention $\,_tJ^0 = I$ (Identity) and the semigroup property
$$ \,_tJ^\mu \, _tJ^\nu =  \,_tJ^\nu \, _tJ^\mu = \,_tJ^{\mu + \nu}\,, \q
\mu \ge 0\,,\; \nu\ge 0\,.\eqno (3.17)$$
The fractional derivative of order $\mu >0$ in the {\it Riemann-Liouville} 
sense  is defined as the operator
$\, _tD^\mu$ which is the
left inverse of
the Riemann-Liouville integral of order $\mu $
(in analogy with the ordinary derivative), 
$$ \,_tD^\mu \, _tJ^\mu  = I\,, \q \mu >0\,. \eqno(3.18) $$
If $m$ denote the positive integer such that 
such that $m-1 <\mu  \le m\,,$  we recognize from Eqs. (3.16) -- (3.18):
$$\, _tD^\mu  \,f(t) :=  \,_tD^m\, _tJ^{m-\mu}  \,f(t)\,. \eqno (3.19)$$
Then,  restricting our attention to a order $\beta$ with $0<\beta \le 1$
(namely $m=1$) the corresponding Riemann-Liouville fractional derivative turns out
$$  
 \,_tD^\beta  \,f(t) = %% \, _tD^m\, _tJ^{m-\mu}  \,f(t) = \,
 \,
 \left\{
  \begin{array}{ll}
  {\ds \frac{d} {dt}}\left[
  {\ds \rec{\Gamma(1-\beta )}\int_0^t
    \frac{f(\tau)\,d\tau}  { (t-\tau )^{\beta}} }\right]\, ,
 &  \;0 <\beta  < 1, \\
   {\ds \frac{d} {dt} f(t)} \,.
&  \;\beta =1.
\end{array}
\right.
\eqno (3.20)$$
Then we get the relationship among the Caputo fractional 
derivative with  the classical Riemann-Liouville fractional integral and derivative:
$$\,_tD_*^\beta  \, f(t) := \,_tJ^{1-\beta}\, _tD^1 f(t) =
   \,_tD^\beta \left[ f(t) - f(0)\right] = \,_tD^\beta f(t) - \frac{f(0)}{\Gamma(1-\beta)\,t^\beta}
    \,,\eqno(3.21)  $$
and,  as a consequence,  the equivalence of (3.1) with the following  problems
 $$  u(x,t) -u(x,0) = \,_tJ^\beta \, _xD_{ \theta}^{\,\alpha} \,u(x,t)\,,
\quad  u(x,0) = \delta (x)\,, \eqno(3.22)$$
 $$\frac{\partial}{\partial t} u(x,t) = \,_tD^{1-\beta}\, _xD_{ \theta}^{\,\alpha} \,u(x,t)\,,
\quad  u(x,0) = \delta (x)\,. \eqno(3.23)$$
%%%%%%%%%%%%%%%%%
%\newpage
%%%%%%%%%%
\section{ Analytic and stochastic   pathways to subordination in space-time fractional diffusion}
Our starting key-point to introduce the analytical and stochastic approaches to subordination
in space-time fractional diffusion processes is 
   the fundamental solution of the space-time fractional diffusion equation
in the Laplace-Fourier domain  given by (3.10). 
%% For reader's convenience we recall it: 
%% $$  \widehat{\widetilde{u}}(\kappa ,s) =  \frac{ s^{\, \ds \beta -1}}
%% {s^{\,\ds \beta} + |\kappa |^{\, \ds \alpha} \,i^{\,\ds \theta \,\sgn \kappa} }\,. $$
   %%%%%%%%%%%%%%%
   \subsection{The analytical interpretation via operational time.}
Separating variables in (3.10) and using the trick to write ${1/(z+a)}$ 
for $\Re (z+a) >0$ as a Laplace integral
$$  \frac{1}{z+a} = \int_0^\infty \e^{-z\rho}\, \e^{-a\rho}\, d\rho 
$$
we have, identifying $\rho  := t_*$  as {\it operational  time}, the following  instructive expression for (3.10): 
$$  \widehat{\widetilde{u}}(\kappa ,s) = \int_0^\infty
\! \left[\exp{\left(-t_* |\kappa |^{\alpha} i^{\theta \,\sgn {\ds\kappa}}\right)}\right]
\, \left[s^{\beta-1}\, \exp \left(-t_* s^\beta\right)\right]\, dt_*\,.\eqno(4.1)$$
We note that the first factor in (4.1) 
$$  \widehat f_{\alpha, \theta}(\kappa, t_*) := \exp{\left(-t_* |\kappa |^{\alpha} i^{\theta \,
\sgn {\ds \kappa}}\right)}\eqno(4.2)$$
is the Fourier transform of a skewed stable density
in $x$, evolving in operational time $t_*$, of
a  process $x=y(t_*)$ along the real axis $x$ happening in operational time $t_*$, that we write as
$$   f_{\alpha,\theta}(x,t_*)\ = t_*^{-1/\alpha}\,L_\alpha^\theta\left(x/t_*^{1/\alpha}\right)
 \,.  \eqno(4.3)$$
% Its Fourier-Laplace transform is
% $$ \widehat{\widetilde\xi}(\kappa,s_*)=
% \frac{1}{s_*+ |\kappa|^\alpha\, i^{\theta \,\sgn \kappa}}\,.$$
We can interpret the second factor
$$\widetilde q_\beta (t_*,s) := s^{\beta-1} \exp (-t_*s^\beta)\eqno(4.4)$$
as Laplace representation of the probability density in $t_*$ evolving in $t$ of a %% {\it Mittag-Leffler}
 process
$t_* =t_*(t)$, generating the operational time $t_*$
 from the physical time $t$, that is expressed via a fractional integral of a skewed L\'evy density as
$$   q_\beta(t_*,t) = t_*^{-1/\beta} \, _tJ^{1-\beta} \, L_\beta^{-\beta} (t/t_*^{1/\beta})
= t^{-\beta}\, M_\beta(t_*/t^\beta)
\,,\eqno(4.5)$$
see Eq. (2.26).
To prove that  $q_\beta(t_*,t)$ (surely positive for $t_*>0$) is indeed a probability density
we must further prove that is normalized,
${\ds \int_{t_*=0}^\infty q_\beta(t_*, t)\,dt_* =1}\,.$
For this purpose it is sufficient to prove that its Laplace transform with respect to $t_*$ 
is equal to 1 for $s_*=0$.
To get this Laplace transform $\widetilde q_\beta(s_*,t)$ we proceed as follows.
Starting from    the known Laplace transform with respect to $t$,
$$  \widetilde q_\beta(t_*,s)= s^{\beta-1} \exp (-t_*s^\beta)\,, \eqno(4.6)$$
we apply a second Laplace transformation  with respect to $t_*$ with parameter $s_*$ to get
$${\widetilde{\widetilde q}}_\beta (s_*,s) = \frac{s^{\beta-1}}{s_*+s^\beta}\,, \eqno(4.7)$$
so, by inversion with respect to $t$
$$ \widetilde q_\beta(s_*,t) =  \int_{t_*=0}^\infty \e^{-s_*t} q_\beta(t_*, t) \,dt_* = E_\beta(-s_*t^\beta)\,,
\eqno(4.8)$$
and setting $s_*=0$
$$ \int_{t_*=0}^\infty q_\beta(t_*, t)\, dt_*=E_\beta(0) =1\,.\eqno(4.9)$$
%%%%
Weighting the density of $x=y(t_*)$ with the density of $t_*=t_*(t)$ over
$0\le t<\infty$ yields  the density $u(x,t)$
in $x$ evolving with time $t$.
\vsp
In physical variables $\{x,t\}$, using Eqs. (4.1) -- (4.5), we have the {\it subordination integral formula}
$$ u(x,t)= \int_{t_*=0}^\infty\!\! f_{\alpha,\theta}(x,t_*)\, q_\beta (t_*,t)\, dt_*\,,\eqno(4.10)$$
where $ f_{\alpha,\theta}(x,t_*)$ (density in $x$ evolving in $t_*$) refers to the process 
$x=y(t_*)$ ($t_* \to x$) generating in ``operational time" $t_*$ the spatial position $x$,
and $q_\beta(t_*,t)$ (density in $t_*$ evolving in $t$) refers to the process 
$t_*=t_*(t)$ ($t \to t_*$) generating from physical time $t$ the ``operational time" $t_*$.
\vsp
Our aim is to construct a process $x=x(t)$ whose probability density is $u(x,t)$, density in $x$, evolving in $t$..
We will soon find justification for denoting the variable of integration by $t_*$.
 We will exhibit it as the ``operational time" for our fractional diffusion process,
 and for distinction we will call the variable $t$ its ``physical time"
  In fact
 $f_{\alpha,\theta}(x,t_*)$ is a probability density in $x\in \RR$, evolving in  operational time $t_*>0$
 and $q_\beta (t_*,t)$ is a probability density in $t_*\ge 0$, evolving in physical time $t>0$.
%%%%%%%%%%%%%%%%%%%
\subsection{Stochastic interpretation.}
Clearly $f_{\alpha,\theta}(x,t_*)$ characterizes a stochastic process 
describing a trajectory $x=y(t_*)$ in the $(t_*,x)$ plane, 
that can be visualized as a particle travelling along space $x$, as operational time $t_*$ is proceeding.
Is there also a process $t_*=t_*(t)$, a particle moving along the positive $t_*$ axis,
 happening in physical time $t$? Naturally we want $t_*(t)$ increasing, at least in the weak sense,
$$  t_2> t_1 \Longrightarrow t_*(t_2) \ge t_*(t_1)\,.$$
We answer this question in the affirmative by showing that,  by inverting  the stable process
$t=t(t_*)$ whose probability density (in $t$, evolving in operational time $t_*$) is 
the extremely positively skewed
stable density
$$ r_\beta(t,t_*) = t_*^{-1/\beta} L_\beta^{-\beta} (t/t_*^{1/\beta})\,. \eqno (4.11)$$
In fact, recalling
$$ \widetilde r_\beta(s,t_*) = \exp (-t_* s^\beta)\,, \eqno(4.12)$$
there exists the  stable process $t=t(t_*)$, weakly increasing, with 
density in $t$ evolving in $t_*$ given by (4.11).
We call this process {\it the leading process}.
\vsp
Happily, we  can invert this process. Inversion of a weakly increasing 
trajectory means that horizontal segments are converted to vertical segments and viceversa 
jumps (as vertical segments) to horizontal segments (in graphical visualization).
\vsp
Consider a fixed sample trajectory $t=t(t_*)$ and its also fixed inversion $t_*=t_*(t)$.
Fix an instant $T$ of physical time and an instant $T_*$ of operational time . 
Then, because $t=t(t_*)$ is increasing, we have the equivalence
$$    t_*(T) \le T_* \Longleftrightarrow T \le t(T_*)\,,$$
which, with notation slightly changed by
$$ t_*(T)\to t'_*\,, \; T_* \to t_*\,, \; T\to t\,, \; t(T_*)\to t'\,,$$
implies 
$$ \int_0^{t_*} \!  q(t_*',t)\, dt_*' = \int_t^\infty \! r_\beta(t', t_*) \, dt' \,, \eqno(4.13)$$
for the probability density $q(t_*,t) $ in $t_*$ evolving in $t$.
It follows
$$q(t_*,t)= \frac{\partial}{\partial t_*}\int_t^{\infty}\! r_\beta(t',t_*)\, dt'
= \int_t^{\infty}\!\! \frac{\partial}{\partial t_*} r_\beta(t',t_*)\, dt'\,.$$
We continue in the $s_*$-Laplace domain assuming $t>0$,
$$ \widetilde q_(s_*,t)= \int_t^\infty \left(s_*\widetilde r_\beta(t',s_*)-\delta(t')\right)\, dt'\,.$$
It suffices to consider $t>0$, so that we have  $\delta(t')=0$ in this integral.
Observing from (4.12)
$$ \widetilde{\widetilde r}_\beta(s,s_*) = \frac{1}{s_*+s^\beta} \,, \eqno(4.14)$$
we find
$$   \widetilde r_\beta(t,s_*) = \beta t^{\beta-1} \,E'_\beta(-s_*t^\beta)\,, \eqno(4.15)$$
so that
$$  \widetilde q(s_*,t) = \int_t^\infty \! s_* \beta {t'}^{\beta-1}\, E'_\beta (-s_* {t'}^\beta)\, dt'
  = E_\beta(-s_*t^\beta)\,,\eqno(4.16)$$
  finally
  $$ q(t_*,t) = t^{-\beta}\, M_\beta(t_*/t^\beta)\,, \eqno(4.17)$$
From (4.16) we also see that
$$ \widetilde{\widetilde q}(s_*,s) =
\frac{s^{\beta-1}}{s_*+s^\beta} = \widetilde{\widetilde q}_\beta(s_*,s)\,, \eqno(4.18)$$
implying (4.6) and, see (4.7),
$$  q(t_*,t) \equiv  q_\beta(t_*,t)\,, \eqno(4.19)$$
so that indeed the process $t_*=t_*(t)$ is the inverse to the stable process $t=t(t_*)$ and has density
$q_\beta(t_*,t)$.

\noindent {\bf Remark}.
The process at hand, $t_*=t_*(t)$, which  is referred to as 
the inverse stable subordinator, is honoured with the name "Mittag-Leffler process" 
by Meerschaert at al. \cite{M3_PRE02sub_PRE02,M3_FPP}. 
%% Meerschaert, M.M., Benson, D.A., Scheffler, H.-P. and Baeumer, B.(2002): 
%% Stochastic solution of space-time fractional diffusion equation, Physical Review E 65, 041103 (R), 14.
Honouring this process by the name of Mittag-Leffler can  be justified by the fact that by (4.16)
the Laplace transform  of its density is a Mittag-Leffler type function or by the fact
 that it is 
a properly scaled diffusion limit of the counting function $N(t)$ of the fractional generalization of 
the Poisson process whose residual waiting time probability is the Mittag-Leffler type function  
$E_\beta(-t^\beta)$, see  recent papers of ours\cite{Gorenflo_PALA09,Gorenflo-Mainardi_HONNEF08}.
  In view of its probability density it may also be called the $M$-Wright process.   
\vsp
Stipulating that there exists  a weakly increasing process $t_*=t_*(t)$ 
with density  $q_\beta(t_*,t)$ we can analogously find the density of  its inverse $t=t(t_*)$
 which comes just as $r_\beta (t,t_*)$. 
 However, in the context of our here presented considerations  not being allowed 
 to know that such process $t_*=t_*(t)$  exists, 
 we have taken as a gift from God the process  $t=t(t_*)$ and shown by its inversion
  that there exists a process $t_*=t_*(t)$ with the desired properties.
%%%\noindent{\bf Conclusion}.
 \vsp
  From the density $r_\beta(t,t_*)$ of {\it the leading process}  $t =t(t_*)$
 we have found the density of {\it the  directing process} $t_*=t_*(t)$
as given by the  Laplace transform pair (4.6), that is
$$q_\beta(t_*,t) \div \widetilde q_\beta(t_*,s) = s^{\beta-1}\,\exp (-t_*s^\beta)\,.$$
\vsp
In physical coordinates we have (4.5) and (4.17), 
so also an expression through  an $M$-Wright function, 
$$   q_\beta(t_*,t) = \, _tJ^{1-\beta}\, r_\beta(t, t_*) = t^{-\beta}\, M_\beta(t_*/t^\beta)\,,
\eqno(4.20)$$
see Eq. (2.26).

%%%%%%%%%%%%
\subsection{Evolution equations for the densities $r_\beta(t,t_*)$ of $t=t(t_*)$ and
$q_\beta(t_*,t)$ of $t_*=t_*(t)$.}
%%%%%%%%%%%%%%%%%%
The Laplace-Laplace representation of the density $r_\beta(t,t_*)$ of the  process $t=t(t_*)$ 
is, according to (4.14),
$$ \widetilde{\widetilde r}_\beta(s,s_*) = \frac{1}{s^\beta +s_*}\,. $$
This implies 
$$s_*\, \widetilde{\widetilde r}_\beta(s,s_*)-1 =-s^\beta \,\widetilde{\widetilde r}_\beta(s,s_*)\,,$$
and by inverting the transforms and observing the initial condition $r_\beta(t,t_*=0)=\delta(t)$ 
we arrive at the Cauchy problem
$$ \frac{\partial}{\partial t_*} r_\beta (t, t_*) = - \, _tD_*^\beta r_\beta(t,t_*)\,,
 \quad r_\beta(t,t_*=0)= \delta (t)\,. \eqno(4.21)$$
Because it suffices to consider only $t>0$ where $\delta(t)=0$, we need not introduce a singular term on the 
right-hand side.
\vsp
The Laplace-Laplace representation of the density $q_\beta(t_*,t)$ of the  process $t_*=t_*(t)$ 
is, according to (4.18),
$$ \widetilde{\widetilde q}_\beta(s_*,s) = \frac{s^{\beta-1}}{s_* +s^\beta}\,. $$
This implies 
$$s^\beta\, \widetilde{\widetilde q}_\beta(s_*,s)-s^{\beta-1} =-s_* \,\widetilde{\widetilde q}_\beta(s_*,s)\,,$$
and by inverting the transforms and observing the initial condition 
$q_\beta(t_*,t=0)=\delta(t_*)$ we arrive at the Cauchy problem
$$ _tD_*^\beta q_\beta(t_*,t)=
- \frac{\partial}{\partial t_*} q_\beta (t_*, t)  
 \quad q_\beta(t_*,t=0)= \delta (t_*)\,. \eqno(4.22)$$
Because it suffices to consider only $t_*>0$ where $\delta(t_*)=0$, we can ignore the delta function on the 
right-hand side.
\vsp
{\bf Remark} The fractional differential equations in the above Cauchy problems have the same form.
By replacing $t$ by $t_*$ and $r$ by $q$
 one of them goes over into the other.
 However, in the first problem the delta initial condition refers to the fractional derivative (of order $\beta$), 
 in the second problem to the ordinary (first order) derivative.
 These equations are akin with the time-fractional drift equation treated in (3.14) and (3.15),
 with different coordinates and proper  initial conditions, as explained above. 
 The process $t=t(t_*)$ of the first problem is 
 a positive-oriented (extreme) stable process, whereas the process $t_*=t_*(t)$
 is a fractional drift process, see (3.14)-(3.15) with $x$ replaced by $t_*$.
 The reason for the two evolution equations to have the same form is that the described two processes 
 are inverse to each other, their graphical representations coincide just by interchanging the 
 coordinate axes. The delta initial  condition for each equation is given at value zero 
 of the evolution variable for the variable in which the solution is a density. 

\subsection{The random walks.}
We can now construct the process $x=x(t)$ for the position $x$ of the particle depending on physical time $t$
as follows in two ways. With the variable $t$ (physical time), $t_*$ (operational time), $x$ (position),
we have the processes (i), (ii) and (iii), as follows:
\par\noindent
(i) $t=t(t_*)$ with density $r_\beta(t, t_*)$ in $t$, evolving in $t_*$, {\it the leading process},
\par\noindent
(ii) $ x=y(t_*)$ with density $f_{\alpha,\theta} (x,t_*)$ in $x$, evolving in $t_*$, {\it the parent process},
\par\noindent
(iii) $ t_*=t_*(t)$ with density $q_\beta(t_*,t)$ in $t_*$, evolving in $t$, {\it the directing process}.
\par\noindent
Observing that  the processes (i) and (iii) are inverse to each other, 
and taking account of the subordination integral (4.10), 
we define {\it the space-time fractional diffusion process}
as  {\it the  subordinated process}
$$ x =x(t) = y(t_*(t))\,. \eqno(4.23)$$
Simulation of a trajectory for the subordinated process means: generate in running physical time $t$ 
the operational time $t_*$, then the operational process $y(t_*)$.
\vsp
Now, the Mittag-Leffler (or $M$-Wright) process $t_*=t_*(t)$ is non-Markovian and not so easy to simulate. 
The alternative 
(we call it "parametric subordination") is to produce in dependence of the operational time $t_*$ 
the processes (i) and (ii) and then eliminate $t_*$ from the system
$$ t=t(t_*)\,, \quad  x=y(t_*)\,, \eqno (4.24)$$
to get $x=x(t)$  from $x=y(t_*)$ by change of time from $t_*$ to $t$.
\vsp
We can produce a sequence of precise snapshots of $t=t(t_*)$
 and $x=y(t_*)$ in the $(t_*,t)$ plane and the $(t_*,x)$ plane by setting, with a step-size $\tau_*>0$,
 $$t_{*,n} = n \tau_*\,, \; \bar t_n= T_{*,1} + T_{*,2} \cdots + T_{*,n}\,, \; 
 \bar x_n= X_{*,1} +X_{*,2} + \cdots + X_{*,n}\,, \eqno(4.25)$$
 taking for $k=1,2, \dots,  n$ each  $T_{*,k}$  as a random number with density
 $ \tau_*^{-1/\beta} \, L_\beta^{-\beta} (t/\tau_*^{1/\beta})$ 
 and each $X_{*,k}$ as a random number   with density 
 $ \tau_*^{-1/\alpha} \, L_\alpha^{\theta} (x/\tau_*^{1/\alpha})\,,$
corresponding by self-similarity to the step $\tau_*$. 
\vsp
We can do this by taking random numbers $T_{k}$ and $X_{k}$ with density
$L_\beta^{-\beta}(t)$ and $L_\alpha^\theta(x)$, respectively, and then
with $\tau=\tau_*^{1/\beta}$ and $h=\tau_*^{1/\alpha}$, setting
$T_{*,k}= \tau\,T_k$, $X_{*,k}= h\,X_k$.
In other words: we produce (a renewal process at equidistant times with reward)
a positively oriented random walk on the half-line $t\ge 0$ and a random walk on
$-\infty<x<+\infty$ with jumps at equidistant operational time instants $t_*=n\tau_*$.
We recognize the scaling relation $\tau^\beta/h^\alpha \equiv 1$, 
analogous to that used by us in earlier papers of ours on well-scaled passage 
to the diffusion limit in $CTRW$ under power law regime, see
\cite{Gorenflo-Abdel-Rehim_VIETNAM04,Gorenflo-Mainardi_HONNEF08,Gorenflo-Mainardi_JCAM09}. 
%% \vsp
Methods for producing  stable random deviates can be found in the books
 \cite{Janicki_LN96,Janicki-Weron_94}
\vsp
Finally, we transfer into the $(t,x)$ plane the points with coordinates
$\bar t_n, \bar x_n=\bar y_n$ and so obtain  a sequence of precise snapshots 
of a true process $x=x(t)$. 
Finer details of the process $x=x(t)$ become visible  by using  smaller values of 
the operational step-size $\tau_*$.
%%%%%%%%%
%%%%%%%%%
  
\section{Graphical representations and Conclusions}
\begin{figure}[!h]
% \vspace{-0.5truecm}
\begin{center}
\includegraphics[width=.60\textwidth]{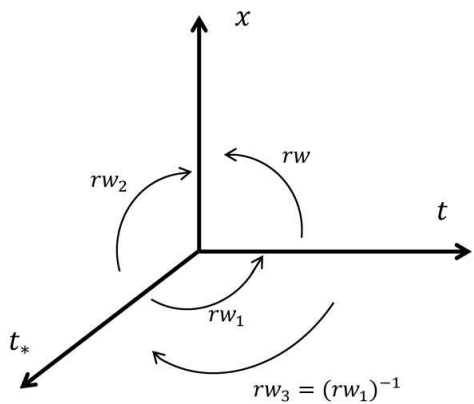}
\end{center}
\vspace{-0.8truecm}
 \caption{Diagram for the connections between the four random walks ($rw_1$),
($rw_2$), ($rw_3$) and ($rw$), related to the leading, parent, directing and subordinated processes, respectively.}
\label{fig:DIAGRAM}
\end{figure}

\begin{figure}
\includegraphics[width=.48\textwidth]{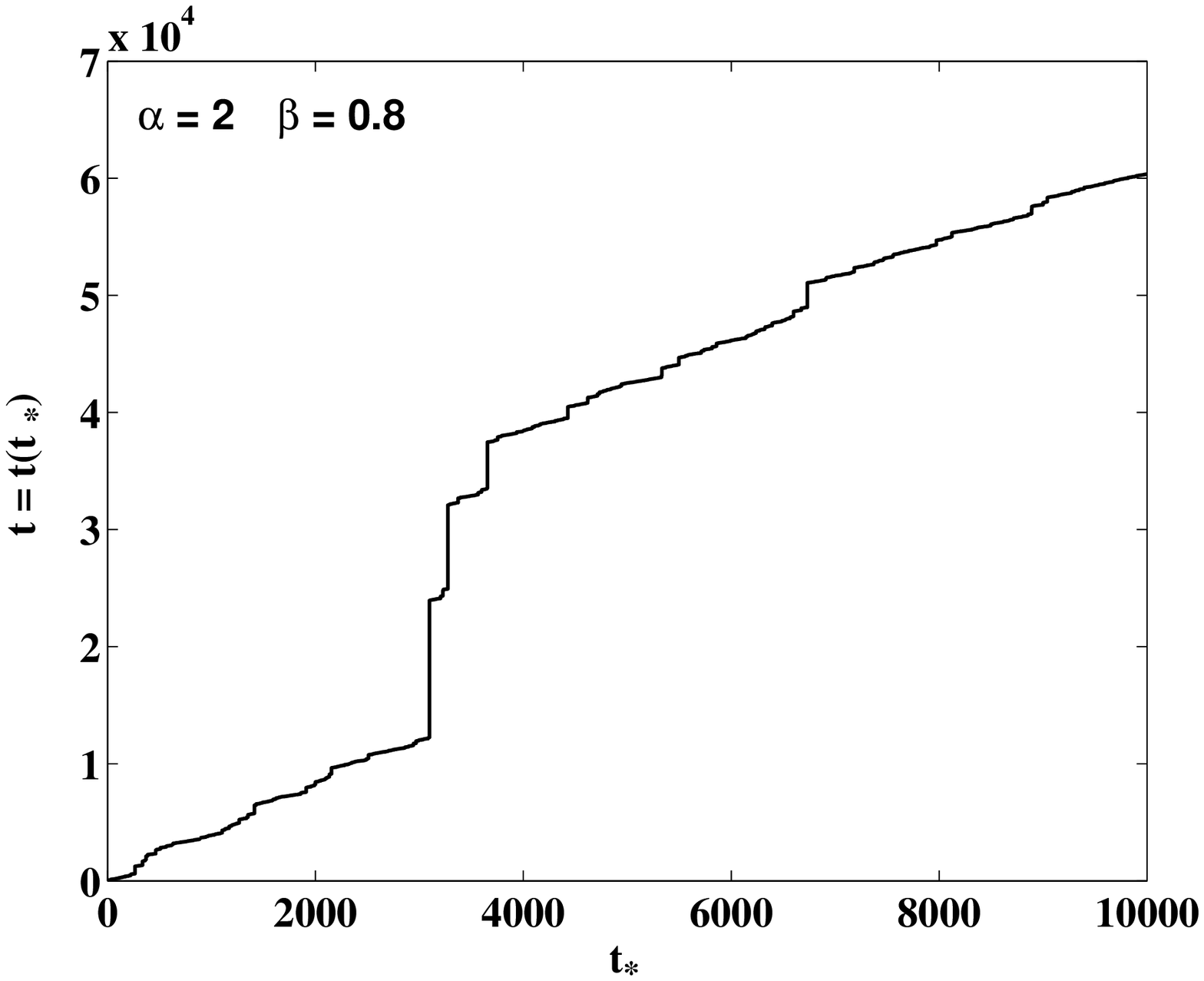}
 \includegraphics[width=.48\textwidth]{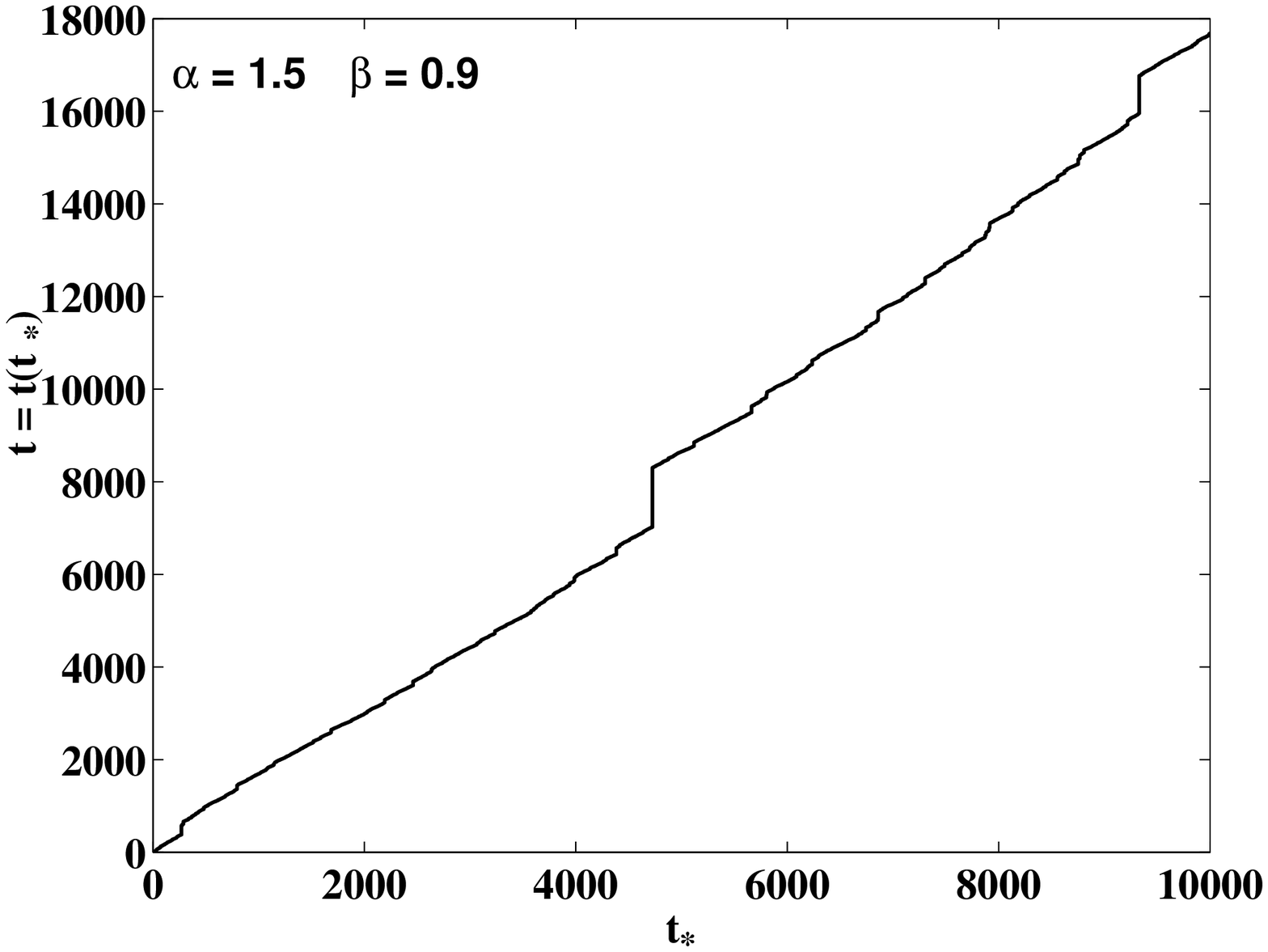}
 \vspace{-0.5truecm}
 \caption{A sample path for ($rw_1$), the leading process $t=t(t_*)$.}
 \centerline{LEFT: $\{\alpha =2\,,\; \beta =0.80 \}$,
      RIGHT: $\{\alpha =1.5\,,\; \beta =0.90 \}$.}
	  \label{fig:rw1}
%%	  \end{figure}
 \vskip 0.50truecm
% \begin{figure}
%% \includegraphics[width=.48\textwidth]{av_a20b80_c.eps}
%% \includegraphics[width=.48\textwidth]{av_a15b90_c.eps}
 \includegraphics[width=.48\textwidth]{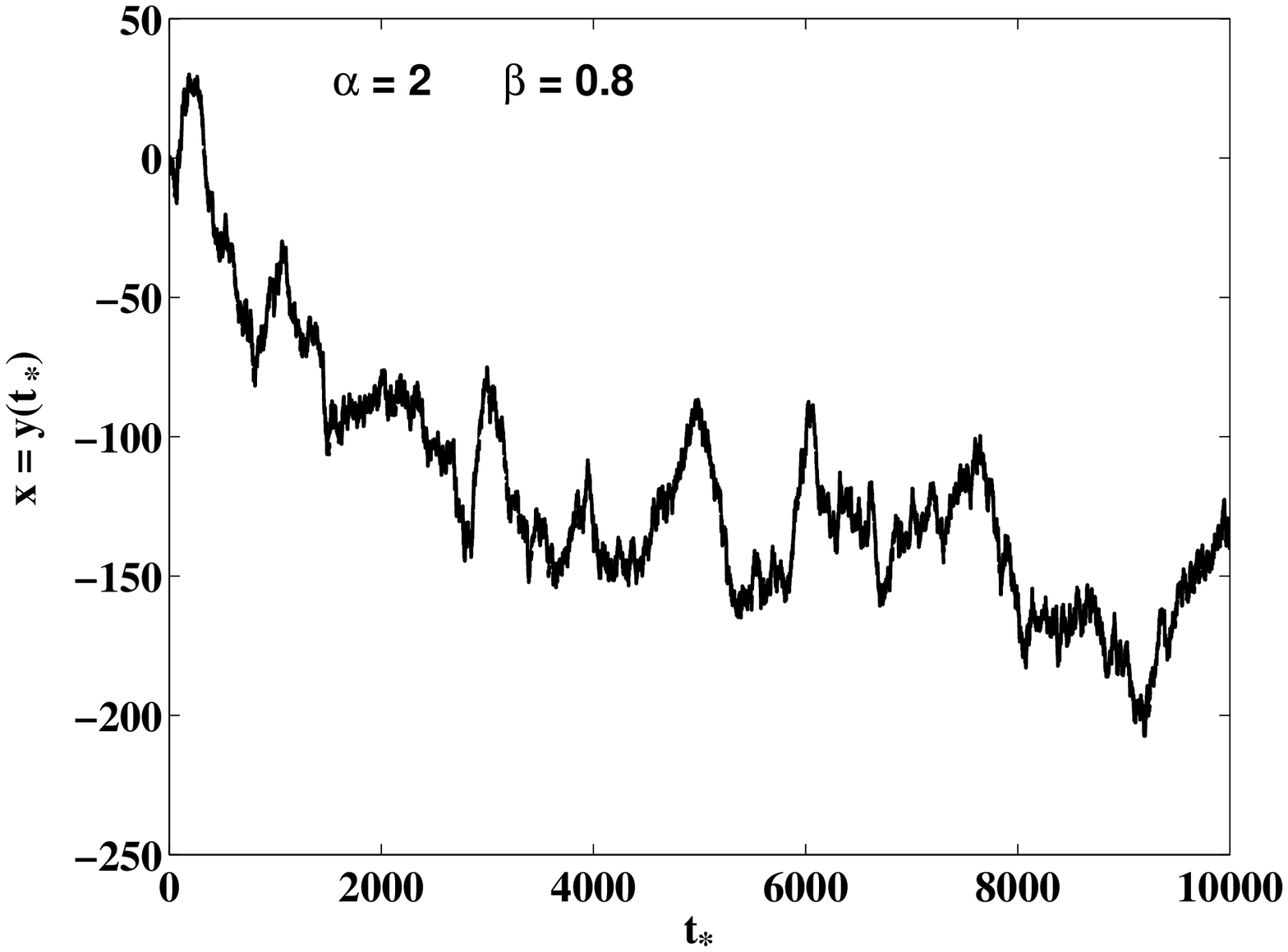}
 \includegraphics[width=.48\textwidth]{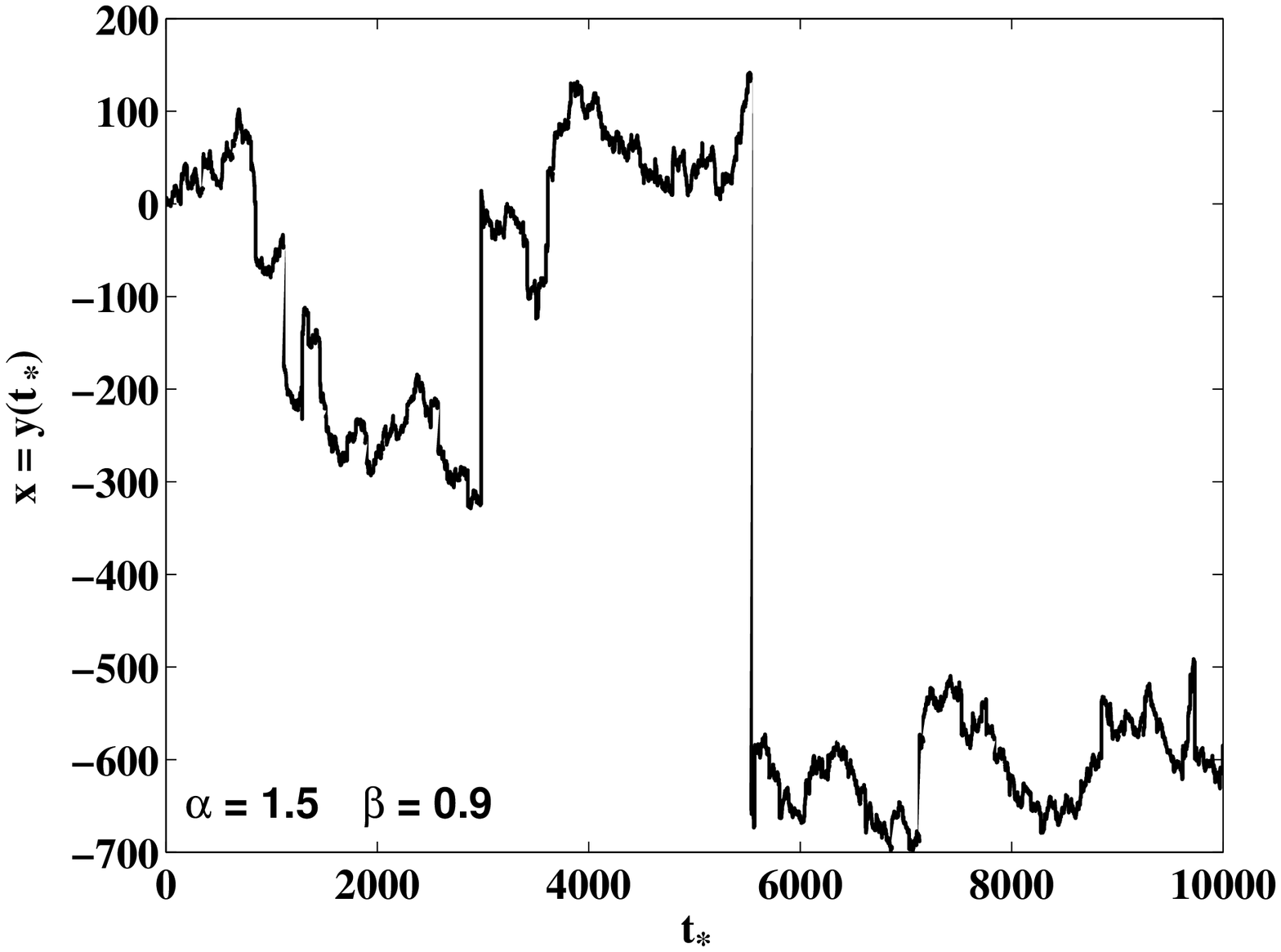}
\vspace{-0.5truecm} 
 \caption{A sample path for ($rw_2$), the parent process $x=y(t_*)$.}
 \centerline{LEFT: $\{\alpha =2\,,\; \beta =0.80 \}$,
      RIGHT: $\{\alpha =1.5\,,\; \beta =0.90 \}$.}
	  \label{fig:rw2}
%%\end{figure}
%%\begin{figure}
 \vskip 0.50truecm
\includegraphics[width=.48\textwidth]{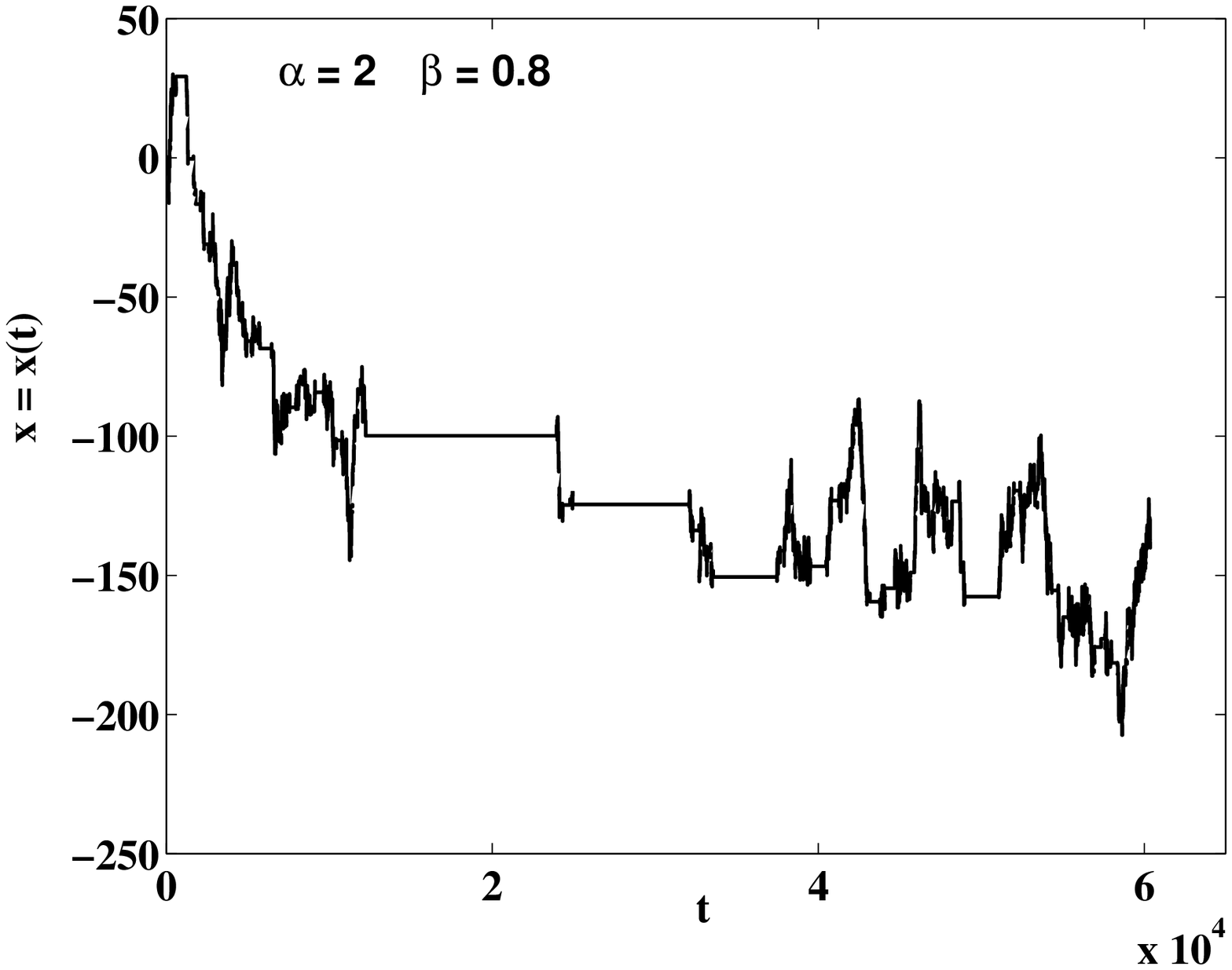}
 \includegraphics[width=.48\textwidth]{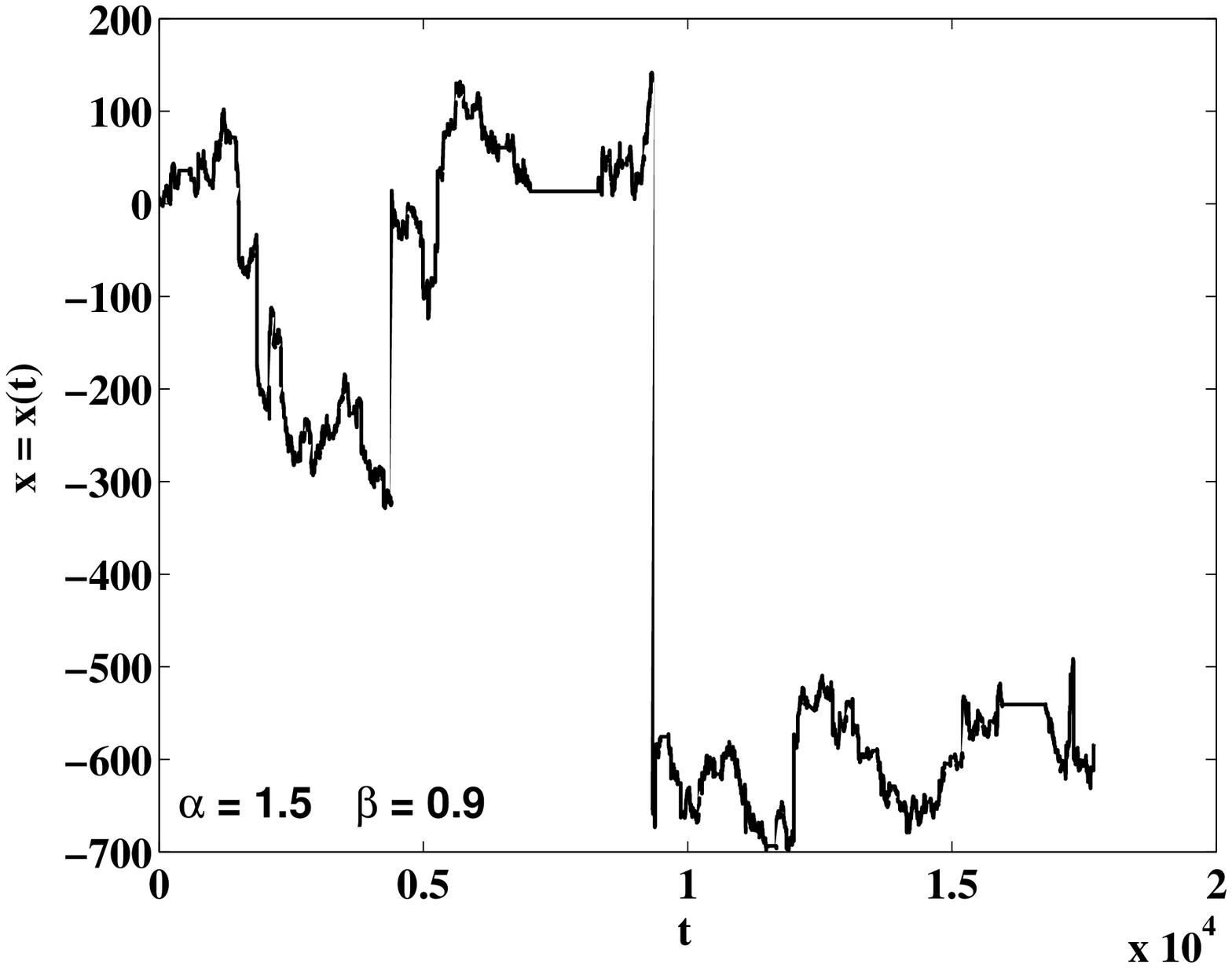}
\vspace{-0.5truecm} 
 \caption{A sample path for ($rw$), the subordinated process $x=x(t)$.}
 \centerline{LEFT: $\{\alpha =2\,,\; \beta =0.80 \}$,
      RIGHT: $\{\alpha =1.5\,,\; \beta =0.90 \}$.}
	  \label{fig:rw}
\end{figure}

We  recall that,  denoting   the physical time with $t$,
  the operational time with $t_*$, the physical space with $x$
%% we have distinguished the following process:
%% $t=t(t_*) \;$ {\it the leading process} , 
%%%% governed by the probability density $r_\beta(t,t_*)$ in $t$ evolving in $t_*$
%% \\ 
%% $t_*=t_*(t) \;$ {\it directing process},
%%%%  governed by the probability density $q_\beta(t_*,t)$ in $t_*$ evolving in $t$
%% \\
%% $x=y(t_*)= x(t(t_*))  \;$ {\it the parent process}, 
%% governed by the probability density  $f_{\alpha,\theta}(x,t_*)$ in $x$ evolving in $t_*$, 
%%\\
%% and finally 
%%$x=x(t) = y(t_*(t)) \;$ {\it the fractional diffusion process}.  
%%\vsp  the probability density
the density of the fractional diffusion process %% , $u(x,t)$ in $x$ evolving in $t$,  
turns out to be  %% related to those of the parent and directing process 
given by the following subordination integral, see (4.10),   
$$ u(x,t) = \int_{t_*=0}^\infty \!\! f_{\alpha,\theta}(x,t_*)\, q_\beta(t_*,t)\, dt_*\,, \eqno(5.1)$$
where $f_{\alpha,\theta}(x,t_*)$, is the  density (in $x$ evolving in $t_*$) of the parent process
$x=y(t_*)= x(t(t_*))$
and $q_\beta(t_*,t)$ is the density (in $t_*$ evolving in $t$) of the directing process
$t_*=t_*(t)$.
\vsp 
By using the Fourier-Laplace pathway  we  recall the two densities related to the parameters
$\alpha, \theta, \beta$ from  (4.3) -- (4.5),
$$ f_{\alpha, \theta}(x,t_*) = t_*^{-1/\alpha}\, L_\alpha^\theta(x/t^{1/\alpha})\,,\eqno(5.2)$$
$$ q_\beta(t_*,t) = t_*^{-1/\beta}\,_tJ^{1- \beta}\, L_\beta^{-\beta}\left( t/t_*^{1/\beta}\right)=
t^{-\beta}\, M_\beta(t_*/t^\beta)\,,\eqno(5.3)$$
where $L$ refers to the L\'evy stable density and $M$ to the Wright function,
both introduced in Section 2.
But for the  {\it parametric subordination} the relevant density is
$r_\beta(t,t_*)$ governing the {\it leading process}, a density in the physical time $t$
evolving with the operational time $t_*$: it turns out to be the
unilateral L\'evy density of order $\beta$,  namely, see (4.11),
$$ r_\beta(t ,t_*) = t_*^{-1/\beta}\, L_\beta^{-\beta}\left( t/t_*^{1/\beta}\right)\,.\eqno(5.4)$$
We have shown in Section 4.4  that in our approach (referred to as {\it parametric subordination})
 the  process of space-time fractional   diffusion (non-Markovian for $\beta<1$)
 can be  simulated by two Markovian processes governed by stable densities,
 provided by $f_{\alpha, \theta}(x,t_*)$ and $r_\beta(t ,t_*)$,
 as pointed out in our 2007 paper with Vivoli\cite{GorMaiViv_CSF07},
  where we have dealt with the $CTRW$ model. 
%% by Gorenflo, Mainardi and Vivoli (2007).
There, before passing to the diffusion limit, we have two Markov processes 
happening on a discrete set of equidistant
instants (for simplicity the non-negative integers)
$n=0, 1, 2, \dots$, meaning $\tau_*=1$, one of them moving randomly rightwards along  $t \ge 0$,
the other moving randomly on the real line  $-\infty < x < \infty$ with
jumps $T_{n}$ and $X_{n}$, respectively, for $n \ge 1$. By summing the
"waiting times" $T_{k}$ and the jumps $X_{k}$ from 1 to $n$ we obtain
sequences of jump instants $t=t_n$  and positions $x=x_n$,that we display in the 
($t,x$) plane. In fact, CTRW is the virgin form of parametric subordination. 
\vsp
We note that our approach is akin to that based on two stochastic differential equations,
known in Physics as Langevin equations, see%%
\cite{Fogedby_PRE94,Kleinhans-Friedrich_PRE07,Zhang-Meerschaert-Baeumer_PRE08}.
 %% Fogedby (1994), Kleinhans and Friedrich (2007), Zhang, Meerschaert and Baeumer (2008).   
We have indicated these two stochastic differential equations in\cite{GorMaiViv_CSF07}. 
Here now we content ourselves with referring to the above cited papers.
\vsp
In Section 4.4 we have splitted the fractional diffusion process into three
processes (i), (ii), (iii), each of them containing two of the three coordinates: space
$x$, physical time $t$, operational time $t_*$. 
We simulate {\it the leading process} by a random
walk $(rw_1)$, {\it the parent process} by a random
walk $(rw_2)$, and {\it the subordinated process}
(which yields the desired trajectory) by a
random walk $(rw)$. The inversion of $(rw_1)$
gives us a random walk $(rw_3)$ for simulation of
{\it the directing process} $t_* = t_*(t)$. 
Essentially, we need to carry out only ($rw_1$) and ($rw_2$) according to the equations
(4.25). By transferring the points ($\bar t_n, \bar x_n$) into the ($t, x$) plane
we get the random walk ($rw$) 
as visualization of a random trajectory $x = x(t) = y(t_*(t))$ 
according to the subordinated process which is our space-time fractional
diffusion process of interest. 
\vsp
To make transparent the situation we display as a
diagram in Fig. \ref{fig:DIAGRAM} the connections between the four random walks.
%%%%%%%%%%%
%%%%%%%%%%%%%
%% \vsp
It is now instructive to
show some numerical realizations of these random walks
for two case studies %%  of the required fractional diffusion process
of symmetric ($\theta=0$) fractional diffusion processes:
%% $\{\alpha =2,  \, \beta =0.90\}$,
$\{\alpha =2,  \, \beta =0.80\}$,
$\{\alpha =1.5, \,  \beta =0.90\}$.
%% $\{\alpha =1.5, \,  \beta =0.80\}$.
%% \vsp
As explained in a previous sub-section, for each case
we need to construct the sample paths for three  distinct processes,
the leading process $t = t(t_*)$,
the parent process $x= y(t_*)$ (both in the operational time)
and, finally, the subordinated process $x =x(t)$,
corresponding to the required fractional diffusion process.
%%%%%%%%%%%%%%%%%%% INSERTION  \newpage
%%%%%%%%%
%% \vsp
%% The essential walks to be carried out are $(rw_1)$
%% and $(rw_2)$; for their construction see the formulas in (4.25) for $\bar{t}_n$ and
%% $\bar{x}_n$. To obtain at first a rough picture (but already with precise points of
%% the true processes) one can take $\tau_* = 1$ for simplicity.
\vsp
We shall depict the above sample paths in 
Figs. \ref{fig:rw1}, \ref{fig:rw2}, \ref{fig:rw}
%% Fig. 1, Fig. 2 and Fig. 3,
respectively, devoting the left  and the right plates to
the different case studies.
%%%%%%%%%%%%%
\vsp
Plots in Fig.\ref{fig:rw1} (devoted to the leading process, the limit of ($rw_1$))
thus represent sample paths in the $(t_*,t)$  plane
of unilateral L\'evy motions  of order $\beta$. 
By interchanging the coordinate axes we can consider Fig.\ref{fig:rw1}  as representing sample paths of 
the directing process, the limit of ($rw_3$). 
 \vsp
Plots in Fig.\ref{fig:rw2}  (devoted to the parent process, the limit of ($rw_2$))
represent   sample
paths in the $(t_*,x)$ plane,
produced in the way explained above, for  L\'evy motions
of  order  $\alpha$ and skewness  $\theta=0$
(symmetric stable distributions).
%% keeping fixed parameter $\theta=0.$
\vsp
By the indicated method, see (4.25), we have with (for simplicity) $\tau_*=1$,
$\theta=0$ (symmetry) produced 10000 numbers $\bar{t}_n$ and corresponding numbers
$\bar{y}_n$. Plotting the points ($t_{*,n}, \bar{t}_n$) into the ($t_*, t$)
(operational time, physical time) plane, the points ($t_{*,n}, \bar{y}_n$) into the
($t_*, x$) (operational time, position in space)  plane we get Fig.\ref{fig:rw1} 
 and Fig.\ref{fig:rw2} for visualization of ($rw_1$) and ($rw_2$), respectively. 
 Fig. \ref{fig:rw} (as a visualization of ($rw$) is obtained by
plotting the points ($\bar{t}_n,\bar{y}_n$) into the ($t,x$) (physical time and
space) plane.
\vsp
Actually, we have invested a little bit more work in producing the figures. Namely,
to make visible the jumps as vertical segments, we have in Fig.\ref{fig:rw1} connected the
points ($(t_{*,n}, \bar{t}_n$) and ($(t_{*,n+1}, \bar{t}_n$) by a horizontal
segment, the points ($(t_{*,n+1}, \bar{t}_n$) and ($(t_{*,n+1}, \bar{t}_{n+1}$) by a
vertical segment. Analogously in Fig.\ref{fig:rw2} with the indexed $\bar{t}$ replaced by indexed
$\bar{y}$. In Fig.\ref{fig:rw} we have connected the points ($\bar{t}_n,\bar{y}_n$) and
($\bar{t}_{n+1},\bar{y}_n$) by a horizontal segment, 
the points ($\bar{t}_{n+1}, \bar{y}_n$) and ($\bar{t}_{n+1},\bar{y}_{n+1}$) by a vertical segment.
\vsp
Resuming,
we can consider Fig.\ref{fig:rw1} as a representation of ($rw_1$) for the leading process,
or by interchange of axes as one of ($rw_3$) for the directing process,
Fig.\ref{fig:rw2} as one of ($rw_2$) for the parent process, and finally 
Fig.\ref{fig:rw} as a representation of ($rw$) for the subordinated process which is 
our space-time fractional diffusion process.
%%%%%%%%%
\newpage
%%%%%%%%
%%%%%%%%%%%%  REFINEMENT FOR THE 2 LEADING PROCESSES  %%%%%%%%%%
\begin{figure}[!h]
\includegraphics[width=.49\textwidth]{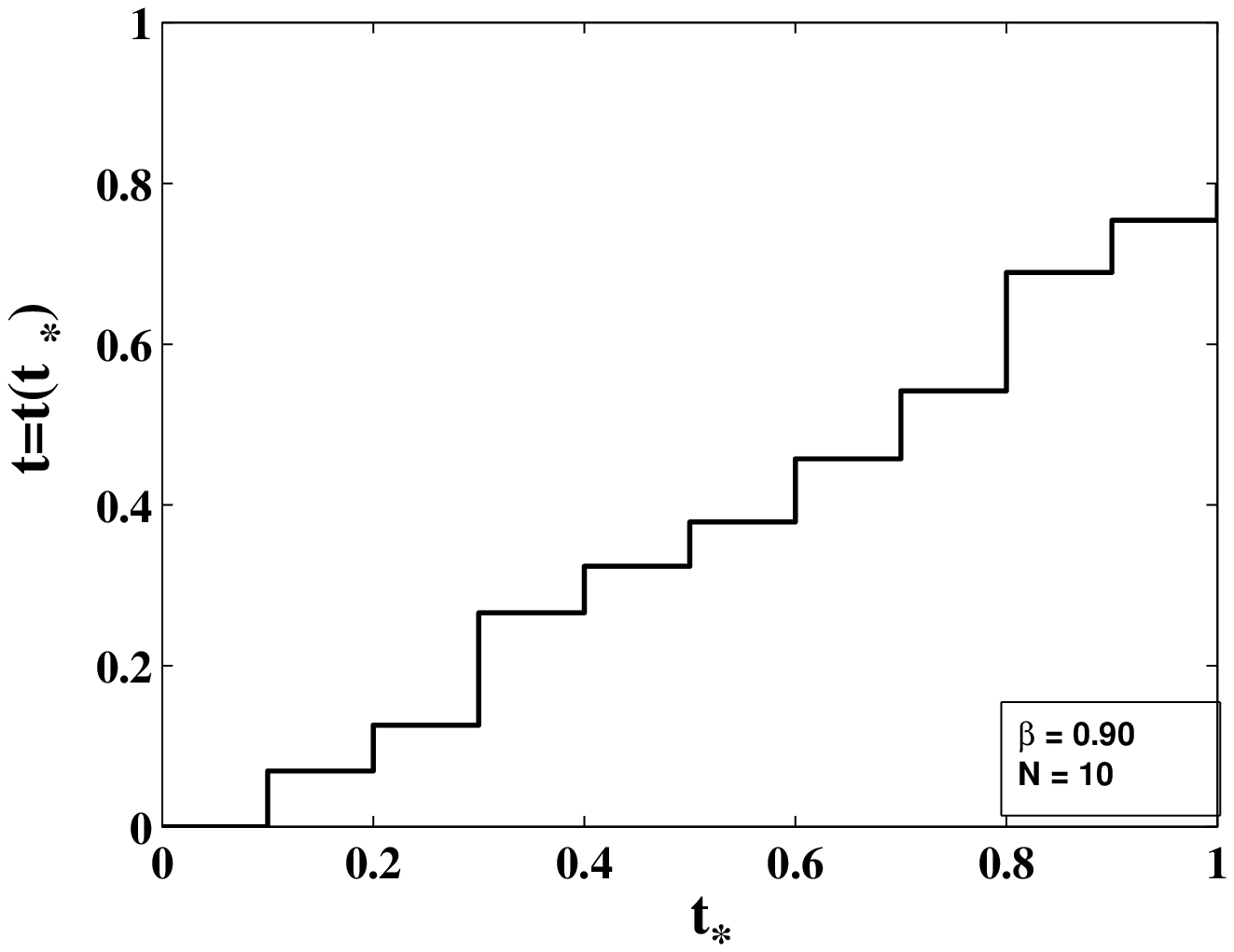} 
\includegraphics[width=.49\textwidth]{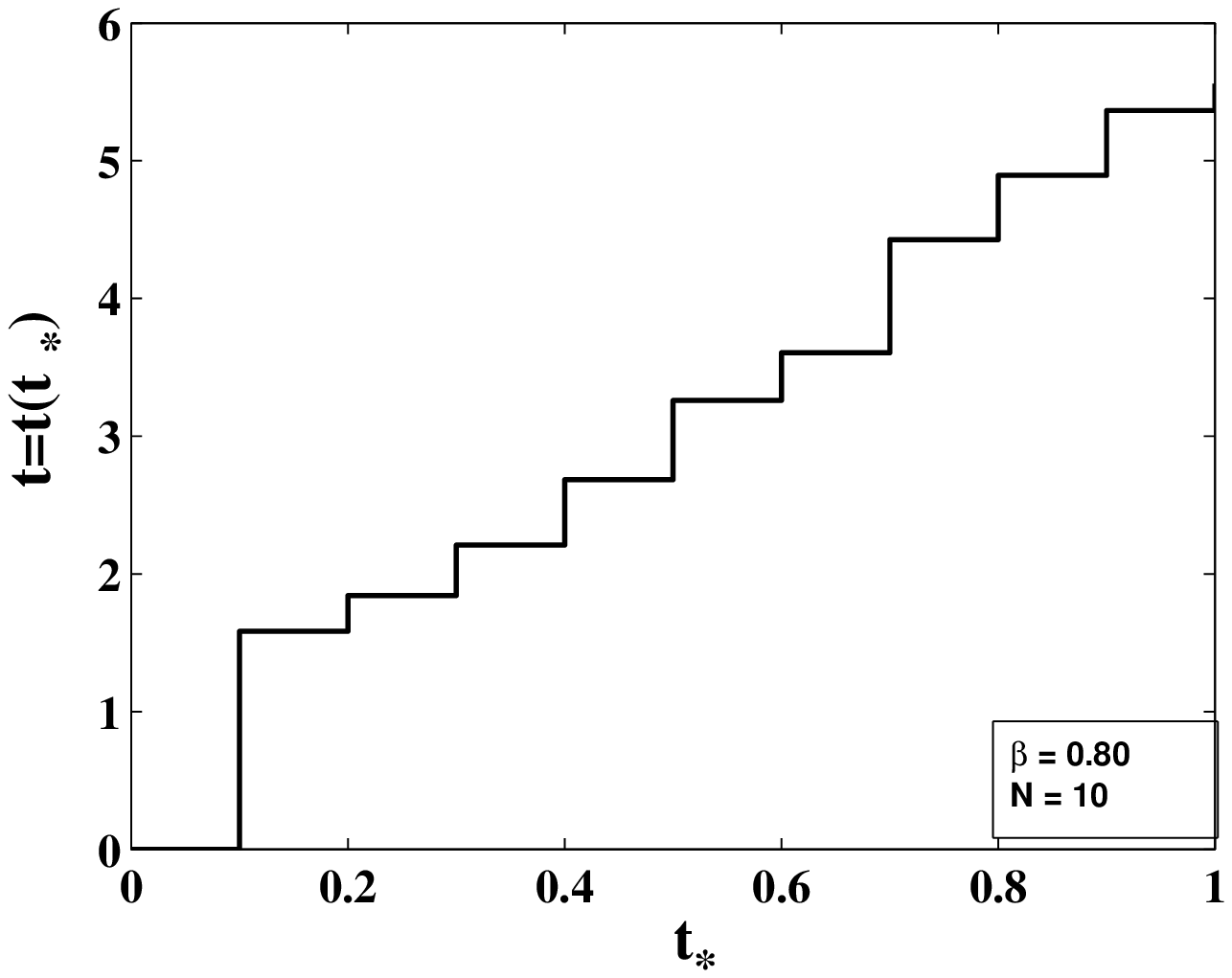}
\vspace{-0.5truecm}
 \caption{A sample path for the leading process $t=t(t_*)$.}
 \centerline{LEFT: $\{\beta  =0.9,\; N = 10^1 \}$, RIGHT: $\{\beta  =0.8,\; N= 10^1 \}$.}
 \label{fig:leading10}
\vskip 0.20truecm
 \includegraphics[width=.49\textwidth]{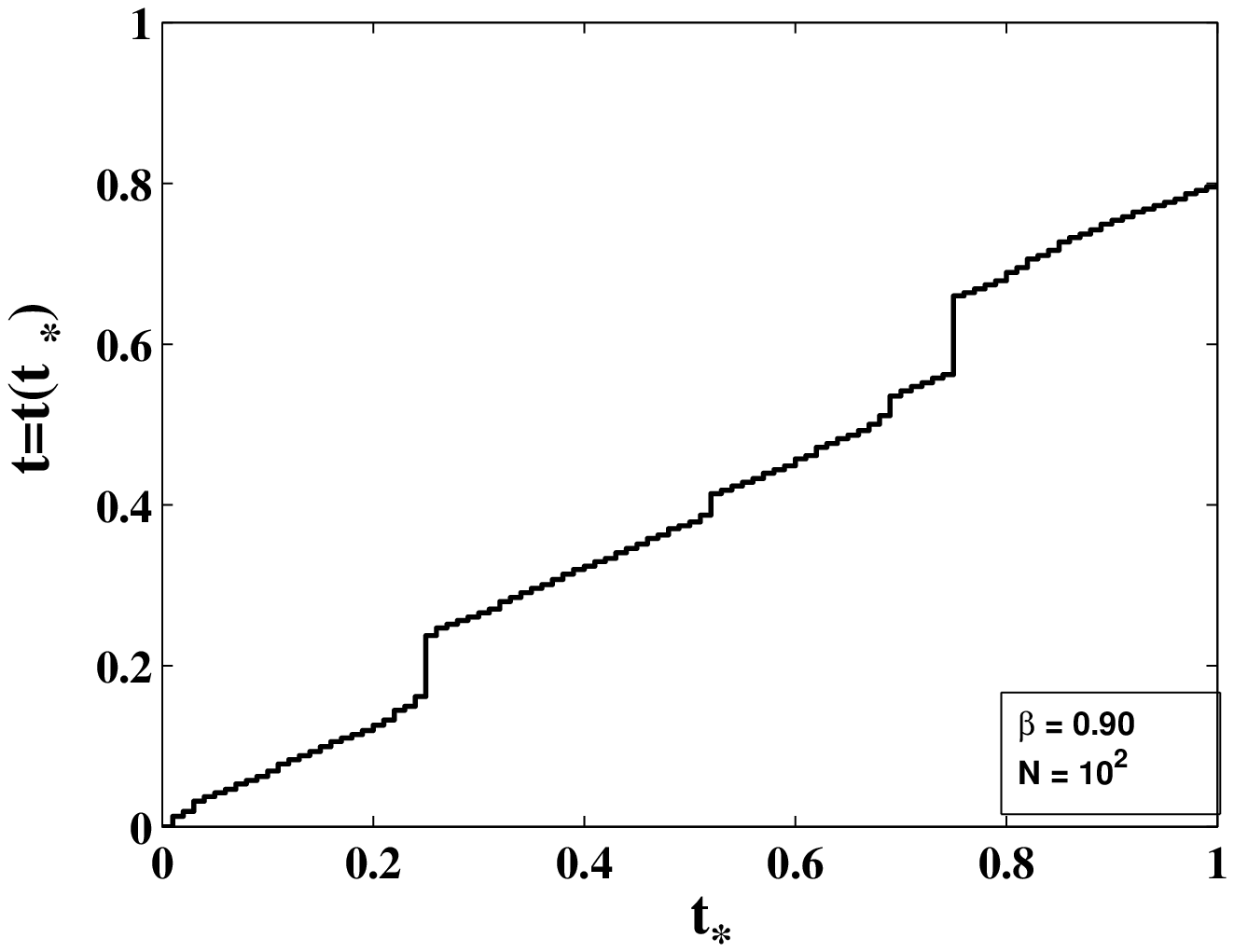} 
 \includegraphics[width=.49\textwidth]{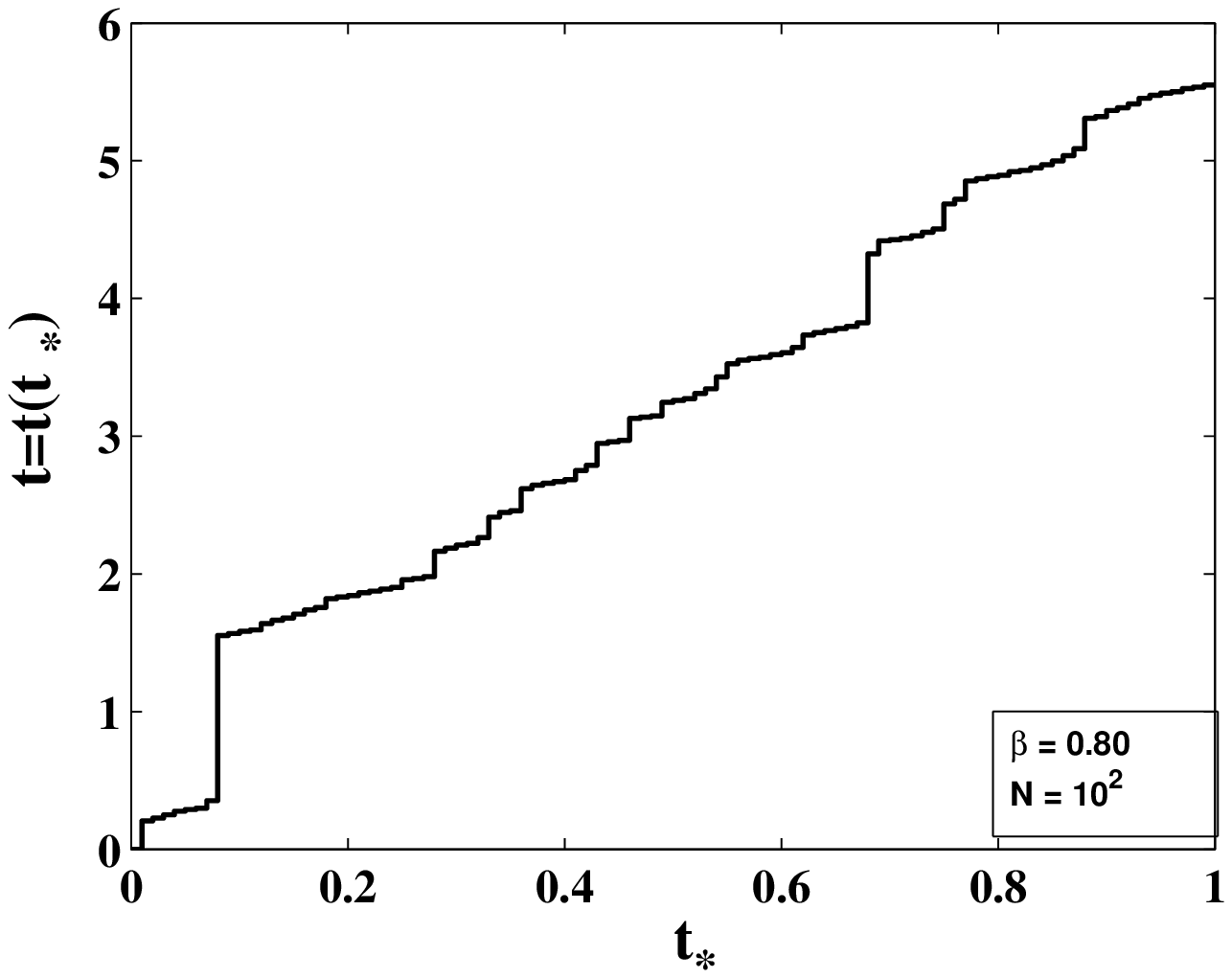}
 \vspace{-0.5truecm}
 \caption{A sample path for the leading process $t=t(t_*)$.}
 \centerline{LEFT: $\{\beta  =0.9,\; N = 10^2 \}$, RIGHT: $\{\beta  =0.8,\; N= 10^2 \}$.}
 \label{fig:leading100}
\vskip 0.20truecm
 \includegraphics[width=.49\textwidth]{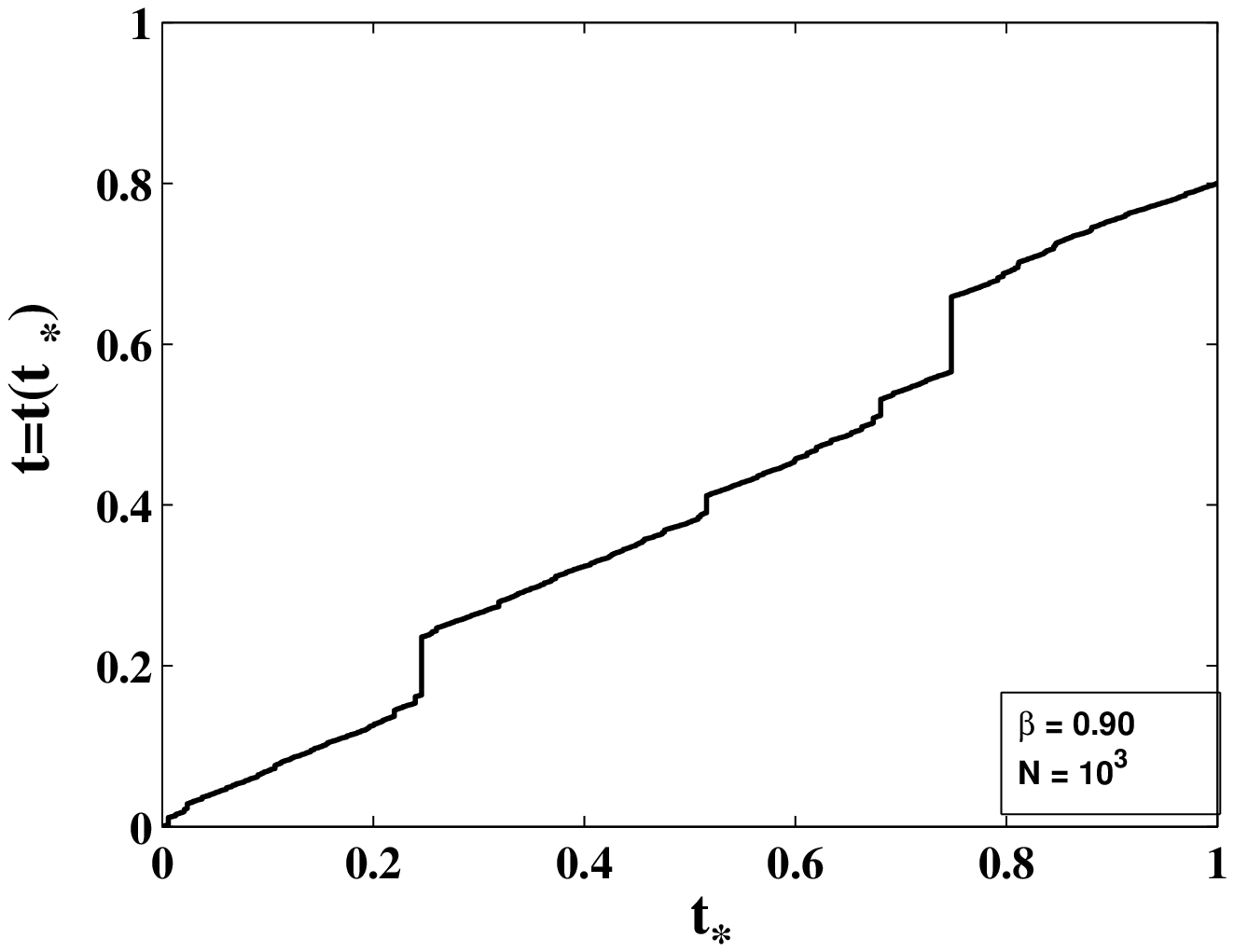} 
 \includegraphics[width=.49\textwidth]{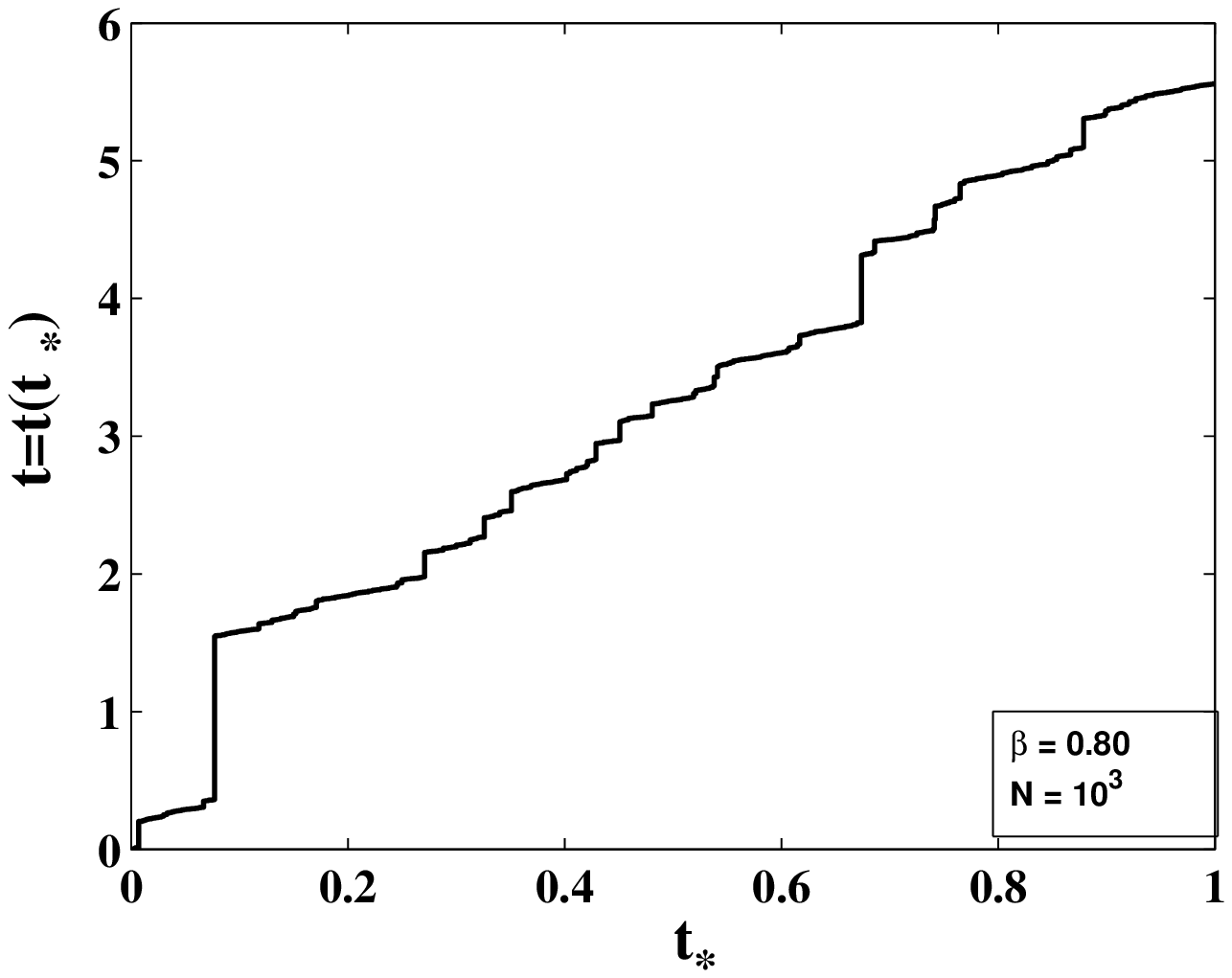}
\vspace{-0.5truecm} 
 \caption{A sample path for the leading process $t=t(t_*)$.}
 \centerline{LEFT: $\{\beta =0.9,\; N=10^3 \}$,
      RIGHT: $\{\beta =0.8,\; N=10^3 \}$.}
	  \label{fig:leading1000}
\end{figure}
%%%%%%%%%%%%  REFINEMENT FOR THE 2 PARENT PROCESSES  %%%%%%%%%%
\newpage
\begin{figure}[!h]
 \includegraphics[width=.49\textwidth]{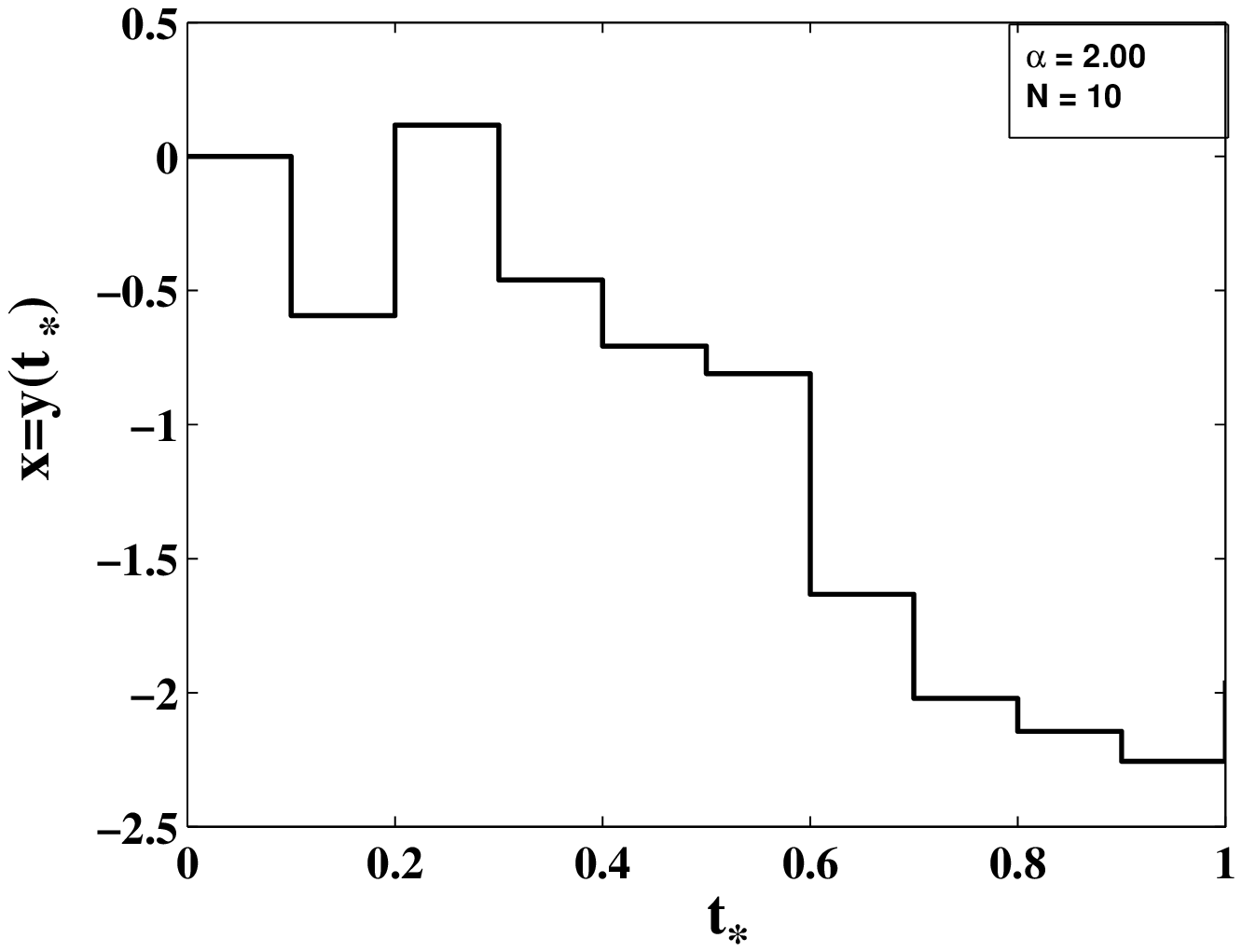}
\includegraphics[width=.49\textwidth]{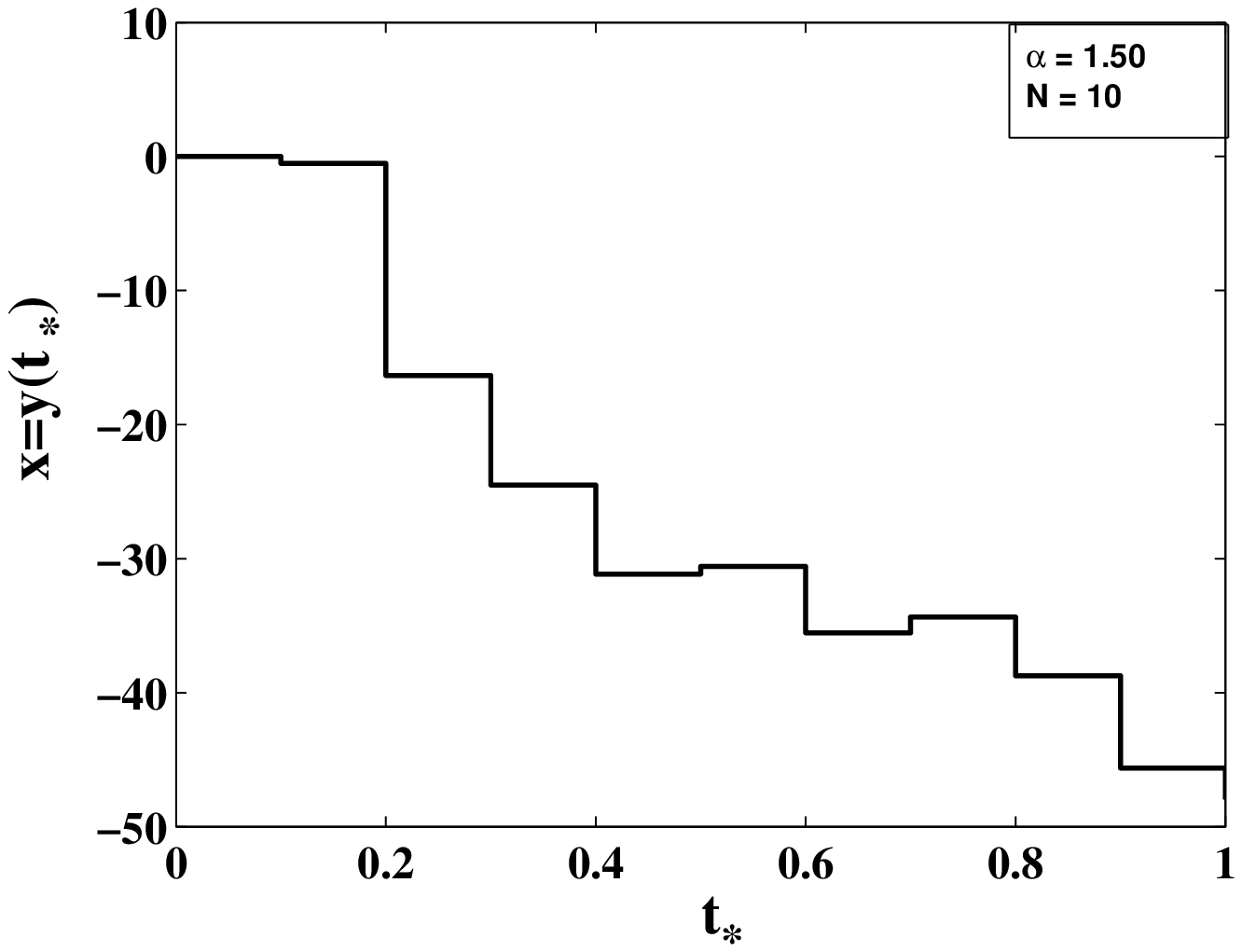}
\vspace{-0.5truecm} 
 \caption{A sample path for the parent process $x=y(t_*)$.}
 \centerline{LEFT: $\{\alpha =2,\; N = 10^1 \}$, RIGHT: $\{\alpha =1.5,\; N= 10^1 \}$.}
 \label{fig:parent10}
\vskip 0.20truecm
 \includegraphics[width=.49\textwidth]{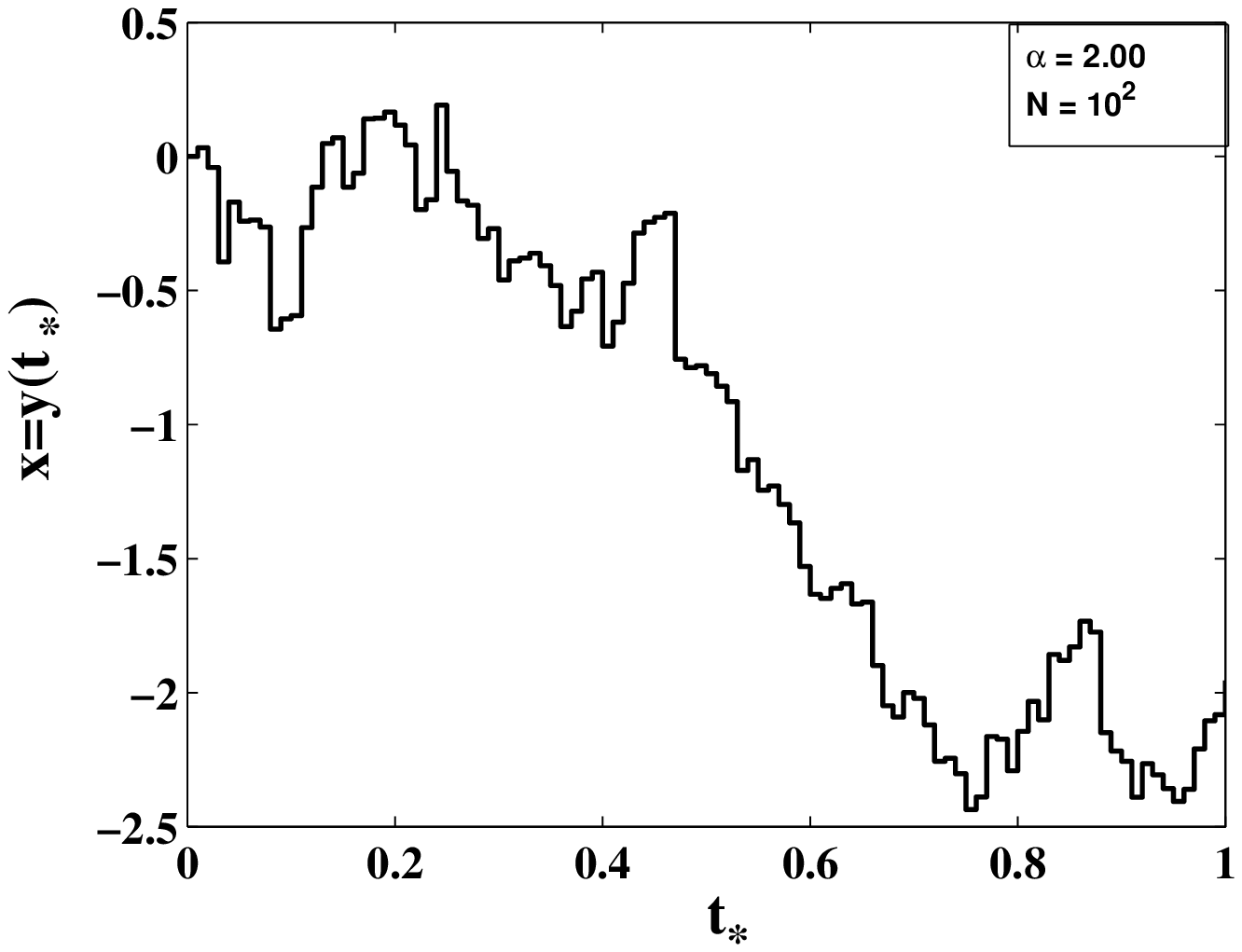}
\includegraphics[width=.49\textwidth]{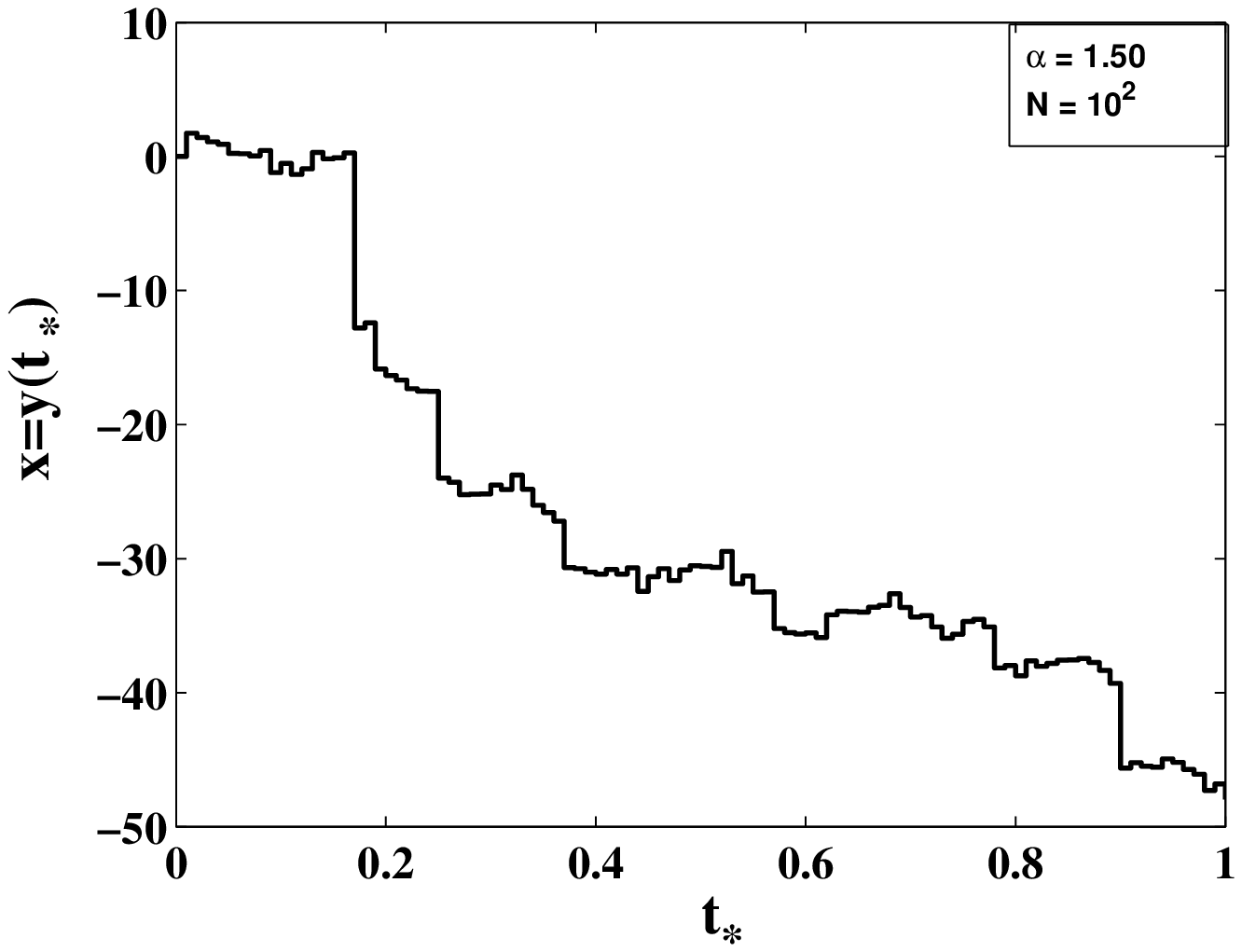}
\vspace{-0.5truecm} 
 \caption{A sample path for the parent process $x=y(t_*)$.}
 \centerline{LEFT: $\{\alpha =2,\; N = 10^2 \}$, RIGHT: $\{\alpha =1.5,\; N= 10^2 \}$.}
 \label{fig:parent100}
\vskip 0.20truecm
 \includegraphics[width=.49\textwidth]{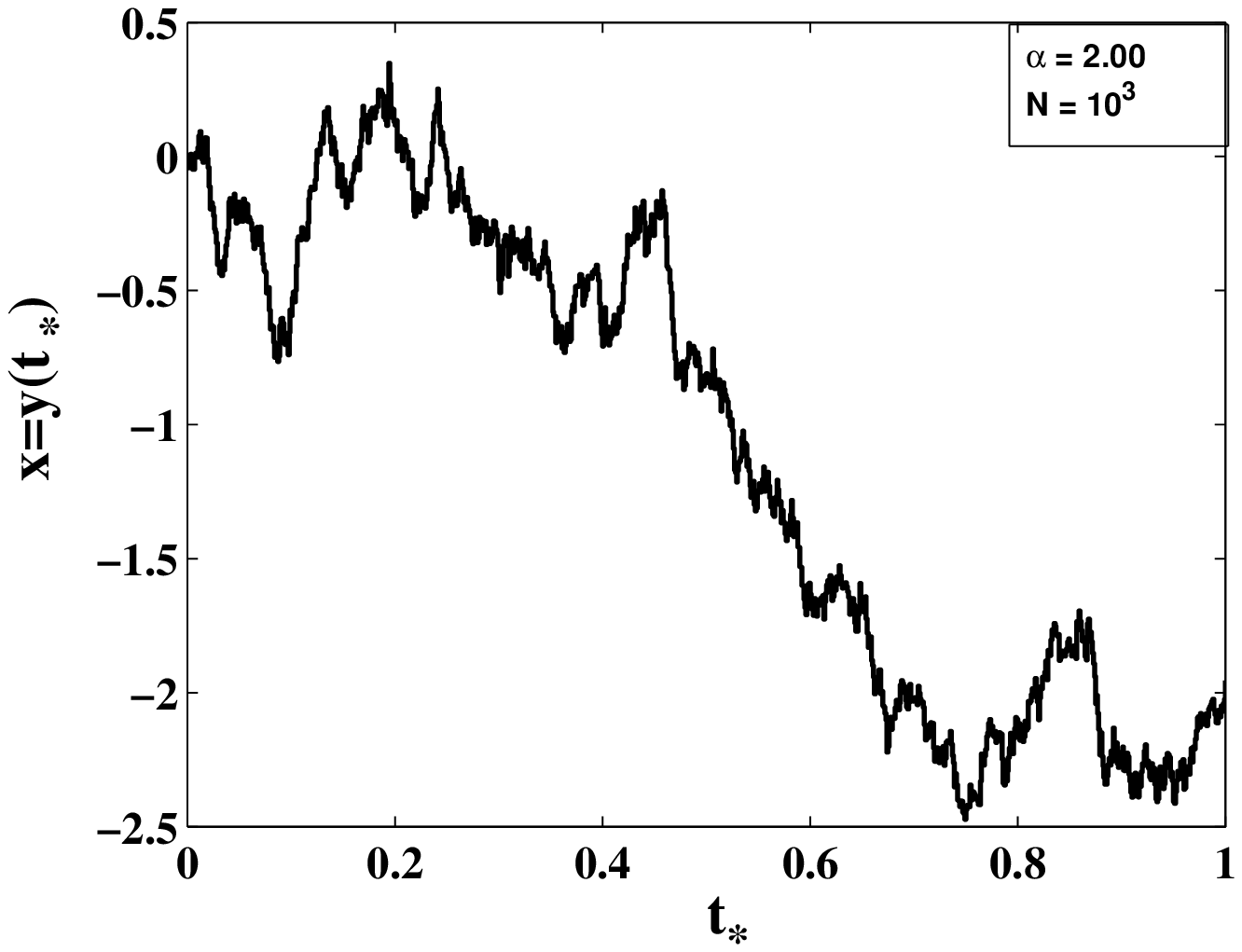}
\includegraphics[width=.49\textwidth]{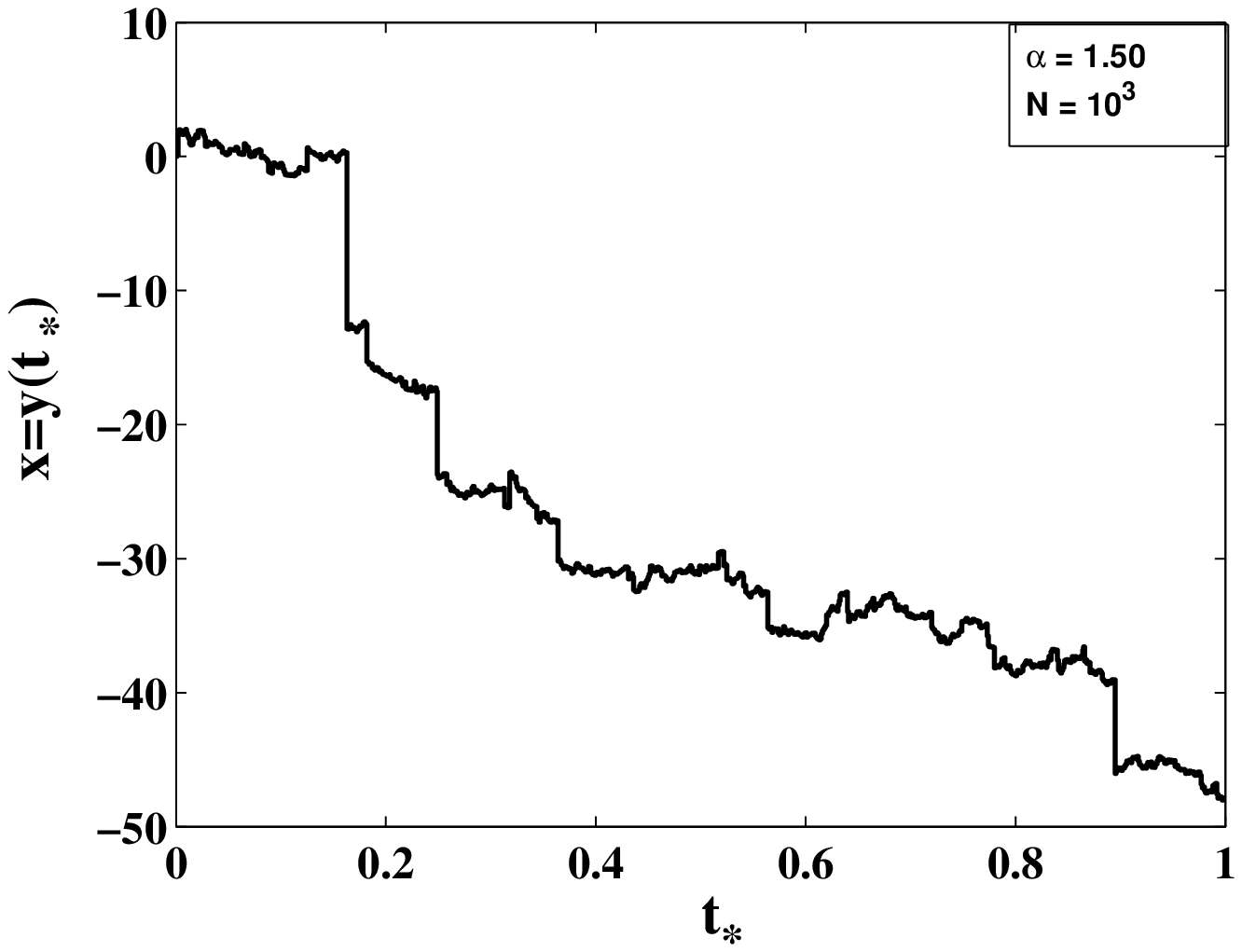}
\vspace{-0.5truecm} 
 \caption{A sample path for the parent process $x=y(t_*)$.}
 \centerline{LEFT: $\{\alpha =2,\; N = 10^3 \}$, RIGHT: $\{\alpha =1.5,\; N= 10^3 \}$.}
\label{fig:parent1000}
\end{figure}
%%%%%%%%%%%%%  REFINEMENT FOR THE 2 SUBORDINATED  PROCESSES
\newpage
\begin{figure}[!h]
 \includegraphics[width=.49\textwidth]{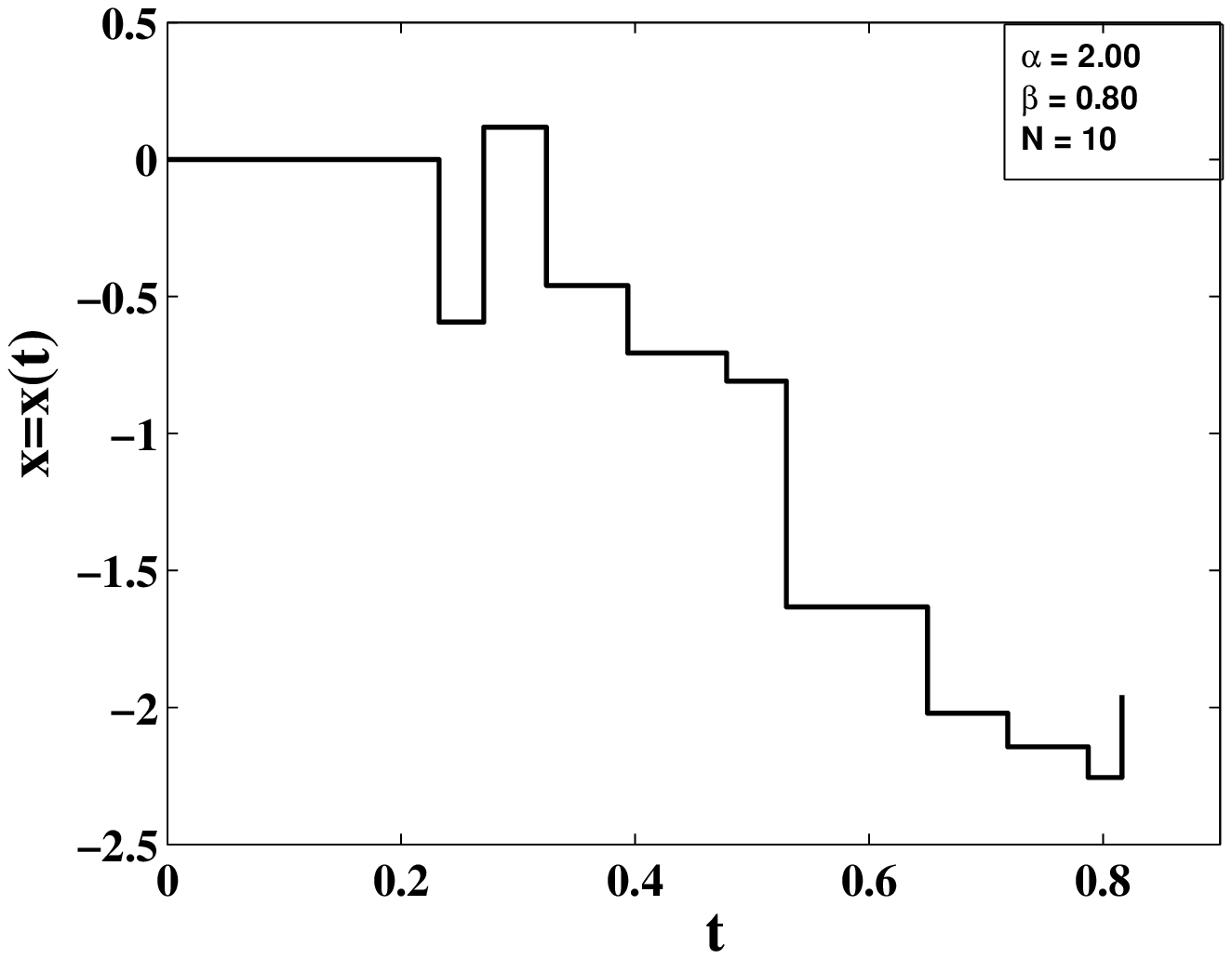}
\includegraphics[width=.49\textwidth]{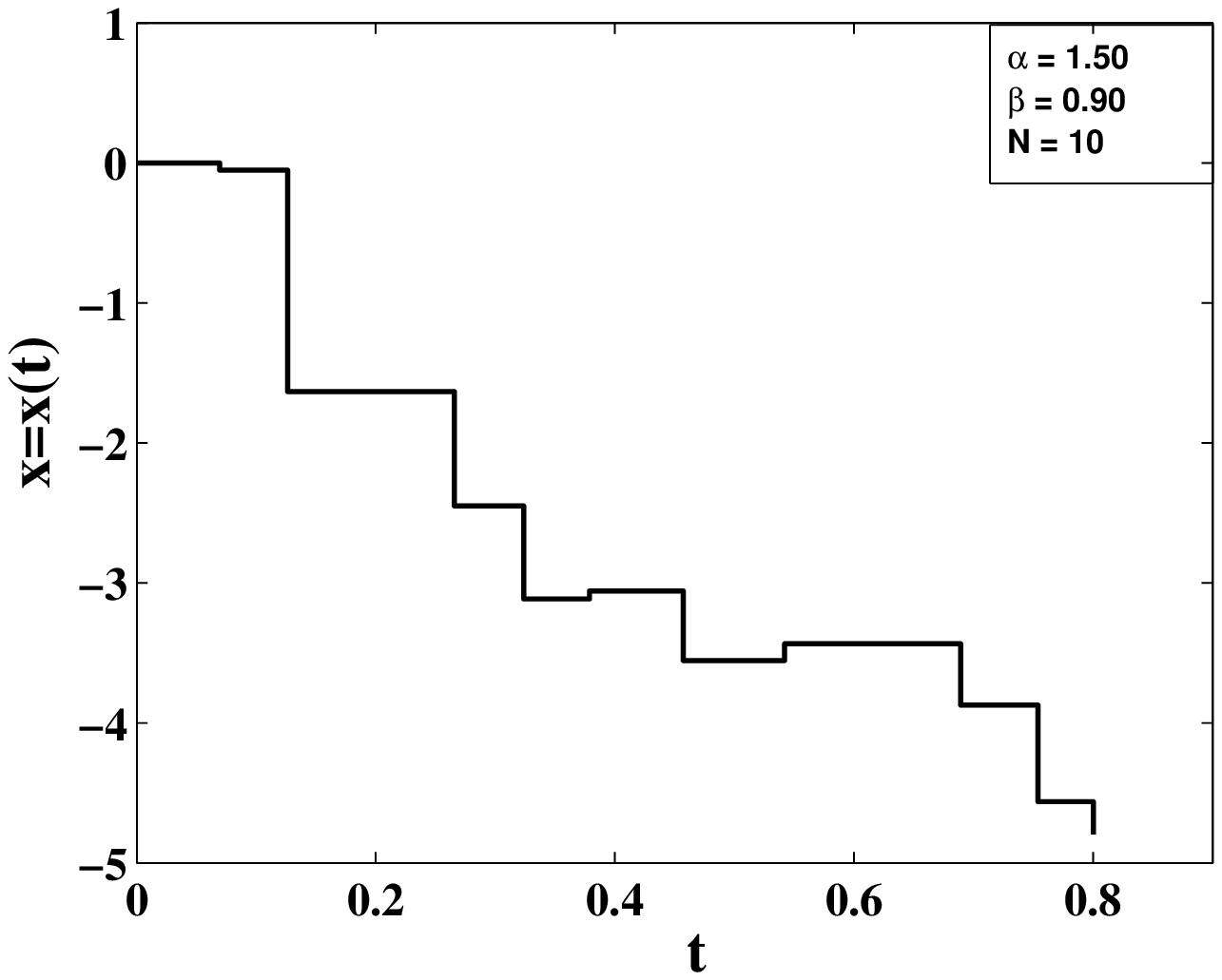}
\vspace{-0.5truecm} 
 \caption{A sample path for the subordinated process $x=x(t)$.}
 \centerline{LEFT: $\{\alpha =2,\; \beta =0.80,\; N=10^1 \}$, RIGHT: $\{\alpha = 1.5,\; \beta =0.90,\; N=10^1 \}$.}
\label{fig:sub10}
\vskip 0.20truecm
 \includegraphics[width=.49\textwidth]{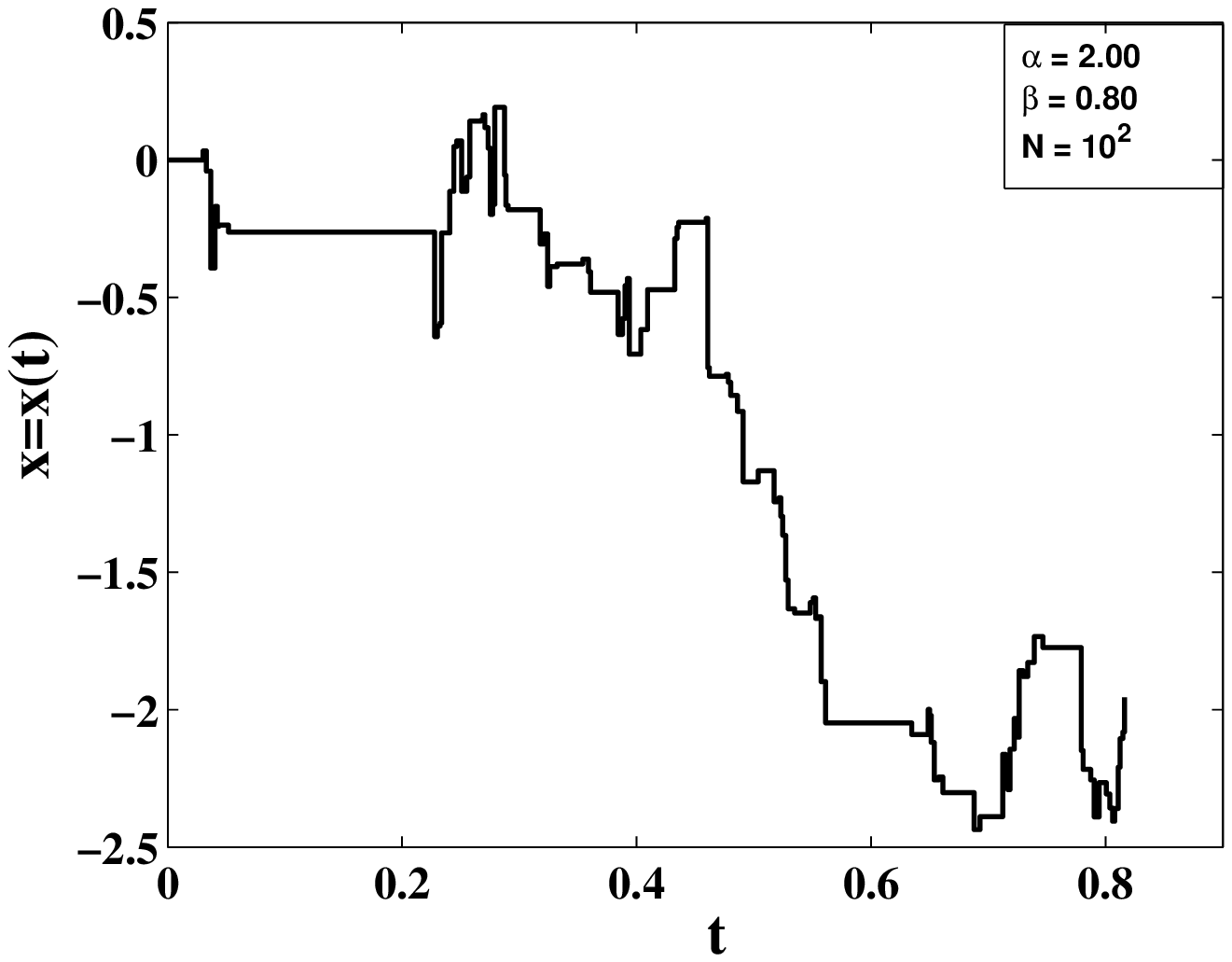}
\includegraphics[width=.49\textwidth]{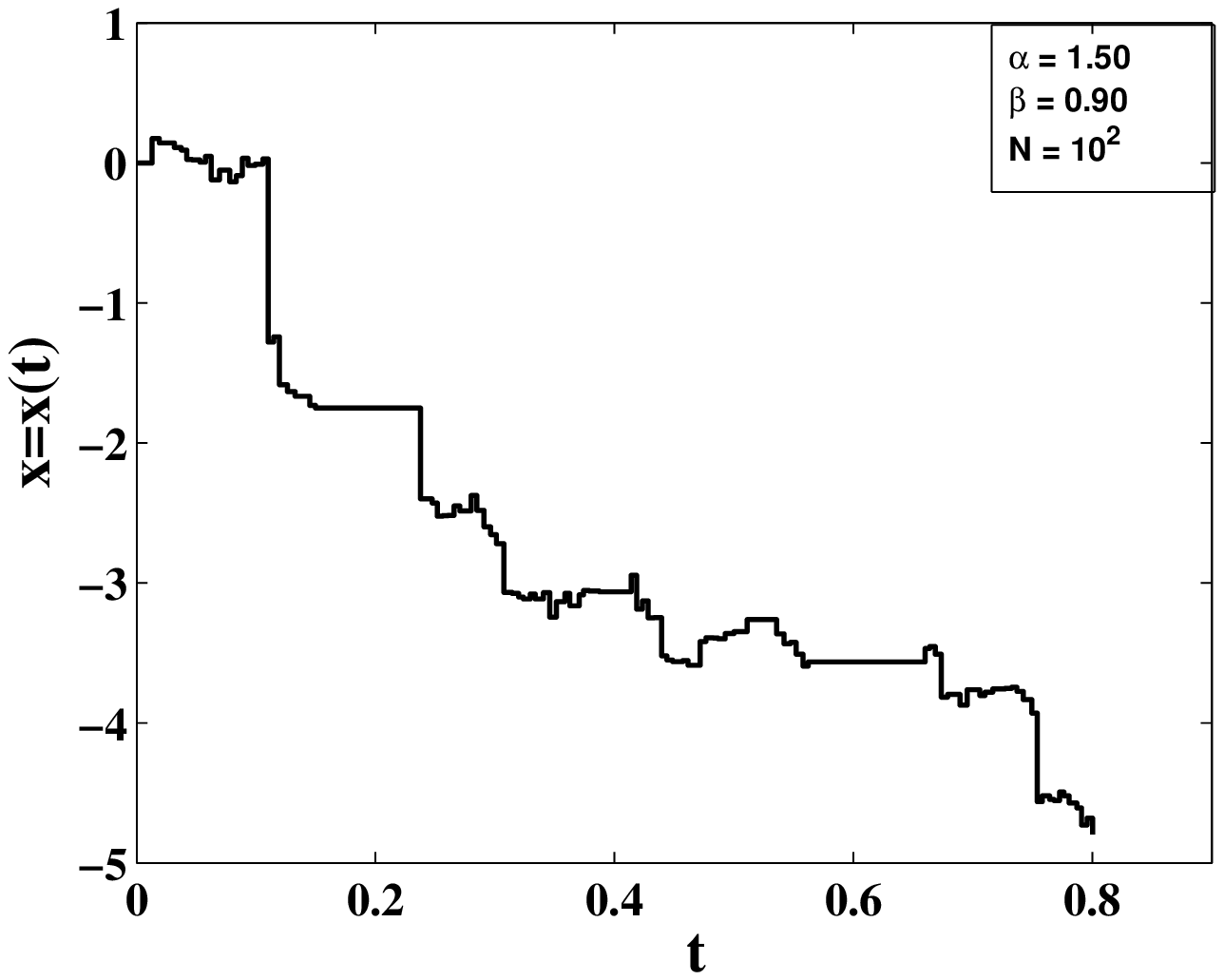}
\vspace{-0.5truecm} 
 \caption{A sample path for the subordinated process $x=x(t)$.}
 \centerline{LEFT: $\{\alpha =2,\; \beta =0.80,\; N=10^2 \}$, RIGHT: $\{\alpha =1.5,\; \beta =0.90,\; N=10^2 \}$.}
\label{fig:sub100}
\vskip 0.30truecm
 \includegraphics[width=.49\textwidth]{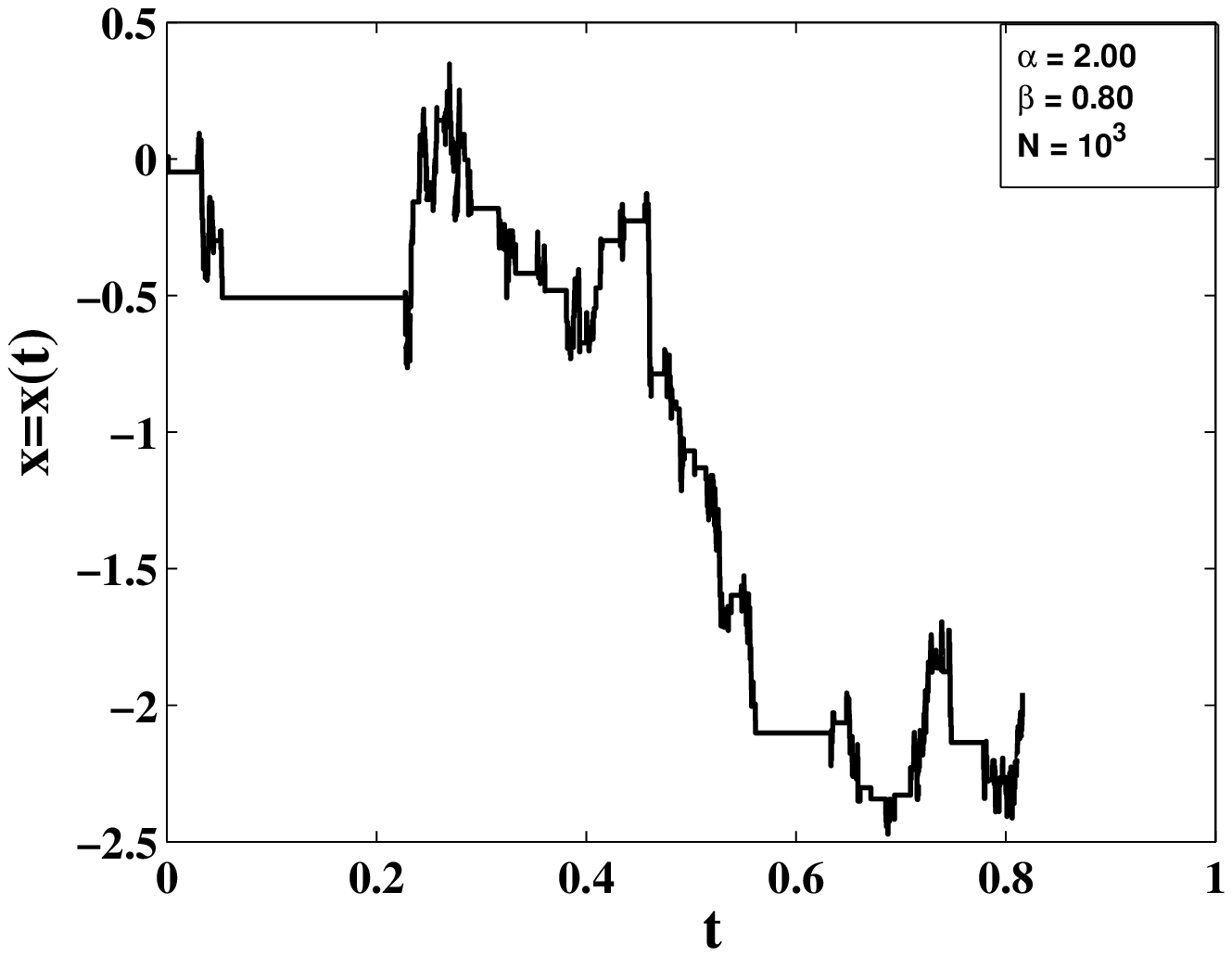}
\includegraphics[width=.49\textwidth]{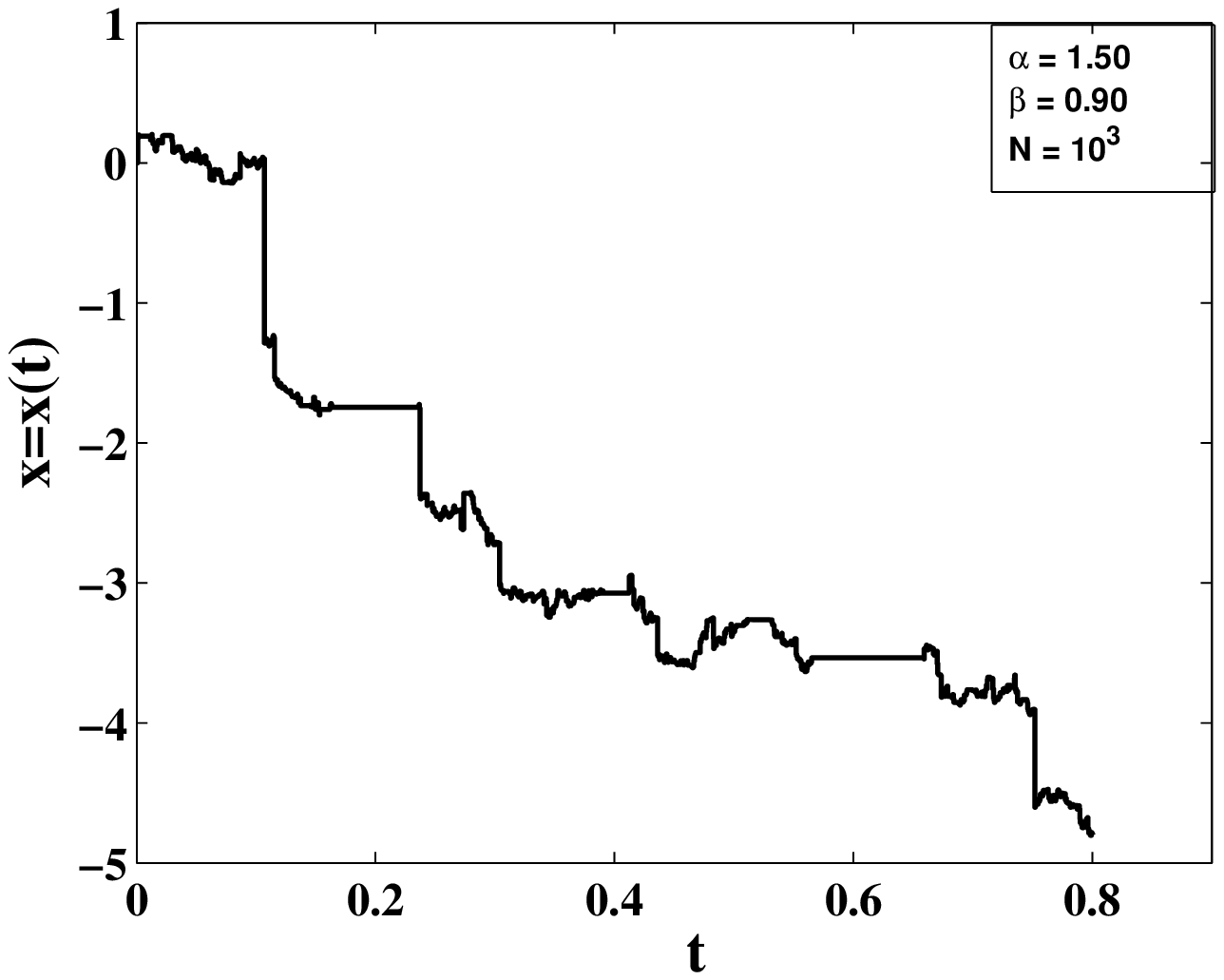}
\vspace{-0.5truecm} 
 \caption{A sample path for the subordinated process $x=x(t)$.}
 \centerline{LEFT: $\{\alpha =2,\; \beta =0.80,\; N=10^3 \}$, RIGHT: $\{\alpha =1.5,\; \beta =0.90,\; N=10^3 \}$.}
\label{fig:sub1000}
\end{figure}
%%%%%%%%%
\vsp
% \section*{Appendix: Finer discretization for path simulations}
We conclude by including additional Figures, showing the effect of taking
smaller step-sizes $\tau_*$, equivalently larger values $N$ of steps following our analysis 
for the $CTRW$\cite{GorMaiViv_CSF07}. 
%%In our CSF paper we have such pictueres with N=10, 100, 1000, making visible finer structures of the
%% processes. We could add the sentence:
%% \vsp
Figs.\ref{fig:leading10},\ref{fig:leading100},\ref{fig:leading1000},
Figs.\ref{fig:parent10},\ref{fig:parent100},\ref{fig:parent1000},
Figs.\ref{fig:sub10},\ref{fig:sub100},\ref{fig:sub1000}
 show the effect of making the operational step-length $\tau_*$
smaller or, equivalently, the number $N$ of operational steps larger
for the sample paths of the leading, parent and subordinated processes, respectively. 
In these
pictures $\tau_*=1/N$, and we have taken $N=10$, $N=100$ and $N=1000$.
Finer details will become visible by choosing in the operational
time $t_*$ the step-length $\tau_*$ smaller and smaller.
%%%%%%%%%%%%%%
%%%%%%%%%%%%%  REFINEMENT FOR THE 2 LEADING  PROCESSES
In the graphs we can clearly see what happens for finer and finer
discretization of the operational time $t_*$, by adopting
$10^1$, $10^2$, $10^3$ of number of steps.
As a matter of fact  there is no visible difference in
the transition  for the successive decades $10^4$, $10^5$, $10^6$
of number of steps as the great majority of spatial jumps and waiting
times even for  very small steps $\tau_*$ of the  operational time. 
 This property also explains the visible  persistence of large jumps and waiting times.
 %% even of very
%% small steps $\tau_*$ of the operational time.
%%%%%%%%%%%%%%%%

\section*{Acknowledgments} 
We  acknowledge the valuable assistance  of  A. Mura and A.  Vivoli 
 in producing the figures. 
We are grateful to the editors of this volume for giving us opportunity 
to present our view of subordination in fractional diffusion  and we have appreciated
their  suggestions.

% \newpage
%%%%%%%%%%%%%%%%%%%%%  BIBLIOGRAPHY
%% GORENFLO-MAINARDI-BIBLIO.tex
%%%%%%%%%
%% \section*{References}

%%%%%%%%

\begin{thebibliography}{99} %
 \bibitem{Barkai_PRE01}
 E. Barkai,
 {Fractional Fokker-Planck equation, solution, and application},
 {\it Phys. Rev. E} {\bf 63}, 046118/1--18 (2001).

\bibitem{Barkai_ChemPhys02}
E. Barkai,
CTRW pathways to the fractional diffusion equation,
{\it Chem. Phys.} {\bf 284}, 13--27 (2002).

\bibitem{Butzer-Westphal_00}
 P. Butzer and U. Westphal,
Introduction to fractional calculus,
 in: H. Hilfer (Editor),
{\it Fractional Calculus, Applications in Physics},
 World Scientific, Singapore (2000), pp. 1--85.

 \bibitem{Diethelm_10}
 K. Diethelm,
 {\it The Analysis of Fractional Differential Equations. 
 An Application Oriented Exposition Using Differential Operators of Caputo Type},  
 Springer, Berlin (2010).
 [Lecture Notes in Mathematics No 2004]

    \bibitem{Erdelyi_HTF}
%%  A. Erd\'elyi (Editor)       Bateman Project,
A. Erd\'elyi, W. Magnus, F. Oberhettinger
 and F.G. Tricomi,
  {\it Higher Transcendental Functions},  Vol. 3, Ch. 18,
 McGraw-Hill, New-York  (1955).
   
 
 
\bibitem{Feller_1971}
 W. Feller, %%  (1971).
{\it An Introduction to Probability Theory and its Applications}, Vol II,
  Wiley, New York (1971).
  %% Second Edition. [First edition (1966)]


 \bibitem{Fogedby_PRE94}
 H.C. Fogedby,
 Langevin equations for continuous time L\'evy flights,
 {\it Phys. Rev. E} {\bf 50}, 1657--1660 (1994).

 \bibitem{Fulger-et-al_PRE08}
D. Fulger, E. Scalas and G. Germano,
 Monte Carlo simulation of uncoupled continuous-time random walks yielding a
 stochastic solution of the space-time fractional diffusion equation,
{\it Phys. Rev. E} {\bf 77}, 021122/1--7 (2008).

\bibitem{Germano-et-al_PRE09}
G. Germano, M. Politi, E. Scalas and R.L. Schilling,
 Stochastic calculus for uncoupled continuous-time random walks,
 {\it Phys. Rev. E} {\bf 79}, 066102/1--12 (2009).

 
\bibitem{Gawronski_1984}
 {W. Gawronski},
  On the bell-shape of stable distributions.
  {\it Ann. Probab.} {\bf 12,}  230--242 (1984).

 

\bibitem{Gorenflo_PALA09}
 R. Gorenflo, Mittag-Leffler waiting time, power laws, rarefaction,
continuous time random walk, diffusion limit, in
 S.S. Pai, N. Sebastian, S.S. Nair, D.P. Joseph and D. Kumar (Editors),
%% CMS Proceedings Pala July 2010, you have been sent this volume.
Proceedings of the National Workshop on  Fractional Calculus and Statistical Distributions,
%% Edited by Shanoja S. Pai, Nicy Sebastian, Seema S. Nair, Dhannya P. Joseph, Dilip Kumar
%% http://www.cmsintl.org/general_announcements/WorkshopOnFractionalCalculus_Proceedings.pdf
%% 25-27 November 2009,
CMS Pala Campus, India (2010), pp.1--22
[E-print: {\tt http://arxiv.org/abs/1004.4413}] %% , 28 pages.




\bibitem{Gorenflo-Abdel-Rehim_VIETNAM04}
R. Gorenflo and E.A. Abdel-Rehim,
 From power laws to fractional diffusion: the direct way,
 {\it Vietnam Journal of Mathematics} {\bf 32} SI, 65--75 (2004).
 [E-print: {\tt http://arXiv.org/abs/0801.0142}]

\bibitem{GorMai_CISM97}  %% [3]
 R. Gorenflo and F. Mainardi,
  Fractional calculus: integral and differential equations of fractional order,
  in: A. Carpinteri and F. Mai\-nardi (Editors),
  {\em Fractals and Fractional Calculus in Continuum Mechanics\/},
  Springer Verlag, Wien  (1997),  pp. 223--276 %% [CISM Lecture Notes Vol. 378]
  [E-print: {\tt http://arxiv.org/abs/0805.3823}]
%%  [Reprinted in NEWS 010101, see {\tt http://www.fracalmo.org}]
%%



\bibitem{GorMai_FCAA98}
R. Gorenflo and	F. Mainardi,
Random walk models for space-fractional diffusion processes.
{\it Fract.  Calc.  Appl. Anal.} {\bf 1}, 167--191 (1998).

%\bibitem{GorMai_INDIA03}    %% [4]
% R. Gorenflo and  F. Mainardi,
% Fractional diffusion processes: probability distributions and continuous time random walk, in:
%    G. Rangarajan and M. Ding (Editors),
%    {\it Processes with Long Range Correlations},
%    Springer-Verlag, Berlin (2003), pp. 148--166. %% [Lecture Notes in Physics, No. 621]
%   [E-print: {\tt http://arxiv.org/abs/0709.3990}]

\bibitem{Gorenflo-Mainardi_HONNEF08}
R. Gorenflo and F. Mainardi, %%  (2008):
Continuous time random walk, Mittag-Leffler waiting  time and fractional diffusion:
mathematical aspects,   in:
R. Klages, G.  Radons,  and I.M. Sokolov  (Editors),
{\it Anomalous  Transport, Foundations and Applications},
 Wiley-VCH Verlag, Weinheim, Germany (2008),  pp. 93--127. %%  ISBN 978-3-527-40722-4.
 [E-print: {\tt http://arxiv.org/abs/0705.0797}]
%% [E-print: {\tt arXiv:cond-mat/07050797}]
 %% Paper presented at the WE-Heraeus-Seminar on Anomalous Transport:
%% Experimental  Results and Theoretical Challenges, Physikzentrum Bad-Honnef (Germany),
%% 12-16 July 2006.

\bibitem{Gorenflo-Mainardi_JCAM09}
R. Gorenflo and F. Mainardi,
Some recent advances in theory and simulation of fractional diffusion processes,
{\it J.  Comp.  Appl. Math.} {\bf 229} No 2, 400--415 (2009).
[E-print: {\tt http://arxiv.org/abs/0801.0146}]
%%  33 pages (Math.PR)R. Gorenflo and F. Mainardi:

\bibitem{Gorenflo-Mainardi_EPJ-ST11}
R. Gorenflo and F. Mainardi,
Subordination pathways to fractional diffusion, 
{\it Eur.  Phys. J. Special Topics} {\bf 193},  119--132 (2011).
[E-print:{\tt  http://arxiv.org/abs/1104.4041}]

\bibitem{GorMaiViv_CSF07} %% [5]
 R. Gorenflo, F. Mainardi and  A. Vivoli,
Continuous time random walk and parametric subordination in fractional diffusion,
   {\it Chaos, Solitons and Fractals} {\bf 34},  87--103 (2007).
   [E-print {\tt http://arxiv.org/abs/cond-mat/0701126}]

\bibitem{Hahn-Umarov_FCAA2011}   %% 
 M. Hahn and S. Umarov,
  Fractional Kolmogorov-Planck type equations and their associated stochastic differential equations,
{\it Fract. Calc. Appl. Anal.} {\bf 14} No. 1, 56--79 (2011).

\bibitem{Hahn-Umarov_JTB2012}
 M. Hahn, K. Kobayashi and S. Umarov,
  SDEs driven by a time-changed Levy process and their associated time-fractional order pseudo-differential equations,
  {\it  J. Theor. Probab.} {\bf 25}, 262--279 (2012).
%%  DOI 10.1007/s10ß59-010-0289-4.

\bibitem{Hilfer_HONNEF08}
 R. Hilfer, Threefold introduction to fractional calculus,
 %% Paper presented at the WE-Heraeus-Seminar on Anomalous Transport:
  %% Experimental Results and Theoretical Challenges, Physikzentrum
%% Bad-Honnef (Germany), 12-16 July 2006.  Published on pages 17-73 in:
%% Rainer Klages, Günter Radons, and Igor M. Sokolov (editors):
in R. Klages, G. Radons and I.M. Sokolov (Editors),
{\it Anomalous  Transport, Foundations and Applications},
Wiley-VCH Verlag, Weinheim, Germany (2008), pp 17--73

 \bibitem{Hilfer_FRACTALS95}
 R. Hilfer,
  Exact solutions for a class of fractal time random walks,
 {\it  Fractals} {\bf 3}, 211--216 (1995).
%%%%%%%%%%%%%%%%%%%%%%
\bibitem{Hilfer-Anton_PRE95}
 {R. Hilfer, L. Anton},
 Fractional master equations and fractal time random walks,
 {\it Phys. Rev. E} {\bf 51}, R848--R851 (1995).


\bibitem{Janicki_LN96}
A. Janicki,
{\it Numerical and Statistical Approximation of
Stochastic Differential Equations with Non-Gaussian Measures},
Monograph No 1, H.  Steinhaus Center for Stochastic Methods
in Science and Technology, Technical University,
Wroclaw,  Poland  (1996). %% pp. 247.
%%%%%%%%%%%%%%%%%%%
\bibitem{Janicki-Weron_94}
 A. Janicki and A. Weron,
 {\it Simulation and Chaotic Behavior of
 $\alpha$--Stable Stochastic Processes},
 Marcel Dekker, New York (1994).


 \bibitem{Kilbas-Srivastava-Trujillo_BOOK06}  %% [6]
 A.A. Kilbas, H.M. Srivastava and  J.J. Trujillo,
{\it Theory and Applications of Fractional Differential Equations},
Elsevier, Amsterdam (2006).
% [North-Holland Series on Mathematics Studies No 204]
%%%%%%%%%
\bibitem{Kleinhans-Friedrich_PRE07}
D. Kleinhans and R. Friedrich,
 Continuous-time random walks: Simulations of continuous trajectories,
 {\it Phys. Rev E} {\bf 76}, 061102/1--6 (2007).


\bibitem{Mainardi_BOOK10}
 F. Mainardi,
{\it Fractional Calculus and Waves in Linear Viscoelasticity},
 Imperial College Press, London (2010)

 \bibitem{Mainardi-Luchko-Pagnini_FCAA01}  %% [7]
 F. Mainardi, Yu. Luchko and  G. Pagnini,
     The fundamental solution of the space-time fractional diffusion equation,
   {\it Fract. Calc.  Appl. Anal.} {\bf 4}, 153--192 (2001)
 [E-print: {\tt http://arxiv.org/abs/cond-mat/0702419}]
%% [Reprinted in NEWS 010401, see {\tt http://www.fracalmo.org}]
% [Paper dedicated to Professor Rudolf Gorenflo for his 70-th birthday]

\bibitem{Mainardi-Mura-Pagnini_IJDE10}
F. Mainardi, A. Mura and G. Pagnini,
The $M$-Wright function in time-fractional diffusion processes: a tutorial survey,
{\it International Journal of Differential Equations} {\bf 2010}, Article ID 104505, 
29 pp (2010).%%  doi:10.1155/2010/104505 .
%%[E-print in {\tt http://www.hindawi.com/a104505.html}
%% Electronic Journal, Hindawi Publishing Corporation, special issue for Fractional
%% Differential Equations http://www.hindawi.com/journals/ijde/contents.html
%% and %% E-print 
[E-print: {\tt http://arxiv.org/abs/1004.2950}]


 \bibitem{Mainardi-Pagnini-Saxena_JCAM05}   %% [8]
 F. Mainardi, G. Pagnini and  R.K. Saxena,
 Fox $H$ functions in fractional diffusion,
 {\it   J. Comp.  Appl. Math.} {\bf 178},  321--331 (2005).
%% Proceedings of the 7-th International Symposium on
%% Orthogonal Polynomials, Special Functions and Applications (OPSFA),
%% Copenhagen (DK) 18-22 August 2003.  URL: www.math.ku.dk/conf/OPSFA2003
%%%%%%%%%%%%%%%%


\bibitem{Mainardi-Tomirotti_GEO97}
F. Mainardi and M. Tomirotti,
  Seismic pulse propagation with constant $Q$ and stable probability
  distributions,
  {\it Annali di Geofisica} {\bf 40},   1311--1328 (1997).
[E-print: {\tt http://arxiv.org/abs/1008.1341}]

\bibitem{Marichev_83}
O.I. Marichev,
{\it Handbook of Integral Transforms of Higher Transcendental Functions,
Theory and Algorithmic Tables}, Ellis Horwood, Chichester (1983).

\bibitem{Mathai-Saxena-Haubold_BOOK-H-2010}
 A.M Mathai, R.K. Saxena and H.J Haubold,
{\it The H-function, Theory and Applications},
Springer Verlag, New York (2010).
%%%%%%%%

\bibitem{M3_PRE02sub_PRE02}  %% [9]
M.M. Meerschaert, D.A. Benson, H.P  Scheffler and B. Baeumer,
Stochastic solutions of space-fractional diffusion equation,
{\it Phys. Rev. E} {\bf  65},  041103/1--4 (2002).
%%A PRESENTATION, GOOD AND SHORT,FOR SUBORDINATION
%%%%%%%%%%%%%%%%%%
\bibitem{M3_PRE02sol_PRE02}  %% [10]
 M.M. Meerschaert, D.A. Benson, H.P  Scheffler and P. Becker-Kern,
  Governing equations and solutions of anomalous random walk limits,
 {\it Phys. Rev. E} {\bf  66},  060102/1--4 (2002).
%%%%%%%%%%%%%%%%%%%%%%
\bibitem{M3_FPP}
M.M. Meerschaert, E. Nane and P. Vellaisamy,
The fractional Poisson process and the inverse stable subordinator,
 {\tt http://arxiv.org/abs/1007.505}   %%%,   22 pp.

\bibitem{Meerschaert-Scheffler_04}
M.M. Meerschaert and H.P  Scheffler,
Limit theorems for  continuous-time random walks with infinite mean waiting times,
{\it J. Appl. Prob.} {\bf  41}, 623--638 (2004).
%%%%%%%%%%%%%%%%%%%%%
\bibitem{Meerschaert-Zhang-Baeumer_CMA10}
M.M. Meerschaert, Y. Zhang and B. Baeumer,
Particle tracking for fractional diffusion with two time scales,
{\it Computers and Mathematics with Applications} {\bf 59}, 1078--1086 (2010).

%%%%%%%%%%%%%%%%%%%%%%
\bibitem{Metzler-Barkai-Klafter_EPL99}
R. Metzler, E. Barkai, and J. Klafter, 
Deriving fractional Fokker-Planck equations from a generalised
master equation, 
{\it Europhys. Lett.} {\bf 46}, 431--436 (1999).


\bibitem{Metzler-Klafter_JPCB00}
R. Metzler and J. Klafter, 
From a generalized Chapman-Kolmogorov equation to the fractional
Klein-Kramers equations, 
{\it J. Phys. Chem. B} {\bf 104}, 3851--3857 (2000).
%%  special issue in honor of Harvey Scher.

\bibitem{Metzler-Klafter-Sokolov_PRE98}
 R. Metzler, J. Klafter, I.M. Sokolov,
 Anomalous transport in external fields:
Continuous time random walks and fractional
diffusion equations extended,
{\em Phys. Rev, E} {\bf 58},  1621-1633 (1998).
%%%%%%%%%%%%%%%%%%%%%%%%%%%%%%%%%%%%%%%%
\bibitem{Metzler-Klafter_PhysRep00}
  R. Metzler and J. Klafter,
 The random walk's guide to anomalous diffusion: a fractional dynamics
 approach, {\it Phys. Reports}  {\bf 339},  1-77 (2000).
%%%%%%%%%%%%%%%%%%%%%


\bibitem{Metzler-Klafter_JPhysics04}  %% [11]
 R. Metzler and J. Klafter,
The restaurant at the end of the random walk: Recent developments
 in the description of anomalous transport by fractional dynamics,
 {\it J. Phys. A. Math. Gen.}  {\bf 37},  R161--R208 (2004).

 \bibitem{Mittag-Leffler_03}
G.M. Mittag-Leffler, 
Sur la nouvelle fonction $E_\alpha (x)$,
{\it C.R. Acad. Sci. Paris} (Ser. II) {\bf 137}, 554--558 (1903).

\bibitem{Mittag-Leffler_04}
G.M. Mittag-Leffler, 
Sopra la  funzione $E_\alpha (x)$,
 {\it Rendiconti R. Accademia Lincei} (Ser. V) {\bf 13}, 3--5 (1904).
%%  [in Italian]

 \bibitem{Mittag-Leffler_05}
G.M. Mittag-Leffler, 
Sur la repr\'esentation analytique d'une branche uniforme
  d'une fonction monog\`ene,
   {\it Acta Math.} {\bf 29}, 101--181 (1905).
 
 
 \bibitem{Podlubny_99}   %% [12]
 I. Podlubny,
  {\it Fractional Differential Equations},
  Academic Press, San Diego (1999).
%%%%%%%%%%%%%
\bibitem{Saichev_PhysA05}
 A. Piryatinska, A.I. Saichev and W.A. Woyczynski,
 Models of anomalous diffusion: the subdiffusive case,
  {\it Physica A}  {\bf 349},  375--420 (2005).
%%%%%%%%%%%%%%%%%%%

%% \bibitem{Rubin_BOOK96}
%%  B. Rubin,
%% {\it Fractional Integrals and Potentials},
%% Addison-Wesley \& Longman, Harlow (1996).
%% [Pitman Monographs and  Surveys in Pure and Applied Mathematics No 82]

\bibitem{Saichev-Zaslavsky_97}
{A. Saichev and G. Zaslavsky},
Fractional kinetic equations: solutions and applications.
{\it Chaos} {\bf 7} (1997), 753-764.


\bibitem{SKM_93}  %5 [13]
 S.G.  Samko, A.A. Kilbas and O.I. Marichev,
{\it Fractional Integrals and Derivatives: Theory  and  Applications},
Gordon and Breach, New York (1993).
% Translation from the Russian edition.
%% Nauka i Tekhnika, Minsk (1987).
%%%%%
\bibitem{Samo-Taqqu_94}
  G. Samorodnitsky and M.S. Taqqu,
  {\em Stable non-Gaussian Random  Processes\/},
  Chapman \& Hall, New York (1994).


\bibitem{Sato_99}
  K-I. Sato,
  {\it L\'evy Processes and Infinitely Divisible Distributions},
  Cambridge University Press, Cambridge (1999).
 %%%%%%%%%%%%%%%%%%%%%%%%%%%
%%%%%%%%%%%%%%%%%%
% \bibitem{SGM_00}
% E. Scalas, R. Gorenflo and F. Mainardi,
% Fractional calculus and continuous-time finance,
% Physica A {\bf 284},  376-384 (2000).
%%%%%%%%%%%%%%%%%%%%%%%%%%%%%%%%%%%%%
\bibitem{Scalas_PRE04} %% [14]
 E. Scalas, R. Gorenflo and F. Mainardi,
    Uncoupled continuous-time random walks:
Solution and limiting behavior of the master equation,
{\it Phys. Rev. E} {\bf 69},  011107/1--8 (2004).
%% (30 January 2004)
%%%%%%


\bibitem{Schneider_LNP86}
 W.R. Schneider,
 Stable distributions: Fox function representation and generalization,
 in S. Albeverio, G. Casati and D. Merlini (Editors),
 {\it Stochastic Processes in Classical and Quantum Systems},
 Springer Verlag, Berlin-Heidelberg (1986), pp.  497--511.
 [Lecture Notes in Physics, Vol. 262]

 \bibitem{Takayasu_FRACTALS}
H. Takayasu,
  {\it Fractals in the Physical Sciences},
 Manchester University Press,
  Manchester and New York (1990). %%  pp. x + 170. [X]
 %% Now published by John Wiley \& Sons, New York.
 

\bibitem{Tomovski-Hilfer-Srivastava_ITSF10}
Z. Tomovski, R. Hilfer and  H.M. Srivastava,
Fractional and operational calculus with generalized fractional
derivative operators and Mittag-Leffler type functions,
{\it Integral Transforms Spec. Funct.} {\bf 21} No 11, 797--814 (2010).


\bibitem{Uchaikin-Zolotarev_99}
{V.V. Uchaikin and V.M. Zolotarev},
 {\it Chance and Stability. Stable Distributions and their Applications},
 VSP, Utrecht (1999).

\bibitem{Weiss_BOOK94} %% [15]
 G.H. Weiss,
 {\it Aspects and Applications of Random Walks},
 North-Holland, Amsterdam (1994).
%%%%%%%%
% \bibitem{West_BOOK03}
% B.J. West, M. Bologna and P. Grigolini.
% {\it Physics of Fractal Operators}, Springer Verlag, New York (2003).
%[Institute for Nonlinear Science]
%%%%%



\bibitem{Wiman_05a}
A. Wiman,
  \"Uber den Fundamentalsatz der Theorie der Funktionen $E_\alpha (x)$,
 {\it Acta Math.} {\bf 29}, 191--201 (1905).




\bibitem{Wright_33}
E.M. Wright,
{On the coefficients of power series having exponential
singularities}, {\em Journal London Math. Soc.\/} {\bf 8}, 71--79 (1933).

\bibitem{Wright_35a}
E.M. Wright, 
The asymptotic expansion of the generalized Bessel function, 
{\em Proc. London Math. Soc. (Ser. II)\/} {\bf 38},  257--270 (1935).

\bibitem{Wright_35b}
{E. M. Wright},
The asymptotic expansion of the generalized hypergeometric function.
{\it Journal London Math. Soc.} {\bf 10},    287--293 (1935).

 \bibitem{Wright_40}
 {E.M. Wright},
 The generalized Bessel function of order greater than one.
 {\em Quart. J. Math., Oxford ser.\/} {\bf 11},  36--48 (1940).



\bibitem{Zhang-Meerschaert-Baeumer_PRE08}
 Y. Zhang, M.M. Meerschaert and B. Baeumer,
Particle tracking for time-fractional diffusion,
{\it Phys. Rev. E.} {\bf 78}, 036705/1--7 (2008).
%%%%%%%%%%

\bibitem{Zolotarev_86}
 {V.M. Zolotarev},    %%  (1982-1986):
 {\it One-dimensional Stable Distributions},
  Amer. Math. Soc.  Providence, R.I. (1986).
  %%  [Translation from the Russian edition (1982)]

\end{thebibliography}
\end{document}